\documentclass[a4paper]{article}

\usepackage[utf8]{inputenc}
\usepackage[T1]{fontenc}
\usepackage{textcomp}
\usepackage[british]{babel}

\usepackage{amsfonts,amssymb,amsthm,amsmath}

\usepackage{geometry}
\usepackage{xcolor}

\usepackage{tikz}
\usepackage{caption}
\usetikzlibrary{
  decorations.pathmorphing,%
  decorations.markings,%
  decorations.pathreplacing,%
  decorations.text,%
  calc,%
  patterns,%
  shapes.geometric,%
  arrows,
  decorations.shapes,
  positioning,
  plotmarks,
  shadings,
  math,
  intersections
}

\usepackage{graphicx}

\definecolor{bleu_sombre}{rgb}{0,0,0.4}
\definecolor{vert_sombre}{rgb}{0,0.3,0}
\definecolor{rouge_sombre}{rgb}{0.3,0,0}
\definecolor{mygray}{rgb}{0.925,0.925,0.925}
\definecolor{mymauve}{rgb}{0.58,0,0.82}
\definecolor{mygreen}{rgb}{0,0.6,0}

\definecolor{myblue}{rgb}{0,0,0.7} 
\usepackage[hyphens]{url}

\usepackage[plainpages=false,colorlinks,linkcolor=myblue,citecolor=rouge_sombre,urlcolor=vert_sombre,breaklinks,bookmarks=true,bookmarksnumbered=true,pdfpagelabels]{hyperref}

\usepackage[capitalize,noabbrev]{cleveref}
\crefname{hypothesis}{Hypothesis}{Hypotheses}

\newtheorem{theorem}{Theorem}
\newtheorem{lemma}[theorem]{Lemma}
\newtheorem{definition}[theorem]{Definition}
\newtheorem{remark}[theorem]{Remark}

\newtheorem{corollary}[theorem]{Corollary}

\numberwithin{theorem}{section}
\numberwithin{equation}{section}
\numberwithin{table}{section}
\numberwithin{figure}{section}

\usepackage{bbm}

\usepackage{siunitx}

\usepackage{empheq}

\usepackage{booktabs}
\usepackage{diagbox}
\usepackage{colortbl}

\usepackage{csquotes}

\usepackage{caption}
\usepackage{subcaption}

\usepackage{placeins}

\newcommand{\nc}{\newcommand}

\makeatletter
\newcommand\RedeclareMathOperator{%
  \@ifstar{\def\rmo@s{m}\rmo@redeclare}{\def\rmo@s{o}\rmo@redeclare}%
}
\newcommand\rmo@redeclare[2]{%
  \begingroup \escapechar\m@ne\xdef\@gtempa{{\string#1}}\endgroup
  \expandafter\@ifundefined\@gtempa
     {\@latex@error{\noexpand#1undefined}\@ehc}%
     \relax
  \expandafter\rmo@declmathop\rmo@s{#1}{#2}}
\newcommand\rmo@declmathop[3]{%
  \DeclareRobustCommand{#2}{\qopname\newmcodes@#1{#3}}%
}
\@onlypreamble\RedeclareMathOperator
\makeatother

\nc{\bfx}{\mathbf{x}} 
\nc{\bfy}{\mathbf{y}} 
\nc{\bfz}{\mathbf{z}} 
\nc{\bfu}{\mathbf{u}} 
\nc{\bfv}{\mathbf{v}} 
\nc{\bfw}{\mathbf{w}} 
\nc{\bft}{\mathbf{t}} 
\nc{\bfb}{\mathbf{b}} 
\nc{\bfn}{\mathbf{n}} 
\nc{\bfr}{\mathbf{r}} 
\nc{\bfc}{\mathbf{c}} 

\nc{\bfA}{\mathbf{A}} 
\nc{\bfB}{\mathbf{B}} 
\nc{\bfC}{\mathbf{C}} 
\nc{\bfR}{\mathbf{R}} 
\nc{\bfD}{\mathbf{D}} 
\nc{\bfI}{\mathbf{I}} 
\nc{\bfM}{\mathbf{M}} 
\nc{\bfK}{\mathbf{K}} 
\nc{\bfP}{\mathbf{P}} 
\nc{\bfU}{\mathbf{U}} 
\nc{\bfZ}{\mathbf{Z}} 
\nc{\bfV}{\mathbf{V}} 
\nc{\bfE}{\mathbf{E}} %
\nc{\bfH}{\mathbf{H}} %
\nc{\bfX}{\mathbf{X}} %
\nc{\bfG}{\mathbf{G}} %
\nc{\bfId}{\mathbf{I}_{\mathrm{d}}} 

\nc{\bbN}{\mathbb{N}} 
\nc{\bbR}{\mathbb{R}} 
\nc{\bbP}{\mathbb{P}} 
\nc{\bbC}{\mathbb{C}} 
\nc{\wH}{\widetilde{H}}
\nc{\calS}{\mathcal{S}} 
\nc{\calD}{\mathcal{D}} 
\nc{\calT}{\mathcal{T}} 
\nc{\calR}{\mathcal{R}} 
\nc{\calP}{\mathcal{P}} 
\nc{\calV}{\mathcal{V}} 
\nc{\calW}{\mathcal{W}} 
\nc{\calJ}{\mathcal{J}} 
\nc{\calE}{\mathcal{E}} 
\nc{\calA}{\mathcal{A}} 
\nc{\calU}{\mathcal{U}} 
\nc{\calK}{\mathcal{K}} 
\nc{\calL}{\mathcal{L}} %
\nc{\calQ}{\mathcal{Q}} %

\nc{\rouge}{\color{red}}
\nc{\bleu}{\color{blue}}
\nc{\cyan}{\color{cyan}}
\nc{\noir}{\color{black}\rm}

\nc\dif{\mathop{}\!\mathrm{d}}   
\DeclareMathOperator{\supp}{supp}   
\DeclareMathOperator{\Span}{span}   
\DeclareMathOperator{\range}{range} 

\DeclareMathOperator{\cond}{cond}
\newcommand{\vertiii}[1]{{\left\vert\kern-0.25ex\left\vert\kern-0.25ex\left\vert #1 
    \right\vert\kern-0.25ex\right\vert\kern-0.25ex\right\vert}} 
\RedeclareMathOperator{\Re}{Re} 
\RedeclareMathOperator{\Im}{Im} 

\nc{\CC}{{C\nolinebreak[4]\hspace{-.05em}\raisebox{.4ex}{\tiny\bf ++}}\ }
\nc{\projout}{\bfP_{{\rm out}}}
\nc{\projcl}{\bfP_{{\rm cl}}}
\nc{\projJ}{\bfP_{\Gamma_j}}

\usepackage[show]{ed}

\newcommand{\cH}{{\cal H}}

\newcommand{\cL}{{\cal L}}

\newcommand{\cS}{{\cal S}}

\newcommand{\cD}{{\cal D}}

\newcommand{\cN}{{\cal N}}

\newcommand{\bx}{x}
\newcommand{\by}{y}
\newcommand{\ba}{\hat{a}}
\newcommand{\br}{\mathbf{r}}

\newcommand{\bv}{\mathbf{v}}
\newcommand{\bu}{\mathbf{u}}

\newcommand{\bff}{\mathbf{f}}

\newcommand{\ta}{\widetilde{a}}

\newcommand{\re}{{\rm e}}
\newcommand{\ri}{{\rm i}}
\newcommand{\rd}{{\rm d}}

\newcommand{\beq}{\begin{equation}}
\newcommand{\eeq}{\end{equation}}
\newcommand{\beqs}{\begin{equation*}}
\newcommand{\eeqs}{\end{equation*}}
\newcommand{\bit}{\begin{itemize}}
\newcommand{\eit}{\end{itemize}}
\newcommand{\ben}{\begin{enumerate}}
\newcommand{\een}{\end{enumerate}}
\newcommand{\bal}{\begin{align}}
\newcommand{\eal}{\end{align}}
\newcommand{\bals}{\begin{align*}}
\newcommand{\eals}{\end{align*}}
\newcommand{\bse}{\begin{subequations}}
\newcommand{\ese}{\end{subequations}}
\newcommand{\bpr}{\begin{proposition}}
\newcommand{\epr}{\end{proposition}}
\newcommand{\bre}{\begin{remark}}
\newcommand{\ere}{\end{remark}}
\newcommand{\bpf}{\begin{proof}}
\newcommand{\epf}{\end{proof}}
\newcommand{\ble}{\begin{lemma}}
\newcommand{\ele}{\end{lemma}}
\newcommand{\bco}{\begin{corollary}}
\newcommand{\eco}{\end{corollary}}
\newcommand{\bex}{\begin{example}}
\newcommand{\eex}{\end{example}}
\newcommand{\bth}{\begin{theorem}}
\newcommand{\enth}{\end{theorem}}

\newcommand{\Rea}{\mathbb{R}}
\newcommand{\Com}{\mathbb{C}}

\newcommand{\Oi}{{\Omega_-}}

\newcommand{\Oe}{{\Omega_+}}

\newcommand{\eps}{\varepsilon}

\newcommand{\pdiff}[2]{\frac{\partial #1}{\partial #2}}

\newcommand{\LtG}{{L^2(\bound)}}

\newcommand{\LtGt}{{\LtG\rightarrow \LtG}}

\newcommand{\HoG}{H^1(\Gamma)}

\newcommand{\tendi}{\rightarrow \infty}
\newcommand{\tendo}{\rightarrow 0}

\def\XXint#1#2#3{{\setbox0=\hbox{$#1{#2#3}{\int}$}
     \vcenter{\hbox{$#2#3$}}\kern-.5\wd0}}

\newcommand*{\N}[1]{\left\|#1\right\|}

\allowdisplaybreaks[4]

\newcommand{\tfa}{\text{ for all }}
\newcommand{\tfor}{\text{ for }}

\newcommand{\ton}{\text{ on }}
\newcommand{\tas}{\text{ as }}
\newcommand{\tand}{\text{ and }}
\newcommand{\tst}{\text{ such that }}

\newcommand{\bound}{\Gamma}

\usepackage{soul}

\definecolor{jwcol}{RGB}{27, 137, 18}  

\definecolor{dalcol}{rgb}{0.8,0,0}

\definecolor{jeffColor}{RGB}{102, 0, 204}

\definecolor{escol}{rgb}{0,0,0.8}
\definecolor{estcol}{rgb}{0,0.5,0}
\definecolor{esnewcol}{rgb}{0,0.5,0}

\newcommand{\QMC}{\epsilon}
\newcommand{\genmatrix}{\bfB}
\newcommand{\matrixD}{\bfA_k'}
\newcommand{\matrixDj}{\bfA_{k_j}'}
\newcommand{\bfbeta}{\boldsymbol{\beta}}
\newcommand{\smallpara}{\nu}

\newcommand{\Capproxone}{{C_{{\rm approx},1}}}
\newcommand{\Capproxtwo}{{C_{{\rm approx},2}}}
\newcommand{\CWeyl}{{C_{{\rm Weyl}}}}
\newcommand{\Ccond}{{C_{{\rm cond}}}}

\title{Applying GMRES to the Helmholtz equation with strong trapping: how does the number of iterations depend on the frequency?}

\author{
    P.~Marchand\thanks{Department of Mathematical Sciences, University of Bath, Bath, BA2 7AY, UK, \tt pfcm20@bath.ac.uk}
    \and
        J.~Galkowski\thanks{Department of Mathematics, University College London, 25 Gordon Street, London, WC1H 0AY, \tt j.galkowski@ucl.ac.uk}
    \and
    A.~Spence\thanks{Department of Mathematical Sciences, University of Bath, Bath, BA2 7AY, UK, \tt A.Spence@bath.ac.uk}
    \and
    E.~A.~Spence\thanks{Department of Mathematical Sciences, University of Bath, Bath, BA2 7AY, UK, \tt E.A.Spence@bath.ac.uk}}
\date{
    \today
}

\begin{document}

\maketitle

\begin{abstract}
We consider GMRES applied to discretisations of the high-frequency Helmholtz equation with strong trapping; recall that in this situation the problem is exponentially ill-conditioned through an increasing sequence of frequencies. 
Our main focus is on boundary-integral-equation formulations of the exterior Dirichlet and Neumann obstacle problems in 2- and 3-d.
Under certain assumptions about the distribution of the eigenvalues of the integral operators, we prove upper bounds on how the number of GMRES iterations grows with the frequency; 
 we then investigate numerically the sharpness  (in terms of dependence on frequency) of \emph{both} our bounds
\emph{and}
various quantities entering our bounds.
This paper is therefore the first comprehensive study of the frequency-dependence of the number of GMRES iterations for Helmholtz boundary-integral equations under trapping.
\end{abstract}

\section{Introduction}

\subsection{Statement of the problem}
\label{sec:statement}

We consider solving the Helmholtz obstacle-scattering problem, where the obstacle traps geometric-optic rays, by the boundary-element method  arising from the Galerkin method applied to boundary-integral-equation formulations of the PDE problem, and then solving the resulting linear systems with the generalised minimum residual method (GMRES). We now give details of each of these aspects.

\subsubsection{The scattering problem.}

Let $\Omega_-\subset \Rea^d$, $d=2,3$, be a bounded Lipschitz open set such that its open complement $\Omega_+:=\Rea^d \setminus \overline{\Omega_-}$ is connected; let $\Gamma:= \partial \Omega_-$, and let $n$ be the outward-pointing unit normal vector to $\Omega_-$. We consider the exterior Dirichlet and Neumann scattering problems. For simplicity, we consider the case when the boundary data comes from an incoming plane wave $u^I(\bx):= \exp(\ri k \bx\cdot \ba)$ for $\ba\in\Rea^d$ with $\|\ba\|_2=1$; i.e. we consider 
\emph{either} the sound-soft \emph{or} the sound-hard plane-wave scattering problem defined by: given $k>0$ and the incident field $u^I$, find the total field $u$ satisfying 
\begin{align}\label{eq:Helmholtz}
    \left\{
        \begin{aligned}
            &\Delta u + k^2 u =0 \hspace{2.5cm}\text{ in }\Omega^+,\\
            &\text{either } u \text{ or } \partial_n u = 0 \hspace{1.65cm}\text{ on }\Gamma, \quad\tand\\
            &\dfrac{\partial u^S }{\partial r} -\ri ku^S = o \left(\frac{1}{r^{(d-1)/2}}\right)  \quad\text{as }r:=|x|\rightarrow \infty, \text{ uniformly in $x/r$},
        \end{aligned}
    \right.
    \end{align}
where \(u^S := u - u^I\) is the scattered field. We are particularly interested in the case when the frequency $k$ is large. 

\subsubsection{Trapping and quasimodes}

We consider domains $\Oi$ such that there exist stable trapped geometric-optic rays in the exterior $\Oe$.
In this situation, the solution operator for the problem \eqref{eq:Helmholtz} can grow exponentially through an increasing sequence of frequencies. This phenomenon can be expressed via the notion of \emph{quasimodes}.

\begin{definition}[Quasimodes]\label{def:quasimodes}
A family of Dirichlet quasimodes of quality $\QMC(k)$
is a sequence $\{(u_j,k_j)\}_{j=1}^\infty\subset H^1_{\rm loc}(\Oe)\times \mathbb{R}$ 
with $u_j=0$ on $\Gamma$ 
such that the frequencies $k_j\tendi$ as $j \tendi$ and there exists a compact subset $\mathcal{K}\subset \Oe$ such that, for all $j$, $\supp\, u_j \subset \mathcal{K}$,
\beqs
\N{(\Delta +k_j^2) u_j}_{L^2(\Oe)} \leq \QMC(k_j), \quad\tand\quad\N{u_j}_{L^2(\Oe)}=1.
\eeqs
The definition of Neumann quasimodes is analogous, with $u_j=0$ on $\Gamma$ replaced by $\partial_n u =0$ on $\Gamma$.
\end{definition}

By the results of \cite[Theorem 2]{Bu:98} (see also \cite{Vo:00}), if a family of Dirichlet or Neumann quasimodes exists, the quality $\QMC(k)$ can be at most exponentially-small in $k$.

For simplicity, in our numerical experiments we focus on the case when $\Oi$ is \emph{either} one of the two ``horseshoe-shaped'' 2-d domains shown in Figure \ref{fig:geometries} (and defined precisely below) \emph{or} certain 3-d analogues (defined in \S\ref{sec:3-dexpts} below); in these cases 
there exist quasimodes with exponentially-small quality, leading to exponential growth of the solution operator -- see Theorem \ref{thm:ellipse} below.

We emphasise however that 
there exist quasimodes with superalgebraically small quality for a much larger class of obstacles (see \cite[Theorem 1]{CaPo:02}, \cite[Theorem 1]{St:00} and the discussion in \S\ref{sec:Weyl_new}) 
and our bound on the $k$-dependence of the number of GMRES iterations (in Theorem \ref{thm:main2})
hold in these more-general situations. 
The existence of quasimodes is linked to the existence of \emph{resonances} (poles of the meromorphic continuation of the solution operator of \eqref{eq:Helmholtz} from $\Im k\geq 0$ to $\Im k<0$);  the relationship between trapping, quasimodes, and resonances is a classic topic in scattering theory; see \cite{StVo:95, StVo:96, TaZw:98, St:99, St:00} and \cite[Chapter 7]{DyZw:19}.

\subsubsection{A particular class of $\Oi$ for which quasimodes exist.}\label{sec:statement:ellipse}

The following theorem is proved by combining \cite[Equation A.16]{BeChGrLaLi:11} and \cite[Theorem 3.1]{NgGr:13} (see \S\ref{sec:Mathieu}).

\begin{theorem}[Quasimodes when $\Oe$ contains part of an ellipse]\label{thm:ellipse}
Let $d=2$. Given $a_1>a_2>0$, let 
\beq\label{eq:ellipse}
E:= \left\{(x_1,x_2) \, : \, \left(\frac{x_1}{a_1}\right)^2+\left(\frac{x_2}{a_2}\right)^2<1\right\}.
\eeq
Assume that $\Gamma$ coincides with the boundary of $E$ in the neighborhoods of the points \((0,\pm a_2)\), 
and that $\overline{\Omega_+}$ contains the convex hull of the union of these neighbourhoods.

Then there exist families of Dirichlet and Neumann quasimodes with 
\beq\label{eq:quality}
\QMC(k)=C_1 \exp( - C_2 k) \quad\tfa k>0.
\eeq
where $C_1, C_2>0$ are both independent of $k$.
\end{theorem}

For $\Oi$ satisfying the assumptions of Theorem \ref{thm:ellipse}, we can compute the frequencies $k_j$ in the quasimodes. Indeed, the functions $u_j$ in the quasimode construction in \cite{BeChGrLaLi:11}/\cite{NgGr:13} are based on the family of eigenfunctions of the ellipse localising around the periodic orbit $\{(0,x_2) : |x_2|\leq a_2\}$ (i.e.~the minor axis of the ellipse); 
when the eigenfunctions are sufficiently localised, the eigenfunctions multiplied by a suitable cut-off function form a quasimode, with frequencies $k_j$ equal to the square roots of the respective eigenvalues of the ellipse. By separation of variables, $k_j$ can be expressed as the 
solution of a multiparametric eigenvalue problem involving Mathieu functions; see Appendix \ref{sec:Mathieu}.
We use the method introduced in~\cite{Wi:06} and the associated MATLAB toolbox to solve these eigenvalue problems for $k_j$.
When giving values of these $k_j$s we give all the digits computed in double precision. 
Note that we are not claiming that all these digits are accurate (see \cite{Wi:06} for some discussion on accuracy), but 
some of the quantities we compute below are very sensitive to the precise values of $k$, and so we give the exact values of $k$ used in our computations.

When giving specific values of these $k_j$, we use the notation from \cite[Appendix A]{BeChGrLaLi:11}, recapped in Appendix \ref{sec:Mathieu}, that $k_{m,n}^e$ and $k_{m,n}^o$ are the frequencies associated with the eigenfunctions of the ellipse that are even/odd, respectively, in the angular variable, with $m$ zeros in the radial direction (other than at the centre or the boundary) and $n$ zeros in the angular variable in the interval $[0,\pi)$. Note that the values of $k_{m,n}^e$ and $k_{m,n}^o$ are different for Dirichlet and Neumann boundary conditions, but we do not indicate this difference in our notation.
The eigenfunctions associated to $k_{m,n}^e$, $k_{m,n}^o$ localise about the minor axis as $m\tendi$ for fixed $n$ (see the proof of Theorem \ref{thm:ellipse} in Appendix \ref{sec:Mathieu}); therefore quasimodes exist for the families of frequencies $\{k^{e/o}_{m,n}\}_{m=1}^\infty$ for fixed $n$.

\subsubsection{Definitions of the ``small cavity'' and ``large cavity'' obstacles $\Oi$.}\label{sec:smalllargegeometries}

Our numerical experiments focus on two specific $\Oi$ satisfying the assumptions of Theorem \ref{thm:ellipse}
with $a_1=1$ and $a_2=1/2$. We define the \emph{small cavity} as the region between the two elliptic arcs 

\begin{figure}[h!]
    \centering
    \includegraphics[width=0.6\textwidth]{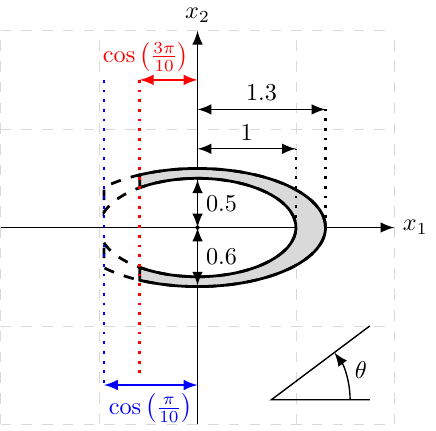}
    \caption{The ``small'' and ``large'' cavities $\Oi$, with the small cavity shaded in grey, and the large cavity equal to the union of the small cavity and the dashed region. We write the incident-plane-wave direction $\widehat{a}= (\cos\theta,\sin\theta)$, with the angle $\theta$ measured in the positive direction from the horizontal, as pictured. }\label{fig:geometries}
\end{figure}

\begin{figure}[h!]
    \centering\includegraphics[width=0.8\textwidth]{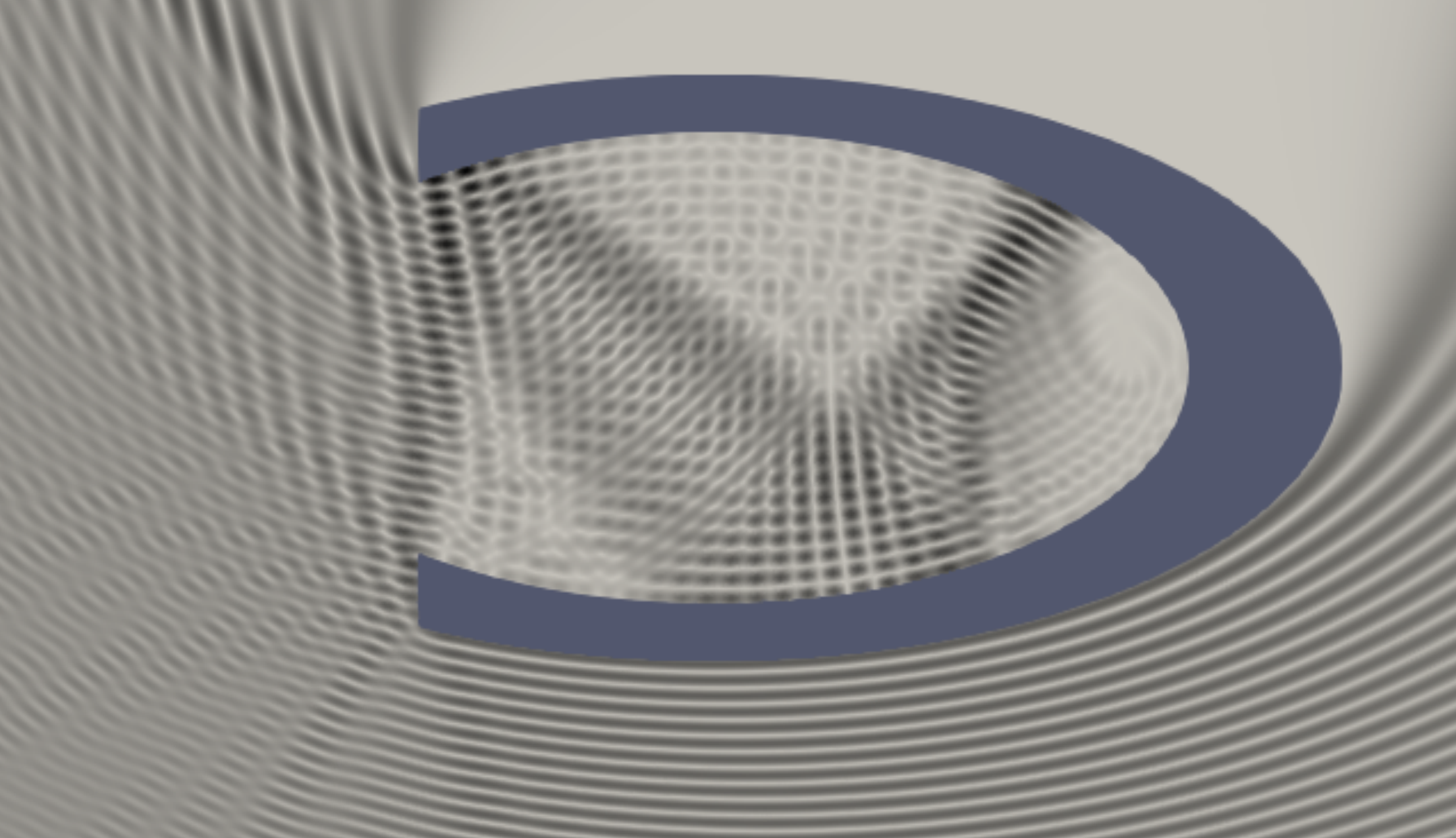}
    \caption{Absolute value of the total field \(u\) defined by \eqref{eq:Helmholtz} with Dirichlet boundary conditions, $\Oi$ the small cavity, $k=100$, \(\hat{a}=(\cos(\theta),\sin(\theta))\), and \(\theta=4\pi/10\). 
    }\label{fig:example}
\end{figure}

\begin{align*}
&(\cos (t), 0.5 \sin (t)), \quad t\in [-\phi_0,\phi_0]
\quad\tand \quad 
(1.3\cos (t), 0.6 \sin (t)),\quad t\in [-\phi_1,\phi_1] \\
&\qquad\qquad\text{ with } \phi_0=7\pi/10\quad\tand \quad \phi_1= \arccos \left(\frac{1}{1.3}\cos (\phi_0)\right);
\end{align*}
this corresponds to the shaded interior of the solid lines in Figure \ref{fig:geometries}. 
We define the \emph{large cavity} as the region between the two arcs now with \(\phi_0=9\pi/10\).
We also consider 3-d analogues of the these cavities, created by rotating them around the $x_1$ axis. 

Figure~\ref{fig:example} plots the absolute value of the total field \(u\) satisfying \eqref{eq:Helmholtz} with Dirichlet boundary conditions with $\Oi$ the small cavity, $k=100$ and \(\hat{a}=(\cos(\theta),\sin(\theta))\) with \(\theta=4\pi/10\); 
this figure was produced by computing the unknown Neumann trace using BEM, and then evaluating the solution given 
in terms of layer potentials by Green's integral representation \eqref{eq:Green}.

\subsubsection{Boundary-integral-equation (BIE) formulations of \eqref{eq:Helmholtz}}\label{sec:BIE}

We are primarily interested in solving \eqref{eq:Helmholtz} by reformulating it as an integral equation on $\Gamma$; recall that this procedure has the advantage of converting a problem posed in an unbounded $d$-dimensional domain (i.e. $\Oe$) to a problem posed on a bounded $(d-1)$-dimensional domain (i.e. $\Gamma$). 
However, the ideas behind our main results are applicable to other methods of solving the Helmholtz equation, and we discuss in \S\ref{sec:FEM} below the standard variational formulation, which is the basis of the finite-element method.

We consider \emph{direct} BIE formulations of \eqref{eq:Helmholtz}, i.e., ones in which the unknown is either the Neumann data (for the Dirichlet problem) or the Dirichlet data (for the Neumann problem); however, our results below also apply to \emph{indirect} formulations (where the unknown has less-immediate physical relevance; see \cite[Page 132]{ChGrLaSp:12}), since there is a close relationship between the integral operators of the direct and indirect formulations; see, e.g., \cite[Remark 2.24, \S2.6]{ChGrLaSp:12}. Once both $u$ and $\partial_n u$ are known on $\Gamma$, the solution in $\Oe$ can be obtained from Green's integral representation (\eqref{eq:Green} below). 
For the Dirichlet problem we find $\partial_n u$ using the standard ``combined-field'' or ``combined-potential'' BIE 
\begin{align}\label{eq:direct_combined_dir}
    A_{k,\eta}' \partial_{n} u = \partial_{n} u^I - \ri \eta u^I 
    \quad\ton \Gamma,
    \quad\text{ where } \quad
    A_{k,\eta}' := \dfrac{1}{2}I+D'_k-\ri \eta S_k,
\end{align}
where $\eta$ is the (arbitrary) ``coupling parameter'' and $S_k$ and $D'_k$ are the single-layer and adjoint-double-layer operators defined by \eqref{eq:SD'} below. 
If $k>0$ and \(\Re(\eta)\neq 0\),  then \(A'_{k,\eta}:L^2(\Gamma)\rightarrow L^2(\Gamma)\) is bounded and invertible (see, e.g., \cite[Theorem 2.27]{ChGrLaSp:12}). There has been much research on the question of how to best choose $\eta$, 
starting from the works \cite{KrSp:83, Kr:85, Am:90} for the case when $\Omega$ is a ball; see the overviews in \cite[Chapter 5]{ChGrLaSp:12}, \cite[\S7]{BaSpWu:16}, \cite[\S6.5]{ChSpGiSm:20}. Roughly speaking, the best choice for large $k$ is $\eta=k$; therefore, in the rest of the paper we take \(\eta=k\), and let \(A'_{k}:=A'_{k,k}\).

For the Neumann problem, the standard ``combined-field'' or ``combined-potential'' BIE is
\begin{align}\label{eq:direct_combined_neu}
    B_{k,\eta} u  = \ri \eta u^I  - \partial_n u^I    \quad\ton \Gamma,
\quad\text{ where }\quad 
    B_{k,\eta} := \ri\eta \left(\dfrac{1}{2}I-D_k\right) +  H_k,
\end{align}
where $D_k$ and $H_k$ are the double-layer and hypersingular operators defined by \eqref{eq:DH} below. In contrast to $A'_k$, $B_{k,\eta}$ is not a bounded operator on $\LtG$ (even when $\Gamma$ is smooth) because of the hypersingular operator $H_k$. If $k>0$, \(\Re(\eta)\neq 0\), and $\Oi$ is Lipschitz,  then $B_{k,\eta}:H^{s+1/2}(\Gamma)\rightarrow H^{s-1/2}(\Gamma)$ is bounded and invertible for $|s|\leq 1/2$ (see, e.g., \cite[Theorem 2.27]{ChGrLaSp:12}). The standard choice of $\eta$ here is also $\eta=k$, and we let \(B_k:=B_{k,k}\).

The fact that $B_k$ is not bounded from $\LtG\rightarrow \LtG$ means that the condition numbers of $h$-version Galerkin discretisation of \eqref{eq:direct_combined_neu} blow up as $h\tendo$ for fixed $k$. 
There has therefore been much research interest in designing alternative Neumann BIE formulations; see, e.g., \cite{StWe:98, AnDa:05, AnDa:07, BrElTu:12, DaDaLa:13}. We use the following BIE, introduced in \cite{BrElTu:12} (which focused specifically on  high-frequency problems)
 based on the idea of Calder\'on preconditioning, 
\begin{align}\label{eq:direct_combined_neu2}
    B_{k,\eta,{\rm reg}} \gamma_+ u  = \ri \eta \gamma^+ u^I  - S_{\ri k} \partial_{n}^+ u^I    \quad\ton \Gamma,
\quad\text{ where } \quad 
    B_{k,\eta,{\rm reg}} := \ri\eta \left(\dfrac{1}{2}I-D_k\right) +  S_{\ri k} H_k.
\end{align}
At least when $\Gamma$ is $C^1$, if $k>0$ and \(\Re(\eta)\neq 0\), then $B_{k,\eta,{\rm reg}}:\LtG\rightarrow \LtG$ is bounded and invertible \cite[Theorem 2.1]{BrElTu:12}. In what follows, we make the same choice for \(\eta\) as in~\cite[Equation 24]{BrElTu:12}, i.e. \(\eta=1/2\), and we let \(B_{k,{\rm reg}}:=B_{k,1/2,{\rm reg}}\).  
We highlight that the idea of combatting the ``bad'' behaviour of the hypersingular operator by composing it with a regularising operator (in this case 
$S_{\ri k}$)  is often called ``operator preconditioning'' (see \cite{Hi:06}).

Because the normality (or not) of an operator is relevant for the analysis of GMRES, we highlight that when $\Oi$ is not a ball, $A_{k}'$ and $B_{k, {\rm reg}}$ are non-normal operators on $\LtG$; this is shown by the plots of the numerical range for $A_k'$ in \cite{BeSp:11} and for $B_{k, {\rm reg}}$ in \cite[\S5]{BoTu:13} (see also \cite{BePhSp:13}).

\subsubsection{The boundary-element method (BEM).}\label{sec:statement:BEM}

We solve the BIEs
\eqref{eq:direct_combined_dir} and \eqref{eq:direct_combined_neu2} with the Galerkin method in $\LtG$.
That is, to solve the BIE \eqref{eq:direct_combined_dir} given a finite-dimensional subspace $V_n\subset \LtG$, we
\beq\label{eq:Galerkin}
\text{ find } v_n \in V_n \tst \big(A_k' v_n, w_n\big)_\LtG = \big(f, w_n\big)_\LtG \quad\tfa w_n \in V_n,
\eeq
where $f$ denotes the right-hand side of the BIE in \eqref{eq:direct_combined_dir}; the Galerkin solution $v_n$ is then an approximation to $\partial_n u$.
We solve the BIE \eqref{eq:direct_combined_neu} via the Galerkin method in $H^{1/2}(\Gamma)$. 
That is, given a finite-dimensional subspace $V_n\subset H^{1/2}(\Gamma)$, we
\beqs
\text{ find } v_n \in V_n \tst \big\langle B_k v_n, w_n\rangle = \big\langle f, w_n\big\rangle \quad\tfa w_n \in V_n,
\eeqs
where $f$ now denotes the right-hand side of the BIE in \eqref{eq:direct_combined_neu}, and $\langle\cdot,\cdot\rangle$ denotes the duality pairing between $H^{-1/2}(\Gamma)$ and $H^{1/2}(\Gamma)$.

Given a basis $\{\phi_j\}_{j=1}^n$ of $V_n$, the Galerkin equations \eqref{eq:Galerkin} are equivalent to the linear system 
\beq\label{eq:1}
\matrixD \bv = \bff
\eeq
where 
\begin{align}\label{eq:Galerkinmatrix}
    (\matrixD)_{i,j}:= \int_{\Gamma} \big(A_k'\phi_j\big) (\bx)\phi_i (\bx) \dif \sigma(\bx) 
    \quad \tand\quad 
    (\bff)_i:=  \int_{\Gamma} f(\bx) \,\phi_i (\bx) \dif \sigma(\bx).
    \end{align}    
    Regarding notation: we put in bold font the $n\times n$ matrices and $n\times 1$ vectors arising from the Galerkin method -- such as $\matrixD$, $\bv$, and $\bff$ in \eqref{eq:1} --  but do not put in bold font the position vectors in $\Rea^d$ -- such as $\bx$ in \eqref{eq:Galerkinmatrix}.
    
We consider the $h$-version of the boundary-element method, and choose $V_n$ to be P1 Lagrange elements (i.e.~continuous piecewise-linear polynomials on the reference elements). 
To maintain accuracy as $k\tendi$, $h$ must be tied to $k$. In applications, one usually chooses $h$ to be proportional to $1/k$, i.e. a fixed number of points per wavelength (see, e.g, \cite{Ma:02}), and we do the same for the numerical experiments in this paper.
At least when $\Oi$ is nontrapping, empirically one sees uniform accuracy as $k\tendi$ with this choice, although this has not yet been proved. 
The current best results proving accuracy of the Galerkin solutions for large $k$ for the Dirichlet problem are in \cite{GaMuSp:19} (following \cite{GrLoMeSp:15}), with these results proving quasioptimality of the Galerkin solution (with quasioptimality constant independent of $k$) 
(i) for smooth and strictly convex $\Oi$ when $hk^{4/3}$ is sufficiently small, and (ii) 
for general nontrapping $\Oi$ when
$hk^{3/2}\log k$ is sufficiently small.
There is almost no analogous theory for the Neumann problem for large $k$; the exception is \cite{BoTu:13} whose results about coercivity 
of the BIE \eqref{eq:direct_combined_neu2} when $\Oi$ is a ball imply a quasioptimality result without any restriction on $h$, albeit with quasioptimality constant growing like $k^{1/3}$.

\subsubsection{Iterative solution of the BEM linear systems via GMRES.} 
\label{sec:GMRES}

A popular way of solving the dense linear systems that arise from the BEM is via iterative methods \cite[Chapter 13]{St:08}, \cite[Chapter 6]{SaSc:11}, \cite[\S4]{RjSt:07}. Since the systems arising from the Helmholtz equation are, in general, non-normal (as highlighted in \S\ref{sec:BIE}), a natural choice of iterative method is the generalised minimum residual method (GMRES) \cite{SaSc:86}.

Given \(\genmatrix \in \bbC^{n\times n}\), \(\bfb\in \bbC^n\), the generalised minimum residual method (GMRES) to find the solution $\bfx$ of $\genmatrix \bfx = \bfb$ is the following. 
Given  \(\bfx_0 \in \bbC^n\), let  \(\bfr_0(\genmatrix,\bfb,\bfx_0):=\bfb-\genmatrix \bfx_0\). 
Let the Krylov space \(\calK_m(\genmatrix,\bfr_0)\) be defined by
\begin{align*}
    \calK_m(\genmatrix,\bfr_0):=\Span \left\{\bfr_0, \genmatrix\bfr_0, \ldots, \genmatrix^{m-1}\bfr_0\right\}.
\end{align*}
The $m$th iterate of GMRES, \(\bfx_m\), is defined as the unique vector 
in $\bfx_0 + \calK_m(\genmatrix,\bfb)$ that minimises the residual 
$\bfr_m:=\bfb-\genmatrix \bfx_m$ with respect to the \(\lVert \cdot \rVert_2\) norm (see, e.g., \cite[\S6.5.1]{Sa:03}).
Observe that, since \(\bfx_m \in \bfx_0 + \calK_m(\genmatrix,\bfb)\), the residual satisfies
\beq\label{eq:residual}
\bfr_m:=\bfb-\genmatrix \bfx_m= p_m(\genmatrix)\bfr_0
\quad\tfor p_m \in \bbP_m \text{ with } p_m(0)=1,
\eeq
where $\mathbb{P}_m$ denotes the set of polynomials of degree $m$. The definition of GMRES therefore implies that
\begin{align}\label{def:gmres}
    \lVert \bfr_m (\genmatrix,\bfb,\bfx_0) \rVert_2= \min_{\substack{p_m\in \bbP_m, \\ p_m(0)=1}}\lVert p_m (\genmatrix)\bfr_0(\genmatrix,\bfb,\bfx_0)\rVert_2.
\end{align}
We apply GMRES to the linear system  \eqref{eq:1} preconditioned by the mass matrix  
      \begin{align}\label{eq:massmatrix}
    (\bfM)_{i,j}:= \int_{\Gamma} \phi_i(\bx) \phi_j(\bx) \dif \sigma(\bx);
\end{align}
i.e.~we solve 
\beq\label{eq:2}
\bfM^{-1}\matrixD \bv = \bfM^{-1}\bff.
\eeq
We solve \eqref{eq:2} instead of \eqref{eq:1} because it is easier to translate information about $A_k'$ to information about $\bfM^{-1}\matrixD $ rather than information about $\matrixD$. This is for the following two reasons.

(a) The eigenvalues of $\bfM^{-1}\matrixD $ approximate the eigenvalues of $A_k'$. Indeed, the eigenvalue problem $A_k' v=\lambda_k v$ is equivalent to the variational problem: find $v\in \LtG$ such that $(A_k' v,w)_\LtG= \lambda_k (v,w)_{\LtG}$ for all $w\in \LtG$, and the Galerkin approximation of this is $\matrixD \bv = \lambda_k \bfM \bv$.

(b) If \(\Gamma\) is \(C^1\), then, given $k$, there exists \(h_0\) and \(\Capproxone\) such that, if \(h\leq h_0\) then
\beq\label{eq:discon2}
(\Capproxone^{-1})\N{A_k'}_{\LtGt}\leq
\N{\bfM^{-1}\matrixD}_2 \leq \Capproxone \N{A_k'}_{\LtGt},
\eeq
and analogous bounds hold for $(A_k')^{-1}$ and $(\bfM^{-1}\matrixD)^{-1}$
(furthermore, if the basis is orthonormal, then $\|\bfM^{-1}\matrixD\|_2\rightarrow \|A_k'\|_{\LtG}$
and $\|(\bfM^{-1}\matrixD)^{-1}\|_2\rightarrow \|(A_k')^{-1}\|_{\LtG}$ as $h\to 0$ for fixed $k$); see Lemma \ref{lem:discon} and Remark \ref{rem:compact} below.
In contrast, in the analogue of \eqref{eq:discon2} with $\bfM^{-1}\matrixD$ replaced by $\matrixD$, the constant $C$ depends on $h$; see \eqref{eq:Galerkinmatrixbounds2}.

We only consider solving the system \eqref{eq:2} with standard GMRES because our goal is to prove rigorous bounds on the number of iterations and the theory of GMRES convergence is most well-developed for standard GMRES.
We note that GMRES is often used with either restarts or restarts with subspace augmentation (see, e.g., \cite{Mo:95, Mo:02, GiGrPiVa:10}) -- this has the advantage of reducing storage and orthogonalisation costs, but with the number of iterations required to obtain a given relative residual necessarily higher than for standard GMRES (although it is difficult to study this increase theoretically).

\subsection{Four features (F1-F4) observed in numerical experiments on the set-up in \S\ref{sec:statement}, and statement of the main goals of this paper}\label{sec:features}

We now highlight four different features one observes from computing approximations to the scattering problem via the set-up in \S\ref{sec:statement} (i.e.~reformulating as a BIE, creating a linear system via the BEM, and solving the linear system using GMRES). We present numerical experiments illustrating each of the features later in the paper.

These features are about, respectively, (1) the accuracy of the Galerkin solution, (2) the condition number of the Galerkin matrix, (3) the number of GMRES iterations, and (4) the accuracy of the GMRES solution. 
\bit
\item[F1] When the incoming plane wave enters the cavity, one needs a larger number of points per wavelength for accuracy of the Galerkin solutions than when the wave doesn't enter the cavity. 

\item[F2] The norm of $(\bfM^{-1}\matrixD)^{-1}$ (i) is very sensitive to whether or not $k=k_j$ for $k_j$ a quasimode frequency, and (ii)
grows exponentially through $k_j$, 
up to some point, and then grows more slowly.

\item[F3] The number of GMRES iterations required to make the residual $\bfr_m$ arbitrarily small 
\bit
\item[(a)] grows algebraically with $k$, with no worse growth through $k=k_j$ than $k\neq k_j$,
\item[(b)] depends on whether $\Oi$ is the small or large cavity, and 
\item[(c)] depends on the direction of the incoming plane wave.
\eit
\item[F4] The GMRES residual being small does not necessarily mean that the error is small, and the relative sizes of the residual and error depend on both $k$ and the direction of the plane wave.
\eit
{\bf The main goals of this paper are to explain F3(a) and, to a certain extent,
F3(b).} 
The following is an outline of the rest of the introduction. In \S\ref{sec:intro_num} we present numerical results about F3(a), F3(b), and F3(c) for the Dirichlet problem. 
In \S\ref{sec:observations} we give plots of the eigenvalues of $\bfM^{-1}\matrixD$.
In \S\ref{sec:theorem} we give a general bound on the number of GMRES iterations when the matrix has a ``cluster plus outliers'' structure. In \S\ref{sec:mainresultH} we apply the bound from \S\ref{sec:theorem} to $\bfM^{-1}\matrixD$, under assumptions based on the eigenvalue plots in \S\ref{sec:observations}; the result is a $k$-explicit bound on the number of GMRES iterations which explains F3(a).
In the last section of the introduction, \S\ref{sec:FEM}, we discuss how the ideas in this paper can be applied to the Helmholtz FEM.
The partial explanation of F3(b) is contained in \S\ref{sec:smalllarge}.

Although our focus is on F3, we still need to be aware of other features; e.g.~the solution obtained by GMRES is useless if the Galerkin solution itself isn't accurate (F1), or if the GMRES solution isn't close to the Galerkin solution (F4).
We therefore make some brief comments here about to what extent the features F1, F2, F3(c), and F4 are rigorously understood; the summary is that F2 is rigorously understood, whereas F1, F3(c), and F4 are not. 

\paragraph{Regarding F1:} 
The fact that the accuracy of the Galerkin solution depends on whether or not the wave enters the cavity makes physical sense, but there is currently no rigorous theory on the subject. Indeed, as discussed in \S\ref{sec:statement:BEM}, the current best analysis of how $h$ must depend on $k$ for the $h$-BEM to be uniformly accurate as $k\tendi$, \cite{GaMuSp:19}, is not sharp in the nontrapping case, and
is therefore very far from proving rigorous sharp results about the trapping case.  
Numerical experiments illustrating F1 are given in Appendix \ref{app:F1F4}.

\paragraph{Regarding F2:} 
The exponential growth of $(\bfM^{-1}\matrixD)^{-1}$ through the sequence of $k_j$ is explained by the following.
The inverses of the boundary integral operators $A'_k$, $B_k$, and $B_{k, {\rm reg}}$ inherit the behaviour of the Helmholtz solution operator, and thus grow when $k=k_j$ for $k_j$ quasimode frequencies. More precisely, if $k_j$ and $\QMC(k)$ are as in Definition \ref{def:quasimodes}, then there exists $C>0$ (independent of $j$) such that 
\beq\label{eq:Ainvblowup}
\N{(A'_{k_j})^{-1}}_{\LtGt} \geq C \left(\frac{1}{\QMC(k_j)} - \frac{1}{k_j}\right)k_j^{1/4} \quad\tfa j
\eeq
\cite[Equation 5.39]{ChGrLaSp:12} \footnote{More precisely, \cite[\S5.6.2, Equation 5.39]{ChGrLaSp:12} proves \eqref{eq:Ainvblowup} with a different power of $k$ on the right-hand side. 
The bound \eqref{eq:Ainvblowup} can be proved by following the same steps as in 
\cite[\S5.6.2]{ChGrLaSp:12}, but using the sharp bound on the single-layer potential from~\cite[Theorem 1.1, Part (i)]{HaTa:15}.}. Therefore, when $\Oi$ is either the small or the large cavity, by \eqref{eq:quality}, 
$\|(A'_{k_j})^{-1}\|_{\LtGt}\geq C_1 \exp (C_2 k_j)$ for some $C_1, C_2>0$ independent of $j$
\cite[Theorem 2.8]{BeChGrLaLi:11}.
We note that, since $\|A'_{k}\|_{\LtGt}$ grows algebraically in $k$ for general Lipschitz domains (see 
Part (i) of Lemma \ref{lem:normbounds} below), the condition number of $A'_{k_j}$ also grows exponentially as $j\tendi$ when $\Oi$ is either the small or the large cavity.

An indication (but not a rigorous proof) of why the growth of $\|(\bfM^{-1}\matrixD)^{-1}\|_2$ through $k=k_j$ stagnates,
and why $\|(\bfM^{-1}\matrixD)^{-1}\|_2$ is very sensitive to whether or not $k=k_j$,
is given
by the recent result of \cite[Theorem 1.1]{LaSpWu:20}. This result shows that, for most frequencies, the Helmholtz solution operator (and hence also $\|(A'_{k_j})^{-1}\|_{\LtGt}$) is bounded polynomially in $k$. More precisely, 
given 
$\smallpara>0$, 
there exists $C_3=C_3(\smallpara)>0$ and a set $J\subset [1,\infty)$ with $|J|\leq \smallpara$ such,
\beq\label{eq:Ainv}
\N{(A'_{k})^{-1}}_{\LtGt}  \leq C_3 k^{5d/2 +1} \quad\tfa k \in [1,\infty)\setminus J
\eeq
(for simplicity we have assumed $k\geq 1$, but an analogous bound holds for any $k\geq k_0>0$).
The bounds \eqref{eq:Ainvblowup} and \eqref{eq:Ainv} then imply that the graph of $\|(A'_{k_j})^{-1}\|_{\LtGt}$ against $k$ consists of a number of ``spikes'' at $k=k_j$, with the heights of the spikes growing exponentially with $k$, but the widths decreasing  with $k$
\footnote{For an illustration of this in a simple 1-d model of resonance behaviour, see \cite[\S23.2]{FeLeSa:64}.}.
Therefore, while
$\|(A'_{k_j})^{-1}\|_{\LtGt}$ grows exponentially through $k_j$, the growth is very sensitive to the precise value of $k$ (with this sensitivity increasing as $k$ increases). 
This result indicates that the growth of $(\bfM^{-1}\matrixD)^{-1}$ through $k=k_j$ stagnates since discretisation error collapses the delicate exponential growth. 

\paragraph{Regarding F3(c):} This feature arises because the GMRES residual $\bfr_m$ \eqref{eq:residual} depends on the right-hand side vector, which depends on the direction of the plane wave (via the right-hand side of the BIE in \eqref{eq:direct_combined_dir}). There are few rigorous results in the literature describing 
the dependence of $\bfr_m$ on the right-hand side vector, but 
in Appendix \ref{app:GMRES_rhs} we describe how the results of  \cite{Titley2014} give a heuristic explanation 
of this feature 
for problems with similar eigenvalue distributions to the Helmholtz problems we consider.

\paragraph{Regarding F4:} 
Numerical experiments illustrating this feature on our problem are given in Appendix \ref{app:F1F4}. 
This feature is poorly understood for non-normal, complex linear systems in general, and thus also for the systems arising from the Helmholtz problems considered here. 
There have been many papers that discuss the convergence of GMRES in the sense of residual reduction; in contrast, there is remarkably little known in the literature about the error $\bfx-\bfx_m$. The most recent (and most relevant) results in this area are given in \cite{Meurant2011} and \cite[\S 5.8]{Meurant2020} and even then the results are stated for real systems. In particular, \cite[Theorem  5.35 and Corollary 5.6]{Meurant2020} gives bounds on $\|\bfx-\bfx_m\|_2$ in terms of $\|\br_m\|_2$ multiplied by a computable expression that requires the existence of $\bfH_m^{-1}$ (where $\bfH_m$ is the square $m \times m$ upper Hessenberg matrix arising at the $m$-th step in the Arnoldi process in the GMRES algorithm) and also depends on other entries of $\bfH_m$. 

\subsection{Numerical experiments about F3}\label{sec:intro_num}

For ease of exposition, we only present here experiments for the Dirichlet problem \eqref{eq:direct_combined_dir}, i.e., involving the operator $A_k'$, in 2-d.
\S\ref{sec:experiments} contains experiments for the two BIEs for the Neumann problem, \eqref{eq:direct_combined_neu} and \eqref{eq:direct_combined_neu2}, and experiments for the Dirichlet problem in 3-d.

All the experiments in this section use ten points per wavelength.
Furthermore, we plot quantities of interest (such as the number of iterations, the condition number) through \emph{either} integer values of $k$ \emph{or} values of $k$ in a quasimode. As described in \S\ref{sec:statement:ellipse}, 
for $\Oi$ the small or large cavities, 
there exists a quasimode with frequencies equal 
$\{k^{e/o}_{m,n}\}_{m=1}^\infty$ for fixed $n$.
Our experiments consider $\{k_{m,0}^e\}_{m=1}^\infty$, but we observe very similar behaviour through 
$\{k_{m,n}^{e/o}\}_{m=1}^\infty$ for $n\neq 0$ fixed.

\paragraph{Experiments illustrating F3(a) (growth of iterations with $k$).}

\begin{figure}[h!]
    \centering
    \includegraphics[width=0.9\textwidth]{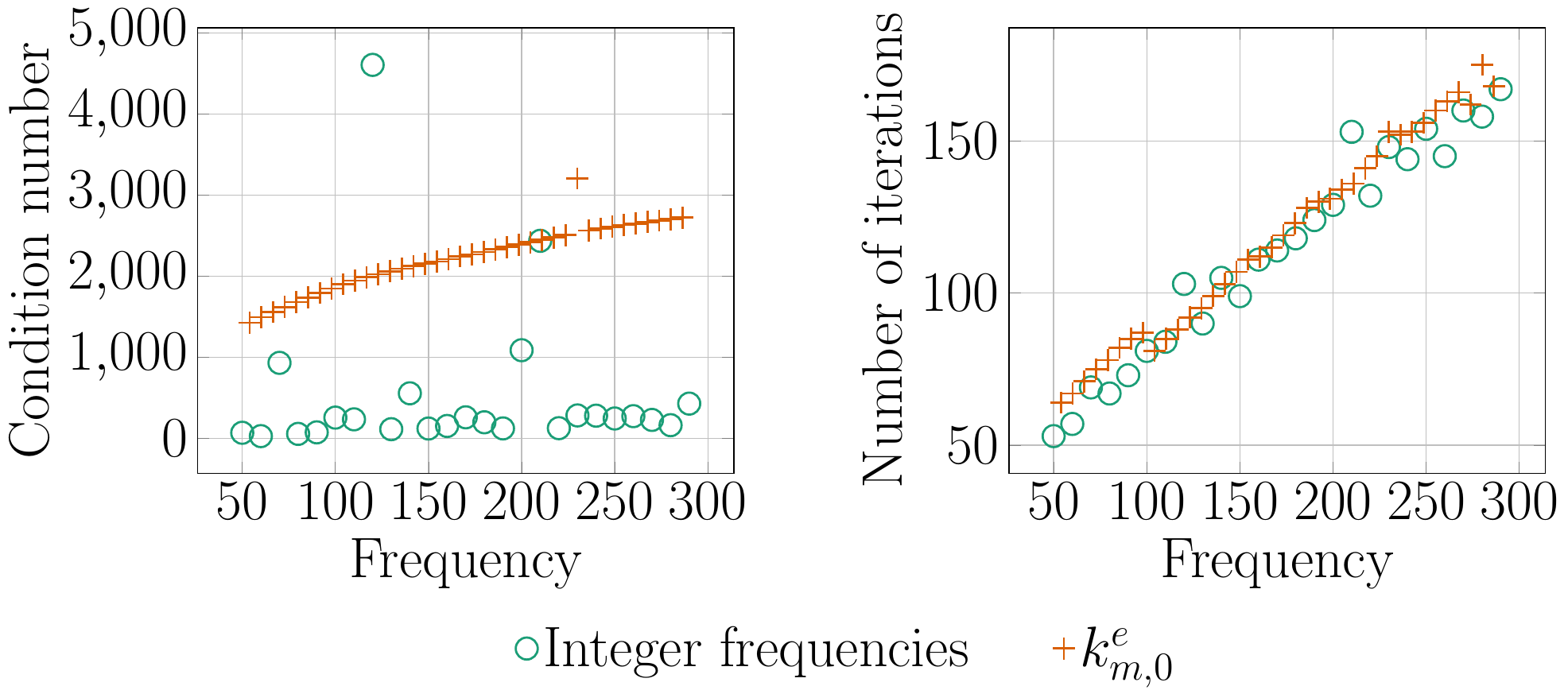}
    \caption{Condition number and number of GMRES iterations for small cavity with incident plane wave at angle $\theta =4\pi/10$ to the horizontal (illustrating F3(a)). 
       }\label{fig:comparison_spectrum_it}
\end{figure}

Figure~\ref{fig:comparison_spectrum_it} 
plots the condition number of $\bfM^{-1}\matrixD$ and the number of GMRES iterations against $k$ for the small cavity.
The direction $\ba$ of the incident plane wave $\exp(\ri k \bx\cdot\ba)$ is chosen as $a=(\cos\theta, \sin\theta)$ with 
$\theta =4\pi/10$; from Figure \ref{fig:geometries} we see that the plane wave is almost vertical and enters the cavity.

The key point from Figure~\ref{fig:comparison_spectrum_it} is that, while the condition number $\bfM^{-1}\matrixD$ is very sensitive to whether or not $k$ 
is near a frequency in the quasimode, the number of iterations is not. This demonstrates the well-known fact that the condition number gives little insight into the behaviour of GMRES for non-normal problems. 

In more detail, the left-hand plot in Figure~\ref{fig:comparison_spectrum_it} shows, via the condition number, the sensitivity of $\|(\bfM^{-1}\matrixD)^{-1}\|_2$ to whether or not $k=k_j$ (i.e., the first point in F2). 
The green circle outlier at $k=120$ is there because, by chance, the integer frequency 120 lies very close to the quasimode frequency
$k^o_{17,7}=119.997615771724$ 
 (note that to 5 significant figures this approximation of the quasimode frequency is equal to $120$).
This plot also shows the growth of $\|(\bfM^{-1}\matrixD)^{-1}\|_2$ through $k=k_j$ stagnating as $k$ increases (i.e., the second point in F2); this was also seen in 
the experiments in \cite[Section IV.H]{BeChGrLaLi:11}  on the small cavity, where, even using 20 points per wavelength, 
  the exponential growth of $\|(\bfM^{-1}\matrixD)^{-1}\|_2$ through $k=k_j$ levelled off after $k=60$.
  This sensitivity of $\|(\bfM^{-1}\matrixD)^{-1}\|_2$ to whether or not $k=k_j$ was also shown    in \cite[Figure 4.7]{LoMe:11}
  for a cavity similar to both our large and small cavities.

The plot of the number of iterations in Figure~\ref{fig:comparison_spectrum_it} is also included in Figure \ref{fig:iterations_main} below, where we see that the number of iterations grows like $k^{0.66}$ (in the range considered).

\begin{figure}[h!!]
    \centering
    \includegraphics[width=0.9\textwidth]{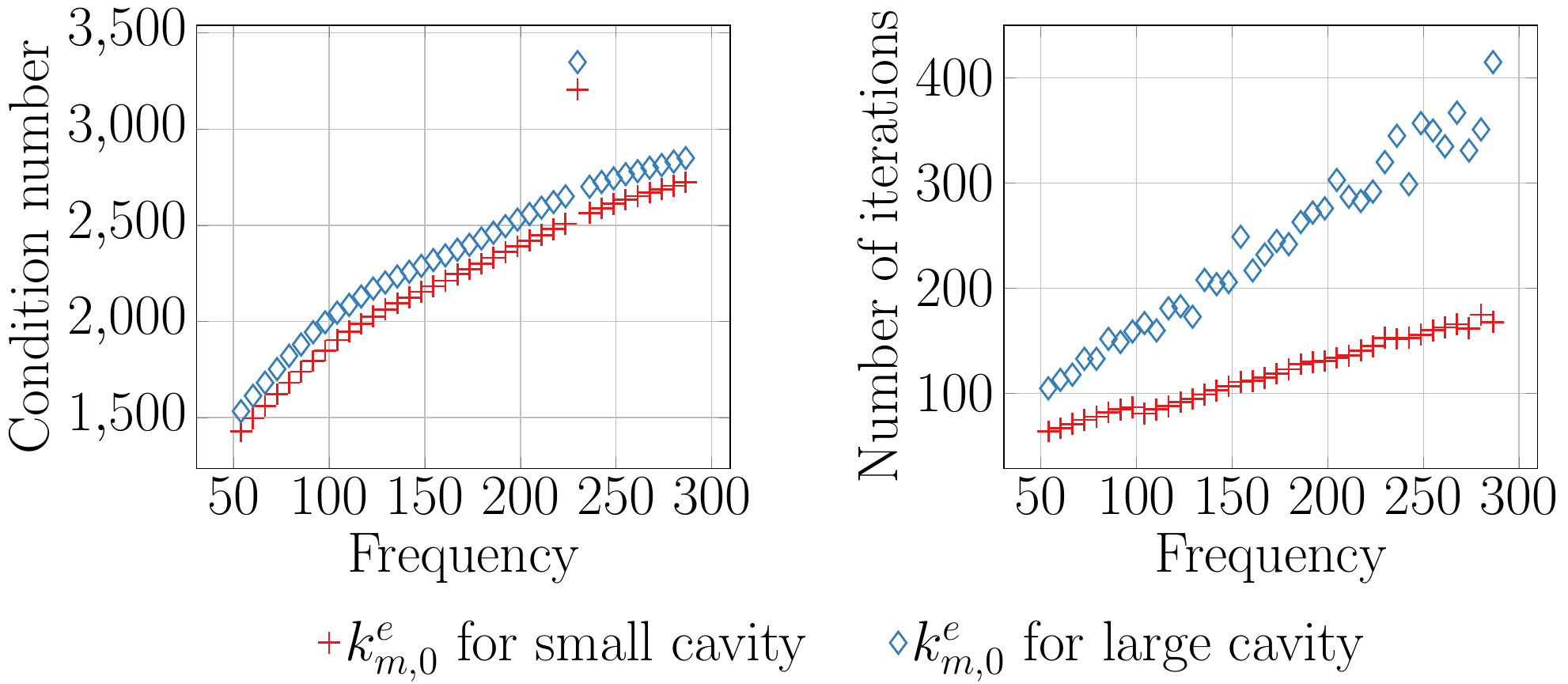}
    \caption{Comparison of the condition number and the number of GMRES iterations for the small and large cavities 
    with incident plane wave at angle $\theta =4\pi/10$ to the horizontal (illustrating F3(b)).}\label{fig:comparison_size_cavit}
\end{figure}
\begin{figure}[h!!]
    \centering
    \includegraphics[width=0.9\textwidth]{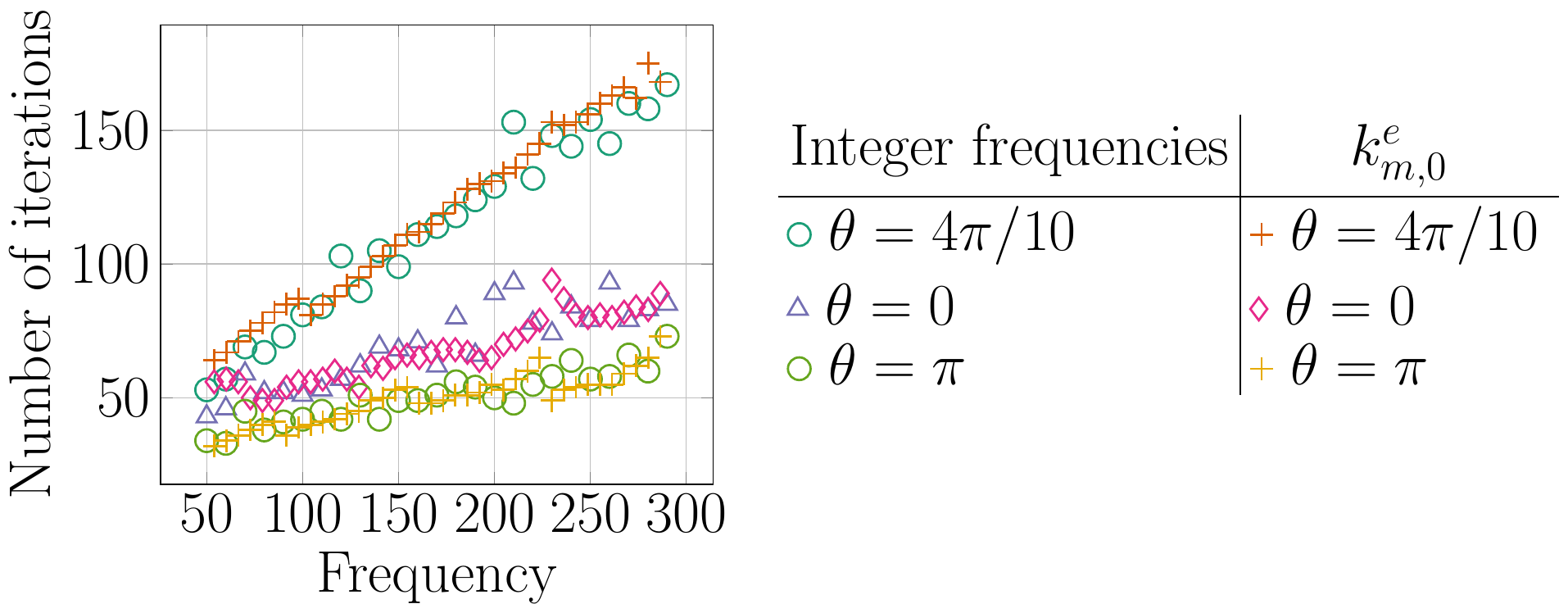}
    \caption{The number of GMRES iterations for the small cavity and three different right-hand sides, corresponding to three different angles $\theta$ of the incident plane wave (illustrating F3(c)).}\label{fig:comparison_rhs}
\end{figure}

\paragraph{Experiments illustrating F3(b) (dependence of iterations on cavity size).}

Figure~\ref{fig:comparison_size_cavit}
plots the condition number of $\bfM^{-1}\matrixD$ and the number of GMRES iterations against $k_{m,0}^e$ for both the small and large cavities.
As in the previous figure, $\ba= (\cos\theta, \sin\theta)$ with 
$\theta =4\pi/10$; i.e. the plane wave is almost vertical and enters the cavity.
While the condition numbers behave very similarly,
the growth in the number of iterations is different, again illustrating the fact that the condition number is not relevant for understanding the convergence of GMRES for non-normal matrices. For the small cavity the number of iterations grows approximately like $k^{0.66}$, and for the large cavity like $k^{0.82}$; see Figure \ref{fig:iterations_main} below.

\paragraph{Experiments illustrating F3(c) (dependence of iterations on plane-wave direction).}

Figure~\ref{fig:comparison_rhs}
plots the number of GMRES iterations against $k$ for the small cavity and varying $\theta$, with the incident plane wave $\ba= (\cos\theta, \sin\theta)$. From Figure \ref{fig:geometries}, we see that when $\theta =4\pi/10$ the  plane wave is almost vertical and enters the cavity, when 
$\theta =0$ the plane wave is horizontal and enters the cavity, and when $\theta=\pi$ the plane wave is horizontal and doesn't enter the cavity.
Physically, we  expect the worst behaviour to occur when $\theta=4\pi/10$, because of the multiple reflections in the cavity, and the best behaviour to occur when $\theta=\pi$, and this is indeed what we see in Figure \ref{fig:comparison_rhs}.

\paragraph{Links with other experiments/results in the literature.}
Both the iterative solution of BEM linear systems 
and solving scattering problem involving cavities
have received a lot of interest in the literature; see, e.g., the books \cite[Chapter 13]{St:08}, \cite[Chapter 6]{SaSc:11}, \cite[\S4]{RjSt:07} for the former, and, e.g., \cite{BaSu:05, WaDuSu:09, GmPh:07, GiGrPiVa:10, Du:11, DaDaLa:13, LaAmGr:14, LaGrON:17} for the latter.
Nevertheless, the features F1-F4 do not appear to have been systemically identified and studied before now. 

We highlight here one previous study where the features F2 and F3(a) are visible in numerical experiments. Indeed, 
\cite{DaDaLa:13} considers solving the Neumann problem with a BIE similar to \eqref{eq:direct_combined_neu2}, but 
with $S_{\ri k}$ replaced by a different regularising operator. 
The figures \cite[Figures 11(a), 18, and 19]{DaDaLa:13} plot the condition number against $k$ when $\Oi$ 
are cavity domains similar to those in Figure \ref{fig:geometries} (although supporting weaker trapping), 
and display spikes; i.e., F2.
The figure \cite[Figure 28(a)]{DaDaLa:13} plots the number of GMRES iterations against $k$ and sees growth with no spikes; i.e., F3(a).

\subsection{Plots of the eigenvalues of $\bfM^{-1}\matrixD$.}\label{sec:observations}

\begin{figure}[h!]
    \centering
    \begin{subfigure}[t]{0.45\textwidth}
        \centering
        \includegraphics[width=0.9\textwidth]{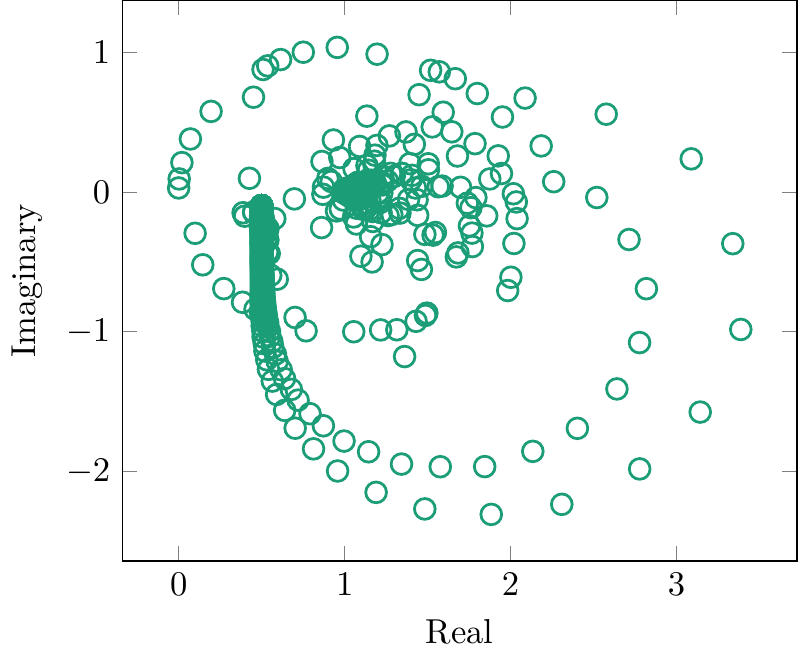}
        \caption{Small cavity at \(k=100\)}
    \end{subfigure}\hfill%
    \begin{subfigure}[t]{0.45\textwidth}
        \centering
        \includegraphics[width=0.9\textwidth]{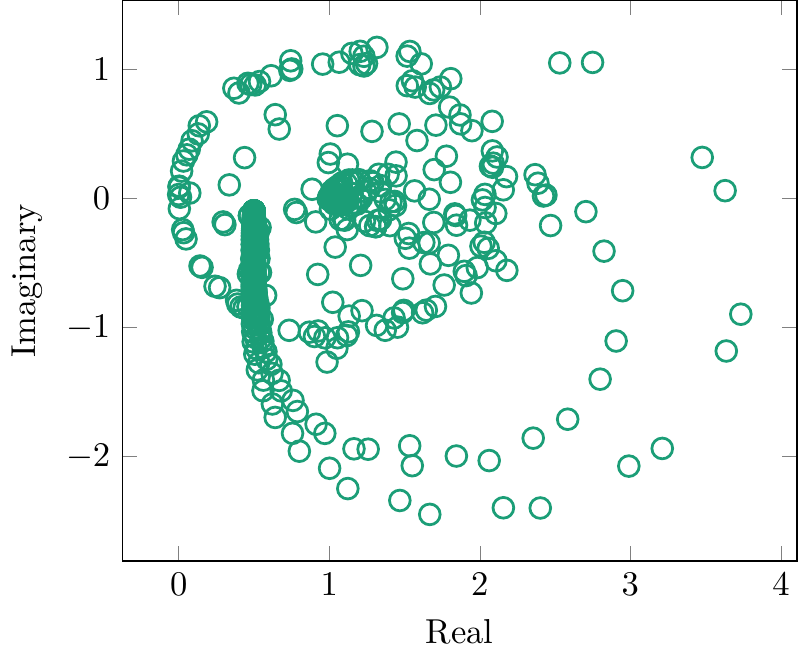}
        \caption{Large cavity at \(k=100\)}
    \end{subfigure}%
    \par\bigskip
    \begin{subfigure}[t]{0.45\textwidth}
        \centering
        \includegraphics[width=0.9\textwidth]{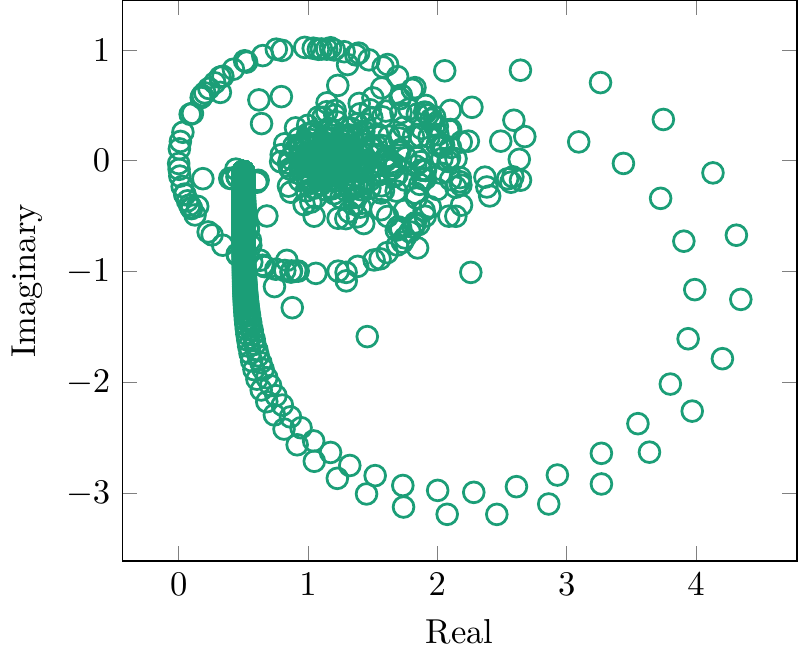}
        \caption{Small cavity at \(k=290\)}
    \end{subfigure}\hfill%
    \hfill
    \begin{subfigure}[t]{0.45\textwidth}
        \centering
        \includegraphics[width=0.9\textwidth]{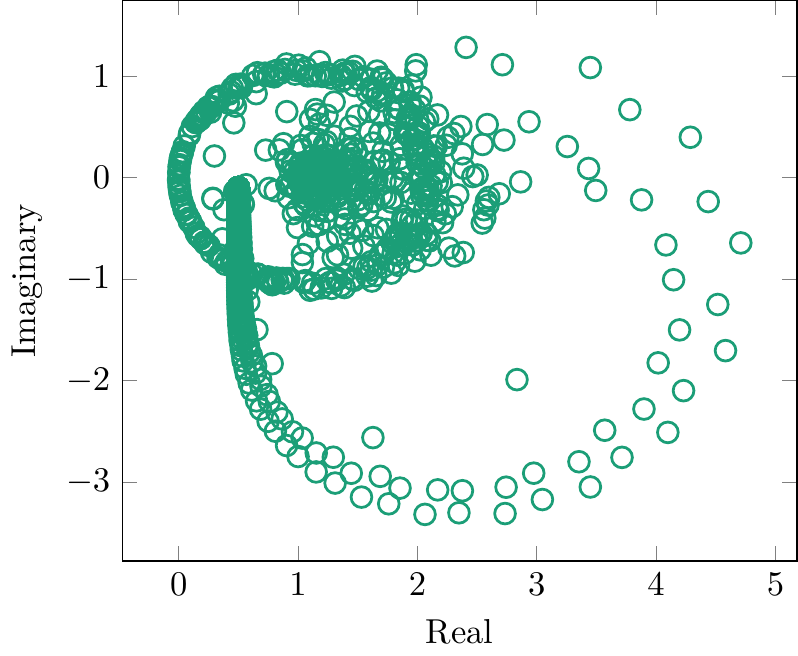}
        \caption{Large cavity at \(k=290\)}
    \end{subfigure}%
    \caption{The eigenvalues of $\bfM^{-1}\matrixD$ for the small and large cavities at $k=100$ and $k=290$.}\label{fig:example_spectrum}
\end{figure}

\paragraph{Summary of the figures.} Figure \ref{fig:example_spectrum} plots the eigenvalues of $\bfM^{-1}\matrixD$ for the small and large cavities at $k=100$ and $k=290$.

\begin{figure}[h!]
    \centering
    \includegraphics[width=0.6\textwidth]{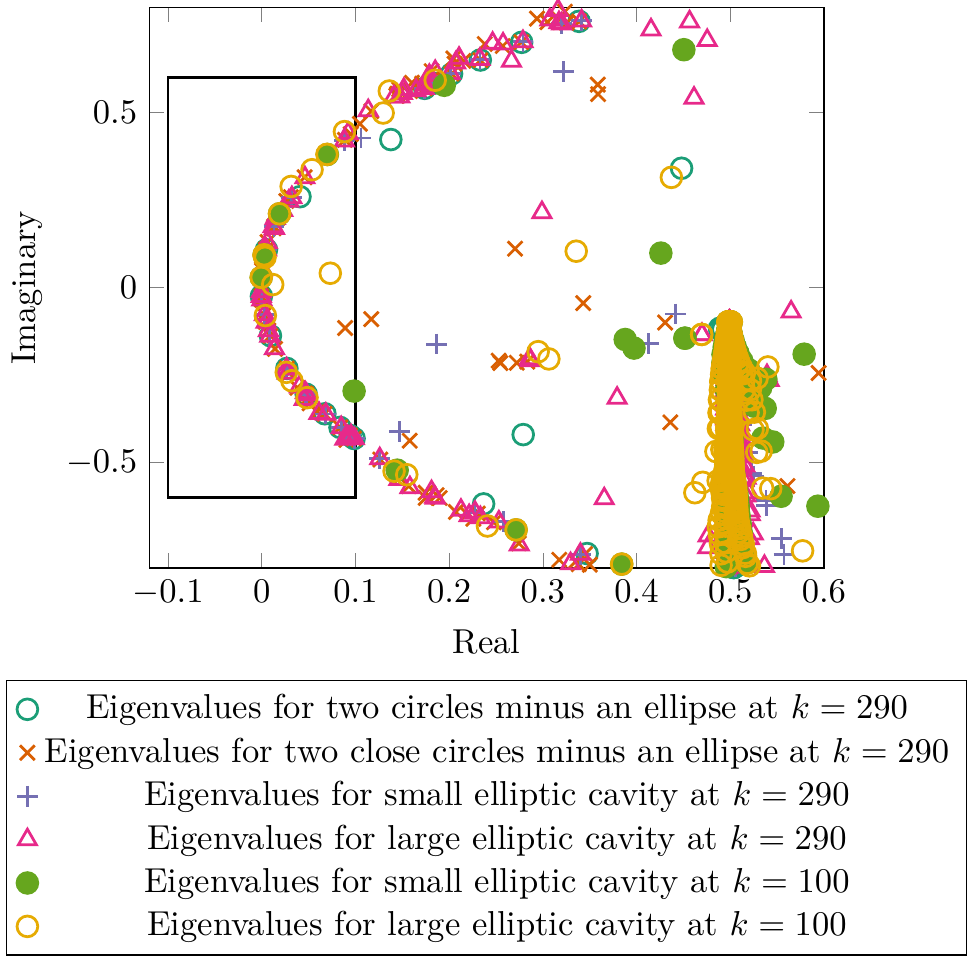}
    \caption{Plots of the near-zero eigenvalues of $\bfM^{-1}\bfA_k'$ for a variety of different domains and frequencies. The black rectangle $[-0.1,0.1]\times[-0.6,0.6]$ is a choice of the set $\cN$ in Observation O1.}  
    \label{fig:outliers}
\end{figure}

\begin{figure}[h!]
    \centering
    \begin{minipage}{0.49\textwidth}
        \centering
        \includegraphics[width=0.9\textwidth]{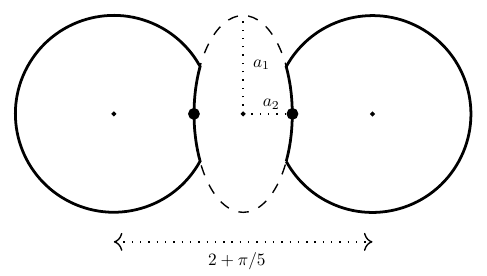}
    \end{minipage}%
    \hfill
    \begin{minipage}{0.49\textwidth}
        \centering
        \includegraphics[width=0.7\textwidth]{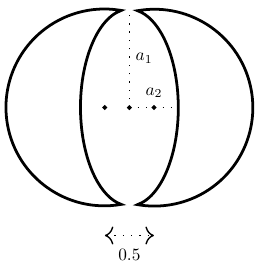}
    \end{minipage}
    \caption{The ``two circles minus ellipse'' and ``two close circles minus ellipse'' $\Oi$ considered in Figure \ref{fig:outliers}
     with $a_1=1$ and $a_2=1/2$.
    }
    \label{fig:circles_minus_ellipse}
\end{figure}

Figure ~\ref{fig:outliers} plots the eigenvalues for the small and large cavities at both $k=100$ and $k=290$, as well as the eigenvalues at $k=290$ for two other $\Oi$ for which Theorem \ref{thm:ellipse} applies; these two other $\Oi$ are plotted in Figure \ref{fig:circles_minus_ellipse}.

\begin{figure}[h!]
    \centering
    \includegraphics[width=0.9\textwidth]{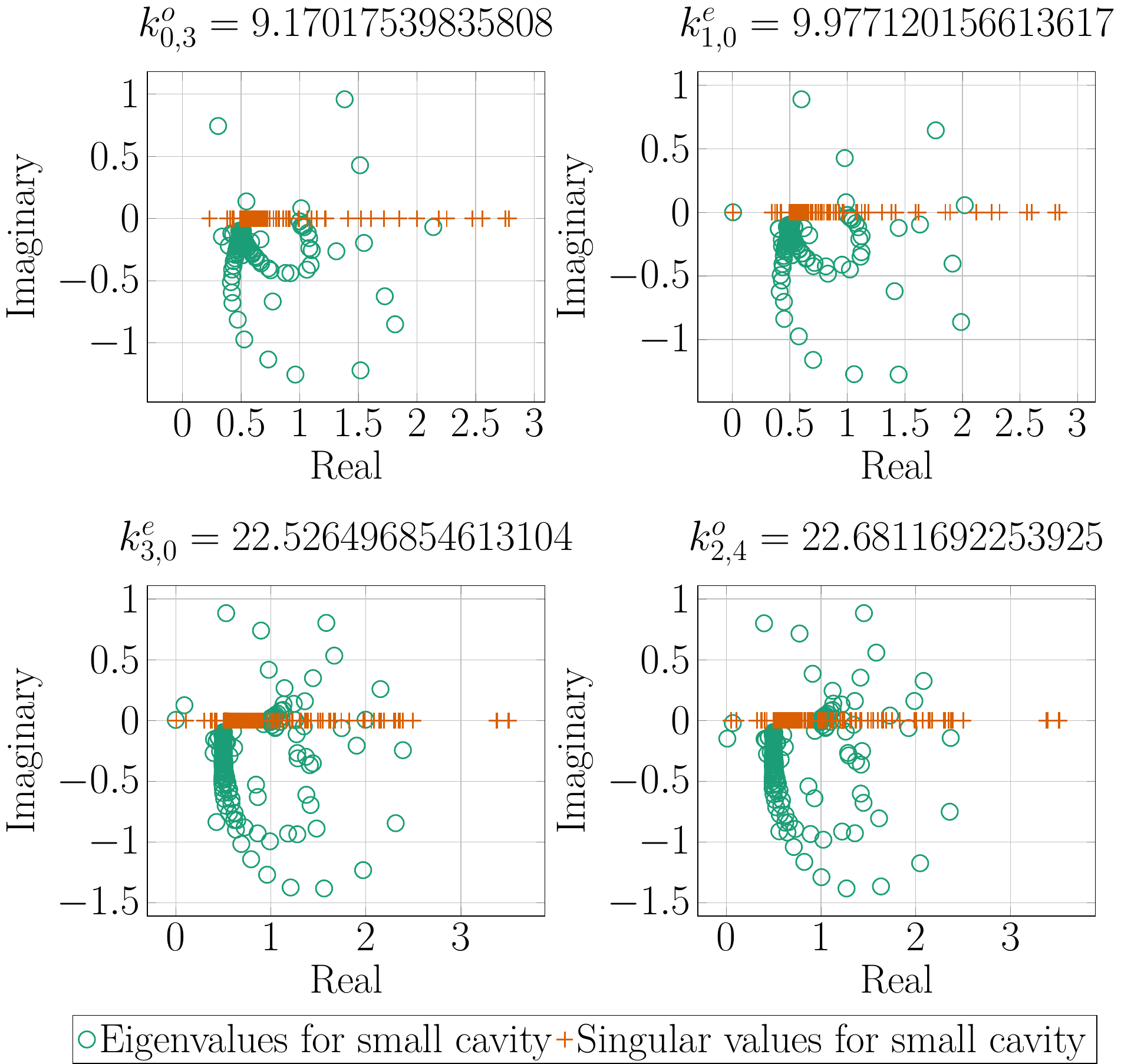}
    \caption{The eigenvalues and singular values of $\bfM^{-1}\matrixD$ for several frequencies $k^{e/o}_{m,n}$ and $\Oi$ the small cavity}\label{fig:spectra_illustrations_elliptic_cavity_2D_dir}
\end{figure}

\begin{figure}[h!]
    \centering
    \includegraphics[width=0.9\textwidth]{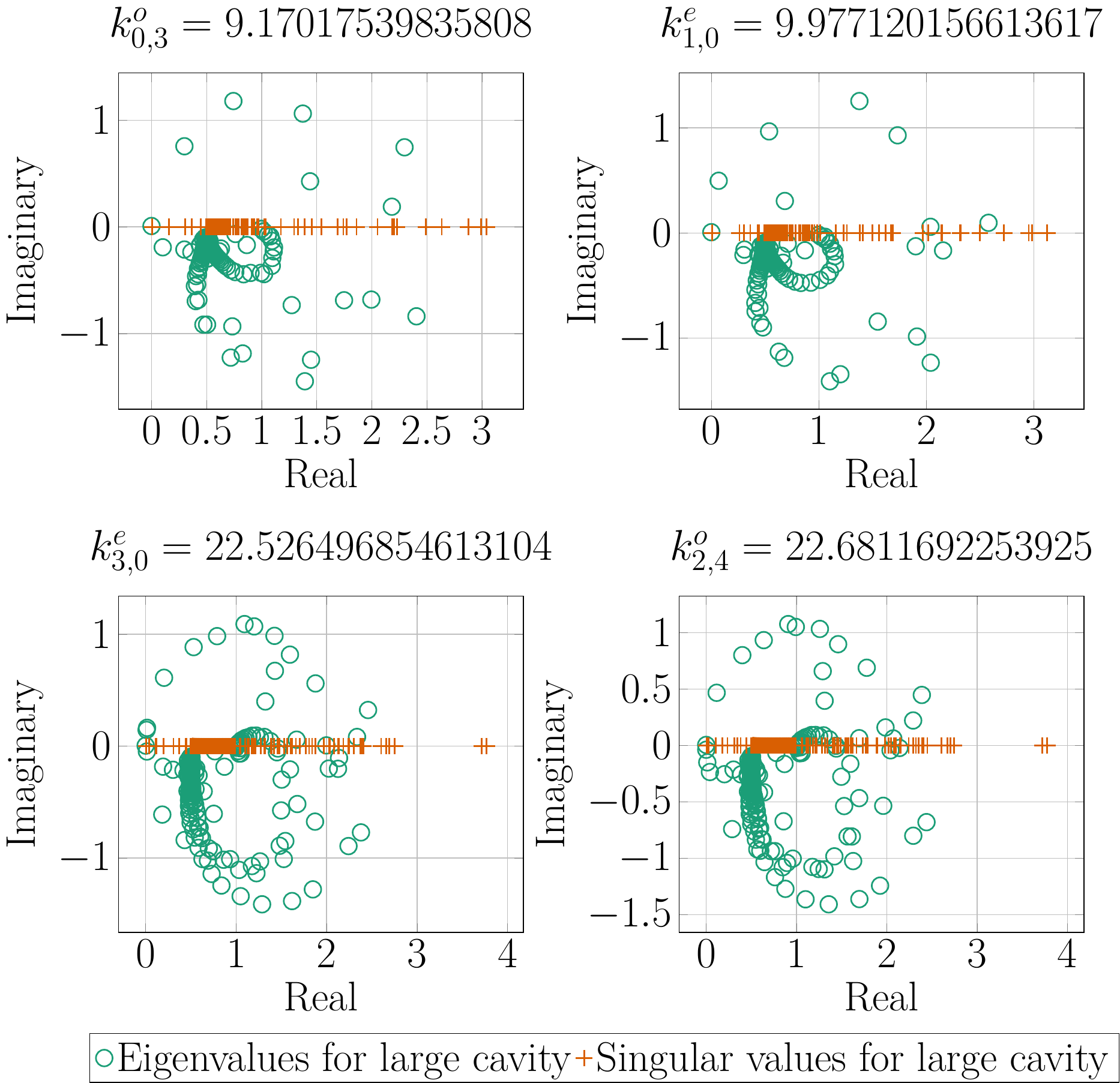}
        \caption{The eigenvalues and singular values  of $\bfM^{-1}\matrixD$ for several frequencies $k^{e/o}_{m,n}$ and $\Oi$ the large cavity}
    \label{fig:spectra_illustrations_elliptic_cavity_bigger_2D_dir}
\end{figure}

Figures \ref{fig:spectra_illustrations_elliptic_cavity_2D_dir} and \ref{fig:spectra_illustrations_elliptic_cavity_bigger_2D_dir}
plot the eigenvalues and singular values of $\bfM^{-1}\matrixD$ for several frequencies $k^{e/o}_{m,n}$ and $\Oi$ the small and large cavities, respectively.

\begin{figure}[h!]
    \centering
    \begin{minipage}{0.49\textwidth}
        \centering
        \includegraphics[width=0.9\textwidth]{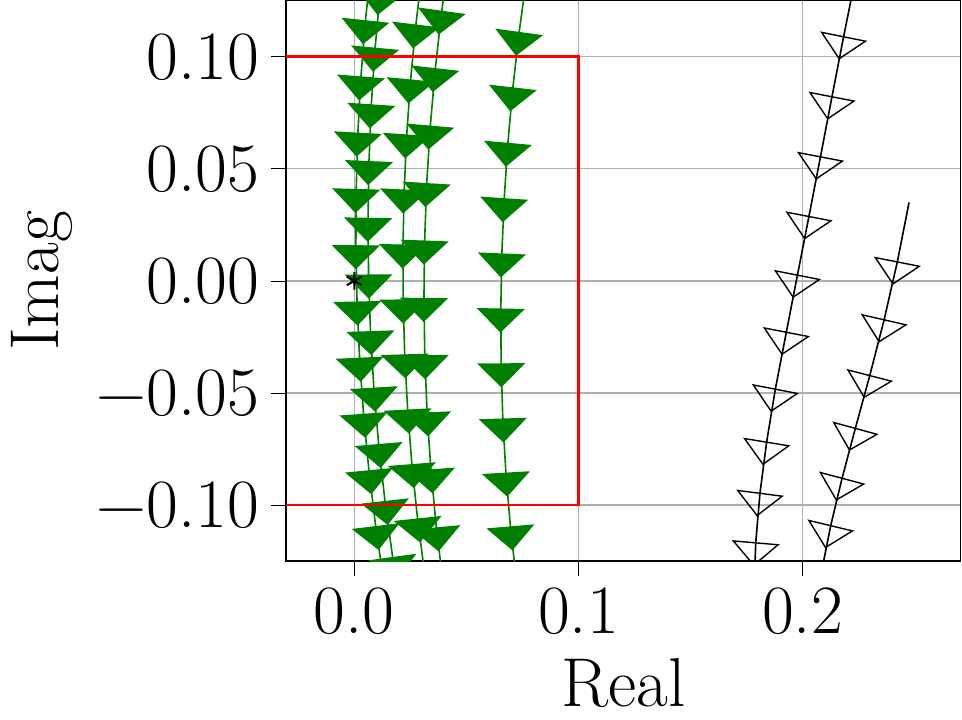}
    \end{minipage}%
    \hfill
    \begin{minipage}{0.49\textwidth}
        \centering
        \includegraphics[width=0.9\textwidth]{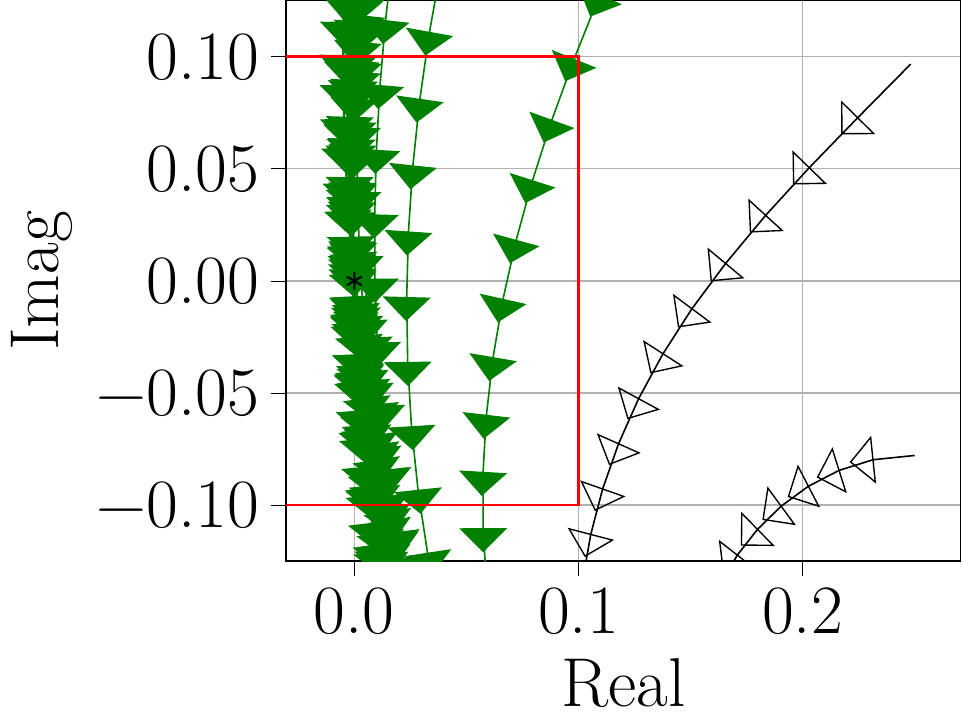}
    \end{minipage}
    \caption{Paths of the eigenvalues for \(k\in (5,15)\) for the small cavity (left) and the large cavity (right) 
    The eigenvalues that enter the  rectangle 
    are coloured green and have shaded arrowheads.}
    \label{fig:flow_eig}
\end{figure}

Figure~\ref{fig:flow_eig} plots the paths of the near-zero eigenvalues as functions of $k$ for $k\in (5,15)$; the spectra are computed 
every \(0.025\), and the arrows placed at these points.

\paragraph{Observations from these figures.} We make five observations from these figures. 
We number them O1, O2(a)-(d), corresponding, respectively, to Assumptions A1 and A2 below, under which we prove a bound on the $k$-dependence of the number of GMRES iterations (Theorem \ref{thm:main2}).
\ben
\item[O1] There exists a bounded open set $\cN$ ($\cN$ for ``near-zero'') containing zero and a closed half-plane $\cH$ not containing zero such that 
(i) all the eigenvalues of $\bfM^{-1}\matrixD$ are contained in $\cN\cup \cH$, and 
(ii) $\cN$ and $\cH$ can both be chosen to be independent of $k$.
\een
Point (i) is clear from Figure \ref{fig:example_spectrum} that 
plots the eigenvalues of $\bfM^{-1}\matrixD$ for the small and large cavities at $k=100$ and $k=290$. 

Point (ii) 
 is shown in Figure \ref{fig:outliers}; 
 indeed, all the near-zero eigenvalues for these two different values of $k$ lie on the same curve, and thus $\cN$ can be taken as the black rectangle in Figure \ref{fig:outliers}. In addition, the curve is the same for the four different $\Oi$ considered; this is perhaps expected since the four different $\Oe$ all contain a neighbourhood of the minor axis of the same ellipse (namely \eqref{eq:ellipse} with $a_1=1$ and $a_2=1/2$) and the near-zero eigenvalues of $\bfM^{-1}\matrixD$ are generated by the trapped ray in this neighbourhood.\footnote{However, with $a_1=1$ and $a_2=1/4$ (i.e., a different ellipse) 
 we see the eigenvalues lying on, by eye, the same curve as Figure \ref{fig:outliers} for analogous small and large cavities at $k=100$ and $k=290$.}
\ben
\item[O2(a)] (Family of quasimodes.) There exists a sequence $0<k_1<k_2<\ldots$, with $k_j\tendi$ as $j\tendi$, such that $\bfM^{-1}\matrixDj$ has a near-zero singular value for  $j$ sufficiently large.
\een
Recall that O2(a) is guaranteed on the continuous level by the lower bound \eqref{eq:Ainvblowup}, and we see small singular values (orange crosses) in three of the four plots in Figure \ref{fig:spectra_illustrations_elliptic_cavity_2D_dir} and all four plots in Figure \ref{fig:spectra_illustrations_elliptic_cavity_bigger_2D_dir}.
\ben
\item[O2(b)] (Quasimode $\implies$ near-zero eigenvalue.) If $j$ is sufficiently large, $\bfM^{-1}\matrixDj$ has a near-zero eigenvalue.
\een 
This can be seen from the fact that near-zero singular values are accompanied by near-zero eigenvalues
in both Figures \ref{fig:spectra_illustrations_elliptic_cavity_2D_dir} and  \ref{fig:spectra_illustrations_elliptic_cavity_bigger_2D_dir}.
\ben
\item[O2(c)] (Near-zero eigenvalues.) 
All the eigenvalues of $\bfM^{-1}\matrixD$ in the set $\cN$ in O1 move at a speed that can be bounded above and below by constants independent of $k$.
\een
In fact, Figure \ref{fig:flow_eig} indicates that the speed of the eigenvalues is independent of $k$ because
the arrows in the  box in Figure \ref{fig:flow_eig} are all evenly spaced; furthermore, 
we observe numerically that the speed is approximately one (at least for that range of $k$).
\ben
\item[O2(d)] The large cavity has more near-zero eigenvalues than the small cavity.
\een
This can be seen from Figure \ref{fig:example_spectrum}. In addition, comparing the top-left plots of Figures~\ref{fig:spectra_illustrations_elliptic_cavity_2D_dir} and~\ref{fig:spectra_illustrations_elliptic_cavity_bigger_2D_dir} we see that 
when $k=k^o_{0,3}$ there is no near-zero eigenvalue for the small cavity, but there is for the large cavity. 
This observation is the reason for the feature F3(b) (the number of iterations is larger for the large cavity than the small cavity).

Observation O2(d) can be partially explained from the fact that a larger number of the Laplace eigenfunctions of the ellipse $E$ \eqref{eq:ellipse} (from which the quasimodes in Theorem \ref{thm:ellipse} are constructed) are localised in the large cavity than in the small cavity. In the FEM case there is a close connection between the functions in the quasimodes and the eigenfunctions of the Galerkin matrix (see \cite[Remark 1.7]{GaMaSp:21}), and thus these localisation considerations immediately explain why the large cavity has more near-zero eigenvalues than the small cavity. However, in the BEM case it is less clear how the eigenvalues of $\bfM^{-1}\matrixD$ (which are discretisations of functions living on $\Gamma$) are connected to the functions in the quasimodes (which live in $\Oe$).

We return to Observation O2(d) in \S\ref{sec:Weyl_new}
where we use heuristics from Weyl asymptotics to estimate how many more Laplace eigenfunctions of the ellipse are 
localised in $\Oe$ 
for the large cavity than for the small cavity. We then compare these heuristics to the number of eigenvalues of $\bfM^{-1}\matrixD$ observed computationally (see \S\ref{sec:numA2}).

\paragraph{Link with other experiments/results in the literature.}

Similar eigenvalue plots for BIEs 
when $\Oi$ is nontrapping or weakly trapping can be found for $A_k'$ in \cite[Figure 9]{Ba:06} and \cite[Figure 3.1 and \S5]{BeSp:11}, and for the indirect analogue of $B_{k, {\rm reg}}$ in \cite[Figure 1]{BrElTu:12}, \cite[Figures 3-5]{BoTu:13}, and  \cite[Figure 2]{ViGrGi:14}.
Furthermore, the analogue of O2(c) (the eigenvalues move at $O(1)$ speed) 
was used to compute large eigenvalues of the Laplacian in \cite{TuSc:07} for BIEs related to $A_k'$  
and \cite{BaHa:14} for a related boundary-based method; see \cite[Figure 5]{TuSc:07} and \cite[Theorem 4.1]{BaHa:14}, respectively.

\subsection{First main result: general bound on number of GMRES iterations for matrix with ``cluster plus outlier'' structure.}
\label{sec:theorem}

For simplicity we consider matrices with simple eigenvalues; the modifications to our assumptions and arguments for matrices with repeated eigenvalues are outlined in Remark \ref{rem:BF1}.

For $\lambda$ a simple eigenvalue,  let $\kappa (\lambda)$ 
be the condition number of $\lambda$ defined by
\begin{align}\label{eq:eigenvaluecondition}
    \kappa (\lambda):= \dfrac{\lVert \bfu \rVert \lVert \bfv \rVert}{\lvert \bfu \cdot \bfv \rvert },
\end{align}
where \(\bfu,\bfv\in\bbC^n\) are the right and left eigenvectors, respectively, corresponding to $\lambda$.
Recall that \(\kappa (\lambda)\geq 1\), and equality holds when \(\bfu\) and \(\bfv\) are collinear (which is guaranteed if the matrix 
 is normal).

\begin{theorem}[Bound on the GMRES relative residual]\label{thm:main1}
    Let \(\genmatrix\in \bbC^{n\times n}\) 
    be diagonalisable with simple eigenvalues. 
     Assume that all the eigenvalues of $\genmatrix$ are contained in $\cN \cup \cH$, where $\cN$ is a bounded open set containing zero and $\cH$ a closed half plane not containing zero. Without loss of generality, let $\cH:=\{\Re z \geq S\}$ for some $S>0$ and assume that $\cH$ contains at least one eigenvalue of $\genmatrix$, so that $S\leq \|\genmatrix\|_2$. 
Let  $\lambda_1,\ldots \lambda_\ell$ be the eigenvalues in $\cN$, and let 
 $\kappa^*$ be the maximum eigenvalue condition number of $\genmatrix$.    
    
Given $L_0, L_1$ with $0<L_0<L_1\leq S$, let 
\beqs
N_{\rm eig}:= \big|\{j \, :\,  \lambda_j \in \{ L_0 < \Re z < L_1 \} \cap\cN\} \big| + 1.
\eeqs
    Let 
 \beq\label{eq:delta}
 \delta:= \frac{L_1-L_0}{4 n\,\kappa^* \, N_{\rm eig}}.
 \eeq
     Let $\beta\in (0,\pi/2)$ be defined by 
     \beq
     \cos \beta =\frac{L_0}{\lVert \genmatrix\rVert_2+\delta}, \quad\text{ and let }\quad
    \label{eq:gammabeta}
    \gamma_\beta := 2\sin \left( \dfrac{\beta}{4-2\beta/\pi}\right)<1.
    \eeq
Then, when GMRES is applied to the equation $\genmatrix \bfx=\bfb$, the $m$th GMRES residual \eqref{eq:residual} satisfies
    \begin{align}\label{eq:Tues1new}
        \dfrac{\lVert \bfr_m (\genmatrix,\bfb,\bfx_0) \rVert_2}{\lVert \bfr_0 (\genmatrix,\bfb,\bfx_0) \rVert_2} 
        \leq \left(\prod_{j=1}^\ell  \dfrac{1}{\lvert \lambda_j \rvert}\right) \,(\|\genmatrix\|_2 + \delta)^{\ell+1}\,3^{\ell+1} \,\delta^{-1} \gamma_{\beta}^{m-\ell}.
            \end{align}
\end{theorem}

\begin{corollary}[Sufficient condition on the number of GMRES iterations for convergence]\label{cor:SCW}
Under the assumptions of Theorem \ref{thm:main1}, given $\eps>0$, if
\begin{align}
 m\geq      \ell+ \left(\log(\gamma_{\beta}^{-1})\right)^{-1} \left(
 \sum_{j=1}^\ell\log\dfrac{1}{\lvert \lambda_j \rvert}+ 
 \log (\varepsilon^{-1})+\log \left(\delta^{-1}\right)+
(\ell+1)\Big(\log\big(\lVert \genmatrix\rVert_2+\delta\big)+ \log 3\Big)\right),
 \label{eq:mlowerbound}
    \end{align}
then, when GMRES is applied to the equation $\genmatrix \bfx=\bfb$, the $m$th GMRES residual \eqref{eq:residual} satisfies
\beqs
\frac{\lVert \bfr_m (\genmatrix,\bfb,\bfx_0) \rVert_2 }{ \lVert \bfr_0 (\genmatrix,\bfb,\bfx_0) \rVert_2}\leq \varepsilon.
\eeqs
    \end{corollary}  

\begin{figure}[h!]
    \centering
    \includegraphics[width=0.35\textwidth]{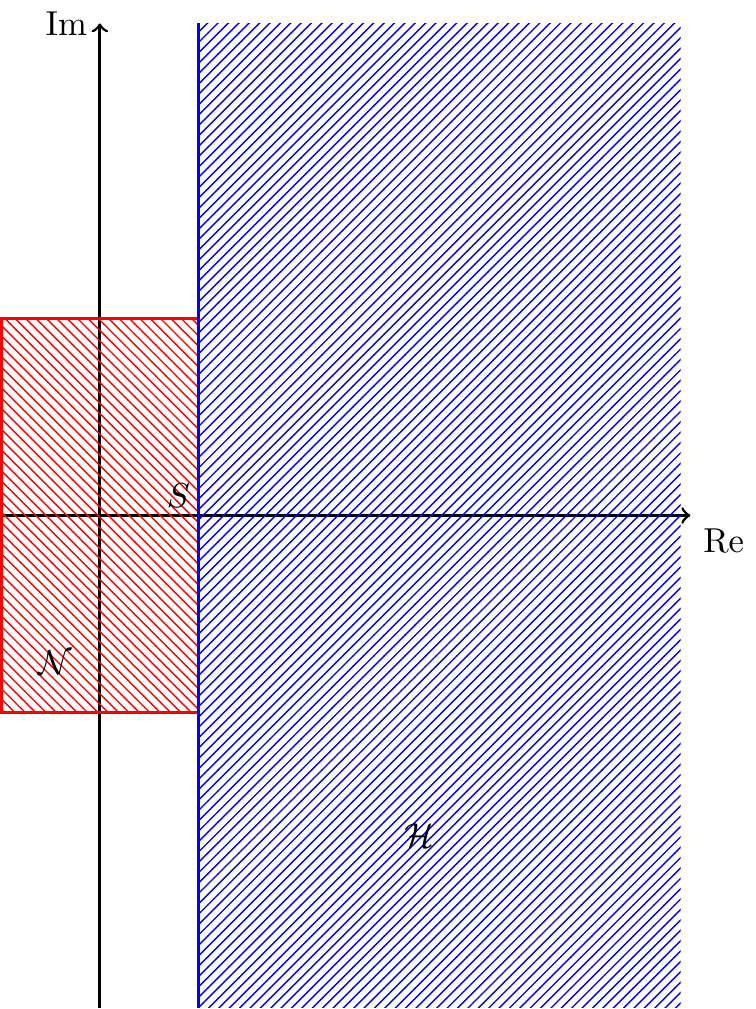}
    \caption{The sets $\cN$ and $\cH$ in the assumptions of Theorem \ref{thm:main1}.}\label{fig:A3A4}
\end{figure}

Figure \ref{fig:A3A4} shows the half-plane $\cH$ in Theorem \ref{thm:main1}, and an example of a possible $\cN$ in Theorem \ref{thm:main1}. 

\bre[The dependence of $\gamma_\beta$ on $\|\genmatrix\|_2$.]
How $\gamma_\beta$ depends on $\|\genmatrix\|_2$ is not immediately clear from the definitions in \eqref{eq:gammabeta}.
However, if $\|\genmatrix\|_2 \gg \delta$ and $\|\genmatrix\|_2 \gg L_0$, then 
\beq\label{eq:Wed1}
\left(\log(\gamma_{\beta}^{-1})\right)^{-1} = \frac{3\sqrt{3}}{4} \left(\frac{\|\genmatrix\|_2}{L_0} \right)\left( 1 + O\left(\frac{\delta}{\|\genmatrix\|_2}\right)\right)
 \left( 1 + O\left(\frac{L_0}{\|\genmatrix\|_2}\right)\right).
\eeq
Indeed, let $\alpha := \pi/2 -\beta$ so that $\cos \beta = \sin \alpha$
and, from the definition of $\beta$,  
\beq\label{eq:Wed2}
\alpha = \frac{L_0}{\|\genmatrix\|_2} \left( 1 + O\left(\frac{\delta}{\|\genmatrix\|_2}\right)\right).
\eeq
The definition of $\gamma_\beta$ then implies that, as $\alpha\tendo$,
\beq\label{eq:Wed3}
\gamma_\beta = 1 - \frac{4\alpha}{3\sqrt{3}} + O(\alpha^2) \quad\text{ so that } \quad -\log \gamma_\beta = \frac{4\alpha}{3\sqrt{3}} + O(\alpha^2).
\eeq
The asymptotics \eqref{eq:Wed1} then follow from combining \eqref{eq:Wed2} and \eqref{eq:Wed3}.
\ere

\paragraph{Interpreting the bound \eqref{eq:mlowerbound}.}
The bound \eqref{eq:mlowerbound} shows that each outlier $\lambda_j$ contributes $C_1 + C_2 \log(1/|\lambda_j|)$ to the number of iterations needed to guarantee a prescribed residual reduction, where $C_1$ and $C_2$ are independent of $\lambda_j$ but depend on $\|\genmatrix\|_2$.
Therefore, if each $|\lambda_j|$ is large, only the number of outliers contributes to the required number of iterations.
If $|\lambda_j|$ is small, its value can have more of an effect on the required number of iterations, but this effect is mitigated by the fact that $|\lambda_j|$ appears in a logarithm.

\paragraph{The ideas behind the proof of Theorem \ref{thm:main1}.}

A convergence theory for GMRES based on modelling the eigenvalues as a ``cluster plus outliers'' was famously used in  \cite{CaIpKeMe:96}, with the idea arising in the context of the conjugate-gradient method \cite{Je:77} and used subsequently, e.g., in \cite{ElSiWa:02}.
This theory in \cite{CaIpKeMe:96} forms the starting point for proving Theorem \ref{thm:main1}; see Lemma \ref{lem:Campbell} below.
 However, a crucial difference is that we are interested in matrices depending on a parameter, namely $k$.
We therefore augment the theory in \cite{CaIpKeMe:96}, with 
(i) the results in \cite{BeGoTy:05} about polynomial min-max problems, and 
(ii) results about pseudospectra appearing in, e.g., \cite{TrEm:05}. 

The result is that when the bound on the number of GMRES iterations \eqref{eq:mlowerbound} is applied with $\genmatrix = \bfM^{-1}\matrixD$, the $k$-dependence of the quantities in the bound (i.e.~$S, L_0, L_1, N_{\rm eig}, \delta, \beta$) is given from \emph{either} 
existing $k$-explicit bounds on the norm and the norm of the inverse of $A_k'$ \emph{or} assumptions about the $k$-dependence of both the number and the condition numbers of the eigenvalues (see Assumptions A2 and A3 below).

We highlight that, in our use of the pseudospectrum, we choose $\delta$ as a function of  $k$ to compensate for the growth of the non-normality with $k$.
This flexibility in choosing $\delta$ is mentioned in~\cite[Page 6]{Em:99} 
when analysing different stages of the GMRES iteration for a single linear system; in contrast, here we use this flexibility applied to a family of linear systems parametrised by $k$.

\subsection{Second main result: $k$-explicit bound on the number of GMRES iterations 
for  Helmholtz BIEs under strong trapping}
\label{sec:mainresultH}

\subsubsection{Statement of assumptions.}

We write $A\lesssim B$ if there exists $C>0$, independent of all parameters of interest (including $h$ and $k$), such that $A\leq C B$, and $A\sim B$ if both $A\lesssim B$ and $A\gtrsim B$.

\ben 
\item[A0] The meshwidth $h$ is chosen as a function of $k$ so that, given \(\Capproxone>1\), \(\Capproxtwo>1\), for all $k$, 

(i) the Galerkin solution exists, is unique, and the relative $L^2(\Gamma)$-error of the Galerkin solution is bounded uniformly in $k$, 

(ii) the second inequality in \eqref{eq:discon2} holds, i.e., $\N{\bfM^{-1}\matrixD}_{2} \leq \Capproxone\N{A_k'}_{\LtGt}$, and

(iii) the eigenvalues of both $A_k'$ and $\bfM^{-1}\matrixD$ are simple and
the eigenvalues of $A_k'$ are approximated by the eigenvalues 
of $\bfM^{-1}\matrixD$ in the following sense:
at a given $k$, let the eigenvalues of $\bfM^{-1}\matrixD$ be $\lambda_1,\ldots,\lambda_n$ (where $n=n(k)$); there exists an injective function $f_k : \{1,\ldots,n\}\to \{\mu : \mu \text{ an eigenvalue of } A_k'\}$ such that 
\beq\label{eq:A0iii}
|\lambda_j|\geq \Capproxtwo |f_k(j)|
\eeq
\item[A1] 
There exists a bounded open set $\cN$ containing zero and a closed half-plane $\cH$ not containing zero such that 
(i) all the eigenvalues of $\bfM^{-1}\matrixD$ are contained in $\cN\cup \cH$, and 
(ii) $\cN$ and $\cH$ are both independent of $k$.
\item[A2] The number of eigenvalues of $\bfM^{-1}\matrixD$ in the set $\cN$ in A1 is $\leq \CWeyl k^{d-1}$, where $\CWeyl>0$ is independent of $k$.
\item[A3] With $\kappa^*(k)$ the maximum eigenvalue condition number of $\bfM^{-1}\matrixD$, there exists $\Ccond, M>0$ and independent of $k$ such that $\kappa^*(k)\leq \Ccond k^M$.
\een

\paragraph{Why do we expect Assumptions A0-A3 to hold?}

As recapped in \S\ref{sec:statement:BEM}, there exist results on which functions $h=h(k)$ ensure A0(i), although they do not appear to be sharp.
A0(ii) and A0(iii) are ensured at least as $h\tendo$ for fixed $k$ when $\Gamma$ is $C^1$, with this regularity of $\Gamma$ ensuring that $A_k'$ (and also $B_{k, {\rm reg}}$) is a multiple of the identity plus a compact operator on $\LtG$ (see Remark \ref{rem:compact}).
Indeed, in this case A0(ii) holds by Lemma \ref{lem:discon} below and stronger results than A0(iii) (showing that the eigenvalues of $\bfM^{-1}\matrixD$ converge to those of $A_k'$, with multiplicity) hold by \cite[Theorems 2, 3]{At:67}, \cite[Theorem Page 214]{At:75} (see also \cite[Theorem 7]{Sp:75}, \cite[Theorem 4.1]{SpTh:83}).

Note that \eqref{eq:A0iii} specifies that, with $\mu(k)$ an eigenvalue of $A_k'$ and $\lambda(k)$ an eigenvalue of $\bfM^{-1}\matrixD$, $|\lambda(k)|\gtrsim |\mu(k)|$. We do not require that $|\lambda(k)|\sim |\mu(k)|$ because we expect at least one $|\mu(k_j)|$ to be exponentially small when $k_j$ is a quasimode frequency (this is proved for the standard variational formulation, i.e., the basis of FEM, in \cite[Theorem 1.5]{GaMaSp:21}), but we expect that $|\lambda(k_j)|$ will be only algebraically small because of the sensitivity in F2. Note also that A0(iii) assumes that this sensitivity does not cause $|\lambda(k)|\ll |\mu(k)|$ for $k$ near $k_j$.

Regarding A1: first note that this corresponds to Observation O1 in \S\ref{sec:observations}. 
When $\Oe$ is nontrapping $\|(A_k')^{-1}\|_{\LtGt}\lesssim 1$ \cite[Theorem 1.13]{BaSpWu:16} and thus the smallest singular value of $\bfM^{-1}\matrixD\sim 1$ in this case.  
This implies that when $\Oi$ is nontrapping the eigenvalues of $\bfM^{-1}\matrixD$ are $\sim 1$ away from zero.
Furthermore, at least for some nontrapping $\Oi$, the eigenvalues are contained in a $k$-independent half-plane away from zero since $A_k'$ is coercive (with constant independent of $k$) \cite{SpKaSm:15}, \cite{BeSp:11}. These facts suggest that the second part of  A1 holds (i.e.~the half-plane $\cH$ is independent of $k$), but are far from a proof.

Regarding Assumption A2: in \S\ref{sec:Weyl_new} we give heuristic arguments backing up this assumption, one based on Weyl-type asymptotics for 
eigenvalues of the Laplacian on bounded domains, and the other based on the Observations O2(a)-(c) in \S\ref{sec:observations} and results about the number of resonances of the exterior Helmholtz problem.

Regarding Assumption A3: we did not give any experiments about this assumption in \S\ref{sec:intro_num}, but these are contained in \S\ref{sec:experiments}, and indicate that A3 holds for BEM discretisations of each of the BIEs \eqref{eq:direct_combined_dir},
\eqref{eq:direct_combined_neu}, and \eqref{eq:direct_combined_neu2}.

\subsubsection{$k$-explicit bounds on the number of GMRES iterations via Theorem \ref{thm:main1}}
\label{sec:Helmholtz}

Before applying Theorem \ref{thm:main1} to $\bfM^{-1}\matrixD$, we recall that GMRES applied to an $n\times n$ matrix converges in at most $n$ iterations (in exact arithmetic). This bound is well known to have ``little practical content'' \cite[Page 270]{TrBa:97} since one never reaches this number of iterations; nevertheless, it does give a theoretical upper bound on the $k$-dependence of the number of iterations. For example, when the $h$-BEM uses a fixed number of points per wavelength, $n\sim k^{d-1}$, and thus there exists $C>0$ such that if $m\geq C k^{d-1}$ then GMRES converges. However, the constant $C$ is 
 both large and dependent on the number of points per wavelength, and this is not what one sees in practice. For example, for the large cavity in 2-d with $k=100$, \(\theta = 4\pi/10\), and ten points per wavelength, $n= 1766$ and GMRES converges to tolerance \(10^{-6}\) in 165 iterations. For twenty points per wavelength, $n= 3528$ and GMRES converges to the same tolerance in 167 iterations.
 
\begin{theorem}[Bound for Helmholtz BIEs]
\label{thm:main2}
Let $\Oi$ be piecewise smooth. Assume that $n\leq C_{\rm{dis}}k^{M'}$ for some $M'>0$ and $C_{\rm{dis}}>0$.  
Consider GMRES  applied to the linear system 
\beqs
\bfM^{-1}\matrixD \bu = \bfM^{-1}\bff.
\eeqs
where $\matrixD$ \eqref{eq:Galerkinmatrix} is the Galerkin matrix from the BEM discretisation of 
the Dirichlet BIE \eqref{eq:direct_combined_dir} and $\bfM$  \eqref{eq:massmatrix} is the mass matrix. 
If Assumptions A0-A3 hold, and $\lambda_1,\ldots \lambda_\ell$ are the eigenvalues of $\bfM^{-1}\matrixD$ in $\cN$,
then there exists $C_j>0$, $j=1,2$ (independent of $k$) such that given $\eps>0$, for all $k\geq 1$,
if 
\beq\label{eq:mlowerbound2}
 m\geq     \CWeyl k^{d-1} + C_1 \N{A_k'}_{\LtGt}\Bigg(  
 \sum_{j=1}^\ell\log\dfrac{1}{\lvert \lambda_j \rvert}
 + C_2 k^{d-1}\log k
+ \log (\varepsilon^{-1})
 \Bigg),
\eeq
then the $m$th GMRES residual satisfies $\|\bfr_m \|_2/\|\bfr_0 \|_2 \leq \eps$. 
Furthermore, $C_1$ only depends on $\Capproxone$, $\Omega_-$, and $d$, and $C_2$ only depends on $\Capproxone, \CWeyl, \Ccond, C_{\rm{dis}},d, M,$ $M'$, and $\Omega_-$ with these constants as defined in Assumptions A0-A3.
\end{theorem} 

An analogous result holds with $A_k'$ replaced by $B_{k,{\rm reg}}$, if $B_{k,{\rm reg}}$ satisfies appropriate analogues of Assumptions A0-A3; see Remark \ref{rem:Neumann} below.

The bound~\eqref{eq:mlowerbound2} gives insight into how the \(k\)-dependence of the number of iterations arises from the eigenvalue distribution of $\bfM^{-1}\matrixD$. Moreover, with $M'$ and $C_{\rm dis}$ fixed, the constants in \eqref{eq:mlowerbound2} are independent of the choice of $n\leq C_{\rm{dis}} k^{M'}$ (provided that Assumptions A0-A3 hold). Therefore, choosing $M'>d-1$ and $C_{\rm dis}>0$, discretizations satisfying $n\leq C_{\rm{dis}} k^{M'}$ include those with an arbitrary number of points per wavelength, and the bound \eqref{eq:mlowerbound2} holds, at least for sufficiently large $k$, uniformly across all of them (which appears consistent with the specific examples of numbers of iterations stated above the theorem).  However, the right-hand side of \eqref{eq:mlowerbound2} contains terms that grow faster than $k^{d-1}$ and so, if $n\sim k^{d-1}$, then the $k$-dependence of \eqref{eq:mlowerbound2} is worse than that of the crude bound that GMRES converges in at most $n$ iterations (in exact arithmetic).
  
\paragraph{Informal explanation of how \eqref{eq:mlowerbound2} arises from \eqref{eq:mlowerbound}.}
When $\genmatrix =\bfM^{-1}\matrixD$, Assumption A1 implies that $S\sim 1$, and thus we can choose $L_0, L_1\sim 1$. Then $\delta$ defined by \eqref{eq:delta} is $\lesssim  1$. By Assumption A0, $\N{\bfM^{-1}\matrixD}_2 \lesssim \N{A_k'}_{\LtGt}$, which grows at most algebraically with $k$ (see Lemma \ref{lem:normbounds} below). Using this, along with \eqref{eq:Wed1}, we find that the bound on $m$ \eqref{eq:mlowerbound} holds if 
\begin{align}
 m -\ell \gtrsim \N{A_k'}_{\LtGt}\left(
 \sum_{j=1}^\ell\log\dfrac{1}{\lvert \lambda_j \rvert}+ 
 \log (\varepsilon^{-1})+\log \left(\delta^{-1}\right)+ \ell \,C \log k
 \right),
 \label{eq:mlowerbound_explain}
    \end{align}
    for some $C>0$ independent of $k$.
 
Assumption A2 is that    $\ell \leq\CWeyl k^{d-1}$, and, by definition, $N_{\rm eig}\leq \ell +1$. Assumption A3 and the bound $N_{\rm eig}\leq \CWeyl k^{d-1} +1$ then imply that $\delta^{-1}$ grows at most polynomially in $k$; the bound \eqref{eq:mlowerbound2} then follows from using these bounds in \eqref{eq:mlowerbound_explain}.
    
\subsubsection{Discussion of the $k$-dependence of the bound in Theorem \ref{thm:main2}, how this bound explains F3(a), and how this bound could be improved.}\label{sec:discussion}

To investigate the sharpness of the bound \eqref{eq:mlowerbound2} in Theorem \ref{thm:main2}, we summarise the results of the numerical experiments from \S\ref{sec:numA2} and \S\ref{sec:experiments} in Table \ref{tab:kdependence}. This table
plots the $k$-dependence for $k\in (50,290)$, for both the small and large cavities, of (i) the number of iterations, (ii) the number of outlier eigenvalues $\ell$ when $\cN:= [-0.1,0.1]\times[-0.6,0.6]$ (i.e., the black rectangle in Figure \ref{fig:outliers}), (iii) $\|\bfM^{-1}\matrixD\|_2$ (and its analogue for the two Neumann BIEs), and (iv) the quantity 
\beq\label{eq:cL}
\cL:= \sum_{j=1}^\ell\log\dfrac{1}{\lvert \lambda_j \rvert}.
\eeq
Each of $\|\bfM^{-1}\matrixD\|_2$ and its Neumann analogues has the same $k$-dependence for both the small and large cavities -- see the top-left plots in Figures \ref{fig:Dirichlet}, \ref{fig:Neumann}, and \ref{fig:reg_Neumann} below -- and so the norm only appears in one column in Table \ref{tab:kdependence}.
The exponents in Table \ref{tab:kdependence} are obtained using 
the nonlinear least-squares Marquardt-Levenberg algorithm (the basis of the `fit' command in gnuplot). 

\begin{table}[h!]
\begin{tabular}{c|| c| c| c| c|c|c|c|c|c|}
   & 
  $\#$ it.~small & $\#$ it.~large 
   & $\ell$ small & $\ell$ large & norm & $\cL$ small 
& $\cL$ large \\
\hline\hline
Dirichlet $A_k'$ & 0.66 & 0.82 & 0.95 & 1.00 & 0.31 & 0.77 & 0.90 \\
Neumann $B_k$ & 0.55 & 0.77 & 0.93 & 0.97 & 0 & 0.78 &0.91 \\
reg.~Neumann $B_{k, {\rm reg}}$ & 0.60 & 0.80 & 0.95 & 0.95 & 0 &0.79 & 0.89
\end{tabular}
\caption{The exponents in how the quantities in the columns vary with $k$ in 2-d (i.e.~if a quantity $\sim k^a$, then $a$ is displayed), as determined by numerical experiments through $k=k^e_{m,0}$ for $k\in (50,290)$.}\label{tab:kdependence}
\end{table}

We structure our discussion around the following points. 

\paragraph{The number of iterations grows slightly less than $k^{d-1}$ for the large cavity.} Table \ref{tab:kdependence} shows that, in 2-d, the number of iterations roughly $\sim k^{0.6}$ for the small cavity and $\sim k^{0.8}$ for the large cavity, for each of the three BIEs. Figure \ref{fig:3d} below shows that for the Dirichlet problem in 3-d the number of iterations 
grows roughly like $k^2$ for both the small and large cavities
 over the range $k\in (20,110.5)$ -- note that this is a smaller range than we consider in 2-d. 
Similarly, Figure \ref{fig:3d_neu} below shows that for the Neumann problem in 3-d the number of iterations for both $B_k$ and $B_{k, {\rm reg}}$ grows roughly like $k^2$ for the small cavity over the range $k\in (20,50)$.
  In both 2- and 3-d, the number of iterations therefore grows with $k$ at roughly the same rate as the number of degrees of freedom, illustrating how difficult a problem this is.
 
\paragraph{The bound \eqref{eq:mlowerbound2} will always give $m\gtrsim k^{d-1}$ because $m\gtrsim \ell$ and $\ell \sim k^{d-1}$.}
Assumption A2 is that $\ell \lesssim k^{d-1}$; in \S\ref{sec:Weyl_new} we present heuristic arguments based on Weyl asymptotics why A2 holds, and the numerical experiments for $d=2$ indicate that $\ell \sim k$. Therefore, since the right-hand side of \eqref{eq:mlowerbound} contains $\ell$, the bound \eqref{eq:mlowerbound2} gives that $m\gtrsim k^{d-1}$ no matter what bounds we obtain on the other quantities in \eqref{eq:mlowerbound}.

\paragraph{The near-zero eigenvalues are, at worst, exponentially small.}
By Assumption A3(iii), $|\lambda_j|^{-1}\lesssim |\mu|^{-1}$, where $\mu:=f_k(\lambda_j)$ is an eigenvalue of $A_k'$.
Since 
$|\mu|^{-1}\leq \|(A_k')^{-1}\|_{\LtGt}$
 and, at least for smooth $\Oi$, $\|(A_k')^{-1}\|_{\LtGt}\lesssim \exp(\alpha k)$ by \cite[Equation 1.35 and Lemma 6.2]{ChSpGiSm:20}, 
$|\mu|^{-1} \lesssim \exp(\alpha k)$. 

\paragraph{The bound \eqref{eq:mlowerbound2} explains F3(a)
because it shows that the actual position of each near-zero eigenvalue is less important than the total number of near-zero eigenvalues.}
Indeed Figures \ref{fig:example_spectrum}, \ref{fig:spectra_illustrations_elliptic_cavity_2D_dir}, and \ref{fig:spectra_illustrations_elliptic_cavity_bigger_2D_dir} indicate that the near-zero eigenvalues are not all simultaneously close to zero. 
Since each $\lambda_j$ enters the bound on $m$ \eqref{eq:mlowerbound2} via the $\log(1/|\lambda_j|)$ in $\cL$ \eqref{eq:cL}, even if one of the $\lambda_j$ is exponentially small (which we expect to be very unlikely by F2), we expect the growth of $\cL$ to still be dominated by the overall number of near-zero eigenvalues, i.e., $\ell \sim k^{d-1}$. This is consistent with the fact that $\ell \sim k$ and  $\cL\sim k^{0.9}$ in the experiments for the large cavity.

Making this argument rigorous and proving that $\cL\lesssim k^{d-1}$ would involve first bounding
\beq\label{eq:boundcL}
\cL:=\sum_{j=1}^\ell \log\frac{1}{|\lambda_j|}\lesssim 
\sum_{\substack{\mu \in \sigma(A_k'),\, \mu \in \widetilde{\cN}}}
 \log\frac{1}{|\mu|},
\eeq
where $\widetilde{\cN}$ is a neighbourhood of $\cN$ (depending on $\Capproxtwo$) containing the images of each $\lambda_j\in \cN$ under $f_k$, and then controlling the number of eigenvalues that can simultaneously be exponentially-close to zero. To our knowledge, the question of whether there exists strong trapping with high multiplicity of quasimodes/resonances is still open.
Even if we knew that all near-zero eigenvalues correspond to localised eigenfunctions of the ellipse, controlling this number
involves understanding the number of eigenvalues exponentially-close together. This could be obtained from proving a Weyl law with $k^{d-2}$ remainder, but this has only been established so far for the torus for $d\geq 5$ \cite{Fr:82}, and is known to be typically false (e.g., on the torus for $d=2$ \cite{Ha:15}). 

\paragraph{The norms of the operators grow with $k$ in the limit $k\tendi$, and so, for $k$ sufficiently large and assuming $\cL\lesssim k^{d-1}$, the bound \eqref{eq:mlowerbound2} will be dominated by $\|A_k'\|_{\LtGt} k^{d-1}\log k$.}
How $\|A_k'\|_{\LtGt}$ depends on both $k$ and the geometry of $\Oi$ is now well-understood thanks to \cite{ChGrLaLi:09}, \cite[Appendix A]{HaTa:15}, \cite{GaSm:15},  \cite[Chapter 4]{Ga:19}, and \cite{GaSp:19}. These results show that, as $k\tendi$, $k^{1/2}\lesssim \|A_k'\|_{\LtGt}\lesssim k^{1/2}\log k$ for both the small and large cavities, and we expect that the ideas behind these results can be used to show that 
$\|B_{k, {\rm reg}}\|_{\LtGt}\gtrsim k^{1/6}$ for these domains;
 see the discussion in \S\ref{sec:numbounds} about the top-left plots.
 Table \ref{tab:kdependence} shows, however, that in the range $k \in(50,290)$ the growth of $\|A_k'\|_{\LtGt}$ is slower than $k^{1/2}$, and $\|B_{k, {\rm reg}}\|_{\LtGt}$ does not grow. For $\|A_k'\|_{\LtGt}$, this discrepancy is explained in \S\ref{sec:numbounds}.

\paragraph{Limitations of the ``cluster plus outliers'' model applied to $\bfM^{-1}\matrixD$ and how it could be improved.}

The limitations of the ``cluster plus outliers'' model where the ``outliers'' are the near-zero eigenvalues are shown in 
Figures \ref{fig:example_spectrum}, \ref{fig:spectra_illustrations_elliptic_cavity_2D_dir}, and \ref{fig:spectra_illustrations_elliptic_cavity_bigger_2D_dir}. Indeed, these plots show that the ``cluster'' of eigenvalues away from zero is itself a cluster with outliers. Furthermore, this ``cluster within the cluster'' appears to be contained in a $k$-independent set (see Figure \ref{fig:example_spectrum}).

We therefore expect that a bound with improved $k$-dependence could be obtained by 
taking this additional structure into account. Indeed, $\|A_k'\|_{\LtGt}$ currently enters the bound \eqref{eq:mlowerbound2} and as a bound on the modulus of the cluster eigenvalues. 
If one could prove that the number of the ``outliers of the cluster'' is $\lesssim k^{d-1}$ and that the ``cluster within the cluster'' is contained in a $k$-independent set, then this would replace the bound \eqref{eq:mlowerbound2} with a bound of the form $m \gtrsim k^{d-1} + \cL$. If one could, in addition, prove that $\cL\lesssim k^{d-1}$ (as discussed above), then this would prove the sharp bound $m\gtrsim k^{d-1}$ (with the omitted constant independent of $k$ and properties of the discretisation).

\subsection{Applicability of the ideas in this paper to Helmholtz FEM.}\label{sec:FEM}

Until now we have focused on solving the scattering problem \eqref{eq:Helmholtz} using BIEs and BEM, however 
the general result of Theorem \ref{thm:main1} can be applied to other Helmholtz discretisations satisfying Assumptions A0-A3
(or suitably modified versions of these).

\paragraph{Location and number of the near-zero eigenvalues for FEM.} As mentioned in \S\ref{sec:observations}, the connection between quasimodes and near-zero eigenvalues of the standard domain-based variational formulation (i.e.~the basis of FEM) is much clearer than for BEM, and this is subject of the companion paper \cite{GaMaSp:21}. Indeed, \cite[Theorem 1.5]{GaMaSp:21} proves that if $k=k_j$, then there exists a near-zero eigenvalue
of the standard domain-based variational formulation, with the distance of this eigenvalue from zero given in terms of the quality $\QMC(k_j)$ of the quasimode. Furthermore, 
\cite[Theorem 1.8]{GaMaSp:21} shows that the eigenvalues inherit the multiplicities of the quasimodes. These results are proved using arguments from microlocal and complex analysis, inspired by the celebrated ``quasimodes to resonances'' results of 
 \cite{TaZw:98}, \cite{St:99} (following \cite{StVo:95, StVo:96}); see also \cite[Theorem 7.6]{DyZw:19}. 

We highlight that while the number of near-zero eigenvalues in the FEM case is the same as for BEM (namely $\sim k^{d-1}$) the number of degrees of freedom for FEM is much larger than that for BEM. Indeed, while BEM is commonly used with a fixed number of points per wavelength, leading to systems of size $\sim k^{d-1}$, the pollution effect means that FEM is used with systems of size $\gg k^d$. 
Therefore, the issue we encountered in \S\ref{sec:Helmholtz} for BEM that the number of near-zero eigenvalues is the same order as total number of degrees of freedom does not occur for FEM.

\paragraph{Location of the other eigenvalues for FEM.} 
While the near-zero eigenvalues for FEM are easier to understand rigorously than those for BEM, one subtlety in the FEM case is that the eigenvalues away from zero need not be in a half-plane (as in A1). If \emph{either} the exact Dirichlet-to-Neumann map \emph{or} an impedance boundary condition is used on the truncation boundary, then the numerical range (and hence the eigenvalues) is contained in the lower-half plane. Furthermore, if the problem is nontrapping, then the eigenvalues are contained in a half-plane an $O(1)$ distance below the real axis \cite[Theorem 1.12]{BaSpWu:16}, suggesting that A1 holds with $\cN$ and $\cH$ independent of $k$ (see also \cite[Theorem 5.1]{LiXiSaDe:20} for stronger results on the eigenvalue distribution of a simple nontrapping problem).

However, if a perfectly-matched layer (PML) is used, then the numerical range of the operator contains elements in the upper-half plane, as can be seen from \cite[Equation after (2.12)]{LiWu:19}.
Nevertheless, we expect that a similar result to Theorem \ref{thm:main1} can be proved under a modified version of A1 by 
replacing the domain $K_\beta\subset \mathbb{C}$ in Lemma \ref{lem:bound_polynomial_on_lens} below by a non-convex domain, such as one of the class introduced in \cite{KoLi:00}; 
see, e.g., the discussion in \cite[\S3.1.2]{LiTi:04}.

\paragraph{Preconditioning FEM discretisations.} The reason 
we have focused on BEM (and not FEM) in this paper is that 
one usually seeks to precondition GMRES applied to the FEM discretisation of the standard variational formulation of the Helmholtz equation. This is because the number of GMRES iterations without preconditioning grows rapidly with $k$ even in non-trapping scenarios. This is in contrast to BEM, where the number of GMRES iterations for discretisations of BIEs \eqref{eq:direct_combined_dir} and \eqref{eq:direct_combined_neu2} enjoy mild growth with $k$ in nontrapping situations; see \cite[Theorem 1.16 and Figure 1]{GaMuSp:19} for \eqref{eq:direct_combined_dir}  and \cite[Tables 1 and 2]{BoTu:13} for \eqref{eq:direct_combined_neu2}.

The design of good preconditioners for GMRES applied to the Helmholtz FEM in nontrapping scenarios is a very active area of research; see the literature reviews in \cite{Er:08, ErGa:12, GaZh:19}, and \cite[\S1.3]{GrSpZo:20}.
Since our theory below only proves bounds for GMRES applied to unpreconditioned matrices, 
our results are less interesting for FEM than for BEM.
Nevertheless, our results still provide insight into the design of preconditioners for trapping problems -- this is discussed in the conclusions \S\ref{sec:conclusions}.

\section{Proofs of Theorems \ref{thm:main1} and \ref{thm:main2}}

\subsection{Definition of spectral projectors}\label{sec:spectral_proj}

Given $\genmatrix \in \mathbb{C}^{n\times n}$, 
let $\lambda_1, \ldots, \lambda_\ell$ be a subset of the eigenvalues of $\genmatrix$ (we later choose this subset to be the eigenvalues in $\cN$ for a matrix satisfying the assumptions of Theorem \ref{thm:main1}, but the results in this subsection hold more generally).
  Let $\Gamma_j$, $j=1,\ldots,\ell$, be a circle enclosing $\lambda_j$ but no other eigenvalue of $\genmatrix$, and let $\Gamma = \cup_{j=1}^\ell \Gamma_j$. Let $\widetilde{\Gamma}$ be a positively-oriented curve enclosing the rest of the spectrum.
Let \(\bfR_k(z):=\left(z \bfI - \genmatrix\right)^{-1}\), i.e.~$\bfR_k(z)$ is the resolvent of $\genmatrix$.

As in, e.g.,~\cite{CaIpKeMe:96,Em:99,Sa:11}, we define the spectral projectors of $\genmatrix$ on \(\Gamma_j\) and \(\widetilde{\Gamma}\) by
\begin{align*}
    \projJ:= \dfrac{1}{2\pi i}\int_{\Gamma_j} \bfR_k(z)
     \dif z,\quad 1\leq j \leq \ell,\quad\text{ and }\quad \projcl:= \dfrac{1}{2\pi i}\int_{\widetilde{\Gamma}} 
     \bfR_k(z)      \dif z.
\end{align*}
Let 
\beq\label{eq:projout}
\projout=\sum_{j=0}^\ell \projJ.
\eeq 
The residue theorem implies that 
\begin{align}\label{eq:projector_relation}
    \projout+\projcl=\bfI
\end{align}
(see, e.g., \cite[Theorem 3.3., Page 53]{Sa:11})
and properties of holomorphic functional calculus imply that
\begin{align*}
\projout \projout =\ \projout,
\quad \projcl \projcl = \projcl,
 \quad \projcl\projout=0,
\quad    \projcl\genmatrix = \genmatrix \projcl,     
\quad    \projout\genmatrix = \genmatrix \projout.
\end{align*}
Let $r_j$ be the index of $\lambda_j$, i.e., the dimension of the largest Jordan block associated with $\lambda_j$. Then 
 \(\range (\projJ)= \ker ((\lambda_j \bfI - \genmatrix)^{r_j})\) where \(r_j\) is the index of \(\lambda_j\), see~\cite[Lemma 3.1]{Sa:11}. 
This last property implies that 
\begin{align}\label{eq:Z_w_minimal_poly}
    \begin{aligned}
        \projJ (\bfI - \lambda_j^{-1}\genmatrix)^{r_j}&=  (\bfI - \lambda_j^{-1}\genmatrix)^{r_j}\projJ = 0,
    \end{aligned}
\end{align}
for $1\leq j\leq \ell$. Finally, let  \(r:= \sum_{j=1}^{\ell}r_j\).

\subsection{The ideas of the proofs}\label{sec:idea}

\paragraph{Idea 1: use the ``cluster plus outliers'' model from \cite{CaIpKeMe:96}.}

The starting point is the bound in the following lemma (proved in \S\ref{sec:proofCampbell} below), which is essentially that in  \cite[Proposition 4.1]{CaIpKeMe:96}.
\begin{lemma}\label{lem:Campbell}
    Let \(\genmatrix\in \bbC^{n\times n}\) and $\widetilde{\Gamma}$  be as in \S\ref{sec:spectral_proj}. 
    If $m\geq r$, then
\begin{align}\label{eq:bound_campbell}
    \dfrac{\lVert \bfr_m (\genmatrix,\bfb,\bfx_0) \rVert_2}{\lVert \bfr_0 (\genmatrix,\bfb,\bfx_0) \rVert_2} 
    \leq \dfrac{1}{2\pi} \lvert \widetilde{\Gamma} \rvert \min_{\substack{p_{m-r}\in \bbP_{m-r}, \\ p_{m-r}(0)=1}} \max_{z \in \widetilde{\Gamma}} \left(\prod_{j=1}^\ell  \dfrac{\lvert \lambda_j - z \rvert }{\lvert \lambda_j \rvert} \lVert \bfR_k(z) \rVert_2   \lvert p_{m-r}(z) \rvert \right).
\end{align}
\end{lemma}

We now  need to use the freedom we have in choosing 
\(\widetilde{\Gamma}\) to bound on this curve the three terms in brackets on the right-hand side of \eqref{eq:bound_campbell}, namely, 
 the distance of the outliers to \(\widetilde{\Gamma}\) (i.e.~$|\lambda_j-z|$), the norm of the resolvent (i.e. \(\lVert \bfR_k(z) \rVert_2\)), and the polynomial \(p_{m-r}\).

\paragraph{Idea 2: choose the shape of $\widetilde{\Gamma}$ so that one can use the min-max result of \cite{BeGoTy:05}.}

To bound the polynomial on the right-hand side of \eqref{eq:bound_campbell}, we use the following result of \cite{BeGoTy:05}. Given a compact set \(K \subset \mathbb{C}\), let
\begin{align}\label{eq:Em}
    E_m(K):=\min_{\substack{p_{m}\in \bbP_{m}, \\ p_{m}(0)=1}} \max_{z \in K} \big| p_{m}(z) \big|
\end{align}

\begin{lemma}[{\cite[Lemma 2.2]{BeGoTy:05}}]\label{lem:bound_polynomial_on_lens}
Given \(\beta\in (0,\pi/2)\), 
    let     \(K_\beta \subset \bbC\)  be defined by 
   \beqs
   K_\beta := \big\{z\,:\, |z|\leq 1\big\}  \cap    \big\{z\,:\, \Re(z) \geq \cos (\beta) \big\}.  
\eeqs
    Let $\gamma_\beta$ be defined by \eqref{eq:gammabeta}.
    Then, for any $m\in \mathbb{Z}^+$,
    \begin{align*}
        E_m(K_\beta) \leq \min \left\{2+\gamma_\beta,\dfrac{2}{1-\gamma_\beta^{m+1}}\right\}\gamma_\beta^m.
    \end{align*}
\end{lemma}

Observe that, for $a\in \Com\setminus\{0\}$, $E_m(aK)= E_m(K)$, since
\beqs
E_m(aK) = \min_{\substack{p_{m}\in \bbP_{m}, \\ p_{m}(0)=1}} \max_{z \in aK} \big| p_{m}(z) \big|
=\min_{\substack{p_{m}\in \bbP_{m}, \\ p_{m}(0)=1}} \max_{\widetilde{z} \in K} \big| p_{m}(a\widetilde{z}) \big|
=\min_{\substack{\widetilde{p}_{m}\in \bbP_{m}, \\ \widetilde{p}_{m}(0)=1}} \max_{\widetilde{z} \in K} \big| \widetilde{p}_{m}(\widetilde{z}) \big|
\eeqs
where $\widetilde{p}_m(\widetilde{z})= p_m(a\widetilde{z})$. 

\begin{figure}[h!]
    \centering
    \includegraphics[width=0.5\textwidth]{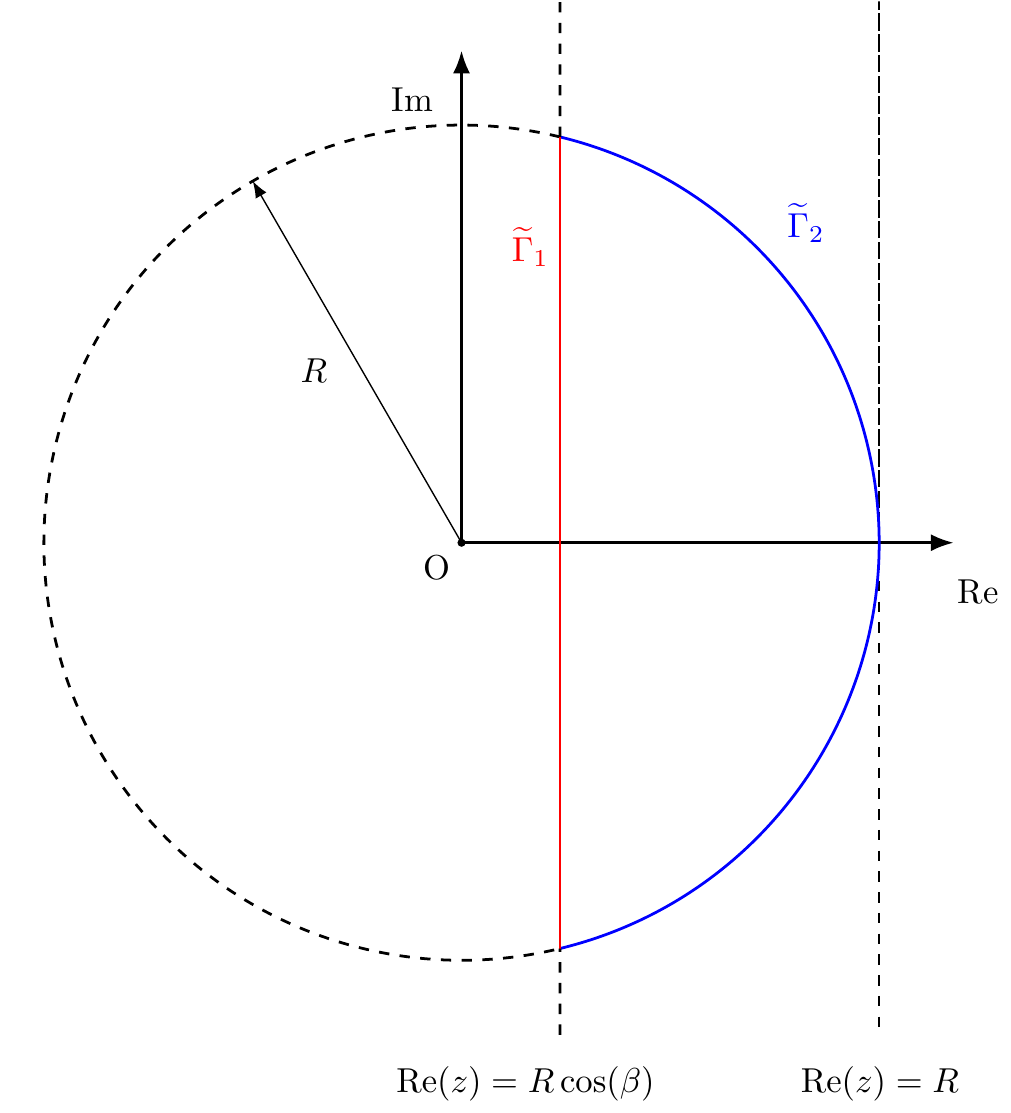}
    \caption{The contour $\widetilde{\Gamma}:= \widetilde{\Gamma}_1 + \widetilde{\Gamma}_2$ in the complex $z$ plane.}\label{fig:lens}
\end{figure}

We therefore choose $\widetilde{\Gamma}$ to be a scaling of $\partial K_\beta$, i.e. $\widetilde{\Gamma} = R(\partial K_\beta)$ for some $R>0$ -- see Figure \ref{fig:lens} --  and use Lemma \ref{lem:bound_polynomial_on_lens} to bound the term involving $p_{m-r}$ in \eqref{eq:bound_campbell}.

\paragraph{Idea 3: choose the parameters defining $\widetilde{\Gamma}$ to control the resolvent.}

We choose $\widetilde{\Gamma}$ to be the boundary of $R K_\beta$ with $R:= \|\genmatrix\|_2+\delta$ for some $\delta>0$. This choice of $R$ ensures that $B(0,R)$ encloses all the eigenvalues of $\genmatrix$. 

We now use the freedom we have in choosing $\delta>0$ and $\beta \in (0,\pi/2)$ to control $\|\bfR(z)\|_2$ on $\widetilde{\Gamma}$. 
Since 
\beqs
\big\|(z \bfI - \genmatrix)^{-1}\big\|_2 \leq \frac{1}{|z| - \|\genmatrix\|_2} \quad\tfor |z| \geq \|\genmatrix\|_2,
\eeqs
we have
    \begin{align}\label{eq:resolvent_bound}
        \lVert \bfR_k(z) \rVert_2\leq \delta^{-1} \quad \tfor z \in \partial B(0,R).
    \end{align}

To bound $\lVert \bfR_k(z)\rVert_2$ on the straight part of \((\lVert \genmatrix\rVert_2+\delta)K_\beta\) we 
choose $\delta$ and $\beta$ so that this straight parts avoids 
the $\delta$-pseudospectrum of $\genmatrix$, $\Lambda_{\delta}(\genmatrix)$, defined by
\begin{align}\label{eq:pseudospectrum}
    \Lambda_{\delta}(\genmatrix):= \big\{z \in \bbC \,\, :\,\, \lVert (z \bfI - \genmatrix)^{-1}\rVert_2\geq \delta^{-1} \big\}.
\end{align}
Avoiding $\Lambda_{\delta}(\genmatrix)$ is possible with $\delta$ sufficiently small because of the following result.

\begin{theorem}[{Bauer-Fike-type theorem \cite[Theorem 52.2]{TrEm:05}}]\label{thm:BF1}
If $\genmatrix \in \Com^{n\times n}$ has $n$ simple eigenvalues, then, for all $\delta>0$,
\beqs
\Lambda_{\delta}(\genmatrix) \subseteq \bigcup_{j=1}^n \Big(\lambda_j+ 
B\big(0, \delta \,n\, \kappa(\lambda_j)\big)\Big).
\eeqs
\end{theorem}

The price we pay for using this general result is that $\delta$ can potentially be small. However, when the resulting bound is applied to Helmholtz BEM, the smallness of $\delta$ is not the limiting factor in the $k$-dependence of the bound.

\subsection{Proof of Lemma \ref{lem:Campbell}}\label{sec:proofCampbell}

As in~\cite{CaIpKeMe:96}, we define the minimal polynomial associated with \((\lambda_j)_{j=1}^\ell \). Let
\begin{align}\label{eq:qr}
    q_r(z):=\prod_{j=1}^\ell  (1 -\lambda_j^{-1} z)^{r_j};
\end{align}
observe that \(q_r\in \bbP_r\) (since \(r:= \sum_{j=1}^{\ell}r_j\)) and \(q_r(0)=1\). 
The significance of $q_r(z)$ is shown by the following lemma.

\ble\label{lem:mainevent}
\beqs
\projout q_r(\genmatrix) = 0.
\eeqs
\ele

\bpf
By the definitions of $q_r$ \eqref{eq:qr} and $\projout$ \eqref{eq:projout}, and then \eqref{eq:Z_w_minimal_poly},
\beqs
\projout q_r(\genmatrix) 
= \left(\sum_{m=1}^\ell \bfP_{\Gamma_m}\right)\left( \prod_{j=1}^\ell \big(\bfI - \lambda_j^{-1}\genmatrix\big)^{r_j} \right) 
= \sum_{m=1}^\ell \left( \prod_{j=1}^\ell \underbrace{
\bfP_{\Gamma_m}\big(\bfI - \lambda_j^{-1}\genmatrix\big)^{r_j}
}_{=0  \text{  when  } j=m}
 \right).
\eeqs
\epf

Let \(p_{m-r}\in \bbP_{m-r}\) be any polynomial of order \(m-r\) such that \(p_{m-r}(0)=1\). Let
\beqs
\overline{p}_m (\genmatrix):= q_r(\genmatrix) p_{m-r}(\genmatrix),
\eeqs
so that \(\overline{p}_m\in \bbP_m\) with \(\overline{p}_m(0)=1\). 
Then, by 
\eqref{eq:projector_relation} and Lemma \ref{lem:mainevent},
\begin{align*}
    \overline{p}_m (\genmatrix)&=\projout q_r(\genmatrix) p_{m-r}(\genmatrix) + \projcl q_r(\genmatrix) p_{m-r}(\genmatrix)
= \projcl q_r(\genmatrix) p_{m-r}(\genmatrix).
\end{align*}
Using this in the characterisation \eqref{def:gmres} of the GMRES residual, we find that
\begin{align}\label{eq:local_bound_on_outliers}
    \lVert \bfr_m (\genmatrix,\bfb,\bfx_0) \rVert_2 &\leq \min_{\substack{p_{m-r}\in \bbP_{m-r}, \\ p_{m-r}(0)=1}}\lVert \projcl q_r(\genmatrix)p_{m-r}(\genmatrix) \bfr_0(\genmatrix,\bfb,\bfx_0)\rVert_2.
\end{align}
By the definition of $\projcl$, 
\begin{align*}
    \dfrac{\lVert \bfr_m (\genmatrix,\bfb,\bfx_0) \rVert_2}{\lVert \bfr_0 (\genmatrix,\bfb,\bfx_0) \rVert_2} &\leq \min_{\substack{p_{m-r}\in \bbP_{m-r}, \\ p_{m-r}(0)=1}} \dfrac{1}{2\pi} \int_{\widetilde{\Gamma}}\lVert \bfR_k(z) \rVert_2 \lvert q_r(z)p_{m-r}(z)  \rvert \dif z, 
\end{align*}
and then \eqref{eq:bound_campbell} follows from the definition of $q_r$.

\subsection{Proof of Theorem \ref{thm:main1}}
We first observe that it is sufficient to prove that 
    \begin{align}\label{eq:Tues1}
        \dfrac{\lVert \bfr_m (\genmatrix,\bfb,\bfx_0) \rVert_2}{\lVert \bfr_0 (\genmatrix,\bfb,\bfx_0) \rVert_2} 
        \leq \left(\prod_{j=1}^\ell  \dfrac{1}{\lvert \lambda_j \rvert}\right) \,(\|\genmatrix\|_2 + \delta)^{\ell+1}\,2^{\ell} \,\delta^{-1} (\gamma_{\beta}+2) \gamma_{\beta}^{m-r}.
    \end{align}
Indeed, \eqref{eq:mlowerbound} then follows by using the inequalities $2\leq 3$ and $\gamma_\beta + 2\leq 3$, and noting that the assumption that $\genmatrix$ has simple eigenvalues implies that $r_j=1$ for $j=1,\ldots,\ell$, so $r=\ell$.

To prove \eqref{eq:Tues1}, we start from the bound \eqref{eq:bound_campbell}, and then follow Ideas 2 and 3 in \S\ref{sec:idea}.
Indeed, we set $\widetilde{\Gamma}:= R K_{\widetilde{\beta}}$, where $R:= \|\genmatrix\|_2 + \delta$, with $\widetilde{\beta}\in (0,\pi/2)$ and $\delta>0$ free parameters to be fixed later.
Let $\widetilde{\Gamma}_1$ and $\widetilde{\Gamma}_2$ be as in Figure \ref{fig:lens}.
        
Since the spectrum of $\genmatrix$ is discrete, for \(\delta\) small enough there exists $L$ with $L_0<L<L_1$ such that that line $\Re z = L$ 
does not intersect the $\delta$-pseudospectrum \(\Lambda_\delta(\genmatrix)\). 
Indeed, combining Theorem \ref{thm:BF1} and the definitions of $N_{\rm eig}$ and $\kappa^*$, we see that this is possible if 
\beqs
2 \delta \, n \, \kappa^*\, N_{\rm eig} < L_1-L_0,
\eeqs
and thus certainly if $\delta$ is given by \eqref{eq:delta}.
With this choice of $\delta$ and the associated $L$, 
let $\widetilde{\beta}$ and $\beta$ be defined so that
\beq\label{eq:widetildebeta}
\cos\widetilde{\beta} = \frac{L}{R} \quad\tand\quad \cos\beta = \frac{L_0}{R},
\eeq
and observe that $K_{\widetilde{\beta}}\subset K_{\beta}$ since $L>L_0$.

In summary, $\widetilde{\Gamma}:= R K_{\widetilde{\beta}}$ with $\delta$ defined by \eqref{eq:delta}, 
$R:= \|\genmatrix\|_2 + \delta$, 
 $L$ is defined so that 
the line $\Re z = L$ 
does not intersect the $\delta$-pseudospectrum \(\Lambda_\delta(\genmatrix)\), and
$\widetilde{\beta}$ is defined by \eqref{eq:widetildebeta}.

Having defined $\widetilde{\Gamma}$, we now bound the quantities appearing on the right-hand side of \eqref{eq:bound_campbell}.
Since  $\widetilde{\Gamma}\subset B(0,R)$,
\beq\label{eq:Tues2}
\max_{z \in \widetilde{\Gamma}}\left(\prod_{j=1}^\ell \lvert \lambda_j - z \rvert \right)\leq ( 2R)^\ell  
\quad\ton\quad\widetilde{\Gamma}.
\eeq
Furthermore, the bound \eqref{eq:resolvent_bound} implies that $\lVert \bfR(z) \rVert_2 \leq \delta^{-1}$ on $\widetilde{\Gamma}_2$ and the
choice of $L$ and the definition of $\Lambda_\delta(\genmatrix)$ \eqref{eq:pseudospectrum} implies that $\lVert \bfR(z) \rVert_2 \leq \delta^{-1}$ on $\widetilde{\Gamma}_2$; therefore
\beq\label{eq:Tues3}
\lVert \bfR(z) \rVert_2 \leq \delta^{-1}  \quad\ton \widetilde{\Gamma}.
\eeq
Using \eqref{eq:Tues2} and \eqref{eq:Tues3} in \eqref{eq:bound_campbell}, we find that 
\begin{align*}
\dfrac{\lVert \bfr_m (\genmatrix,\bfb,\bfx_0) \rVert_2}{\lVert \bfr_0 (\genmatrix,\bfb,\bfx_0) \rVert_2} 
\leq \left(\prod_{j=1}^\ell  \dfrac{1}{\lvert \lambda_j \rvert}\right) \, R\, (2R)^\ell  \, \delta^{-1}  \min_{\substack{p_{m-r}\in \bbP_{m-r}, \\ p_{m-r}(0)=1}}  \max_{z \in \widetilde{\Gamma}}  \lvert p_{m-r}(z) \rvert.
\end{align*}
Using the fact that $\widetilde{\Gamma}= \partial( R \,K_\beta)$, 
the definition of $E_{m-r}(K_\beta)$, and the fact that $K_{\widetilde{\beta}}\subset K_{\beta}$, we have 
\begin{align*}
\min_{\substack{p_{m-r}\in \bbP_{m-r}, \\ p_{m-r}(0)=1}}  \max_{z \in \widetilde{\Gamma}}  \lvert p_{m-r}(z) \rvert &\leq \min_{\substack{p_{m-r}\in \bbP_{m-r}, \\ p_{m-r}(0)=1}}  \max_{z \in R K_{\beta}} \lvert p_{m-r}(z) \rvert = E_{m-r}(R K_{\widetilde\beta})
\leq E_{m-r}(R K_{\beta}).
\end{align*}
The result \eqref{eq:Tues1} then follows from using 
Lemma \ref{lem:bound_polynomial_on_lens}.

\begin{remark}[Removing the assumption that the eigenvalues are simple]
\label{rem:BF1}
We assumed that the eigenvalues of $\genmatrix$ were simple to use Theorem \ref{thm:BF1}.
To remove this assumption, 
one can use the Bauer-Fike theorem (see, e.g, ~\cite[Theorem 2.3]{TrEm:05}) that 
\beqs
\Lambda_{\delta}(\genmatrix)\subset \Lambda(\genmatrix)+B(0,\delta \kappa(\bfV)),
\eeqs
 for \(\genmatrix=\bfV {\bf \Lambda} \bfV^{-1}\) with ${\bf \Lambda}$ a diagonal matrix with eigenvalues on the diagonal and \(\bfV\) the corresponding matrix of eigenvectors. 
 Assumption A3 would then be replaced with an assumption that $\kappa(\bfV)$ grows at most polynomially with increasing $k$. 
\end{remark}

\subsection{Proof of Theorem \ref{thm:main2}}

\ble[Bound on $\|A_k'\|_{\LtGt}$] 
\label{lem:normbounds}
If $\Oi$ is piecewise smooth (in the sense of, e.g., \cite[Definition 1.3]{GaSp:19}, then, given $k_0>0$, there exists $C_1$ (depending on $k_0$, $d$, and $\Omega$) such that 
\beq\label{eq:upper}
\N{A_k'}_{\LtGt} \leq C_1 k^{1/2}\log (k+2) 
\quad\tfa k\geq k_0.
\eeq
\ele

\bpf[References for the proof]
This follows from the bounds in 
\cite[Appendix A]{HaTa:15}, \cite{GaSm:15},  \cite[Chapter 4]{Ga:19} (see also \cite{GaSp:19}). 
\epf

We now prove  Theorem \ref{thm:main2}. 
First observe that, by the definitions of $N_{\rm eig}$ and $\ell$ and Assumption A2, 
\beq\label{eq:proof1}
1\leq N_{\rm eig}\leq \ell +1 \leq \CWeyl k^{d-1} +1.
\eeq
By Assumption A1, $S$ is independent of $k$; we then choose $L_0$ and $L_1$ to be also independent of $k$; i.e. $L_0, L_1\sim 1$.  
By Assumption A3 and the definition \eqref{eq:eigenvaluecondition}, $1\leq \kappa^* \leq \Ccond k^M$. 
Finally, by assumption, $1\leq n\leq C_{\rm dis} k^{M'}$.
Using in \eqref{eq:delta} all these inequalities, we find that 
\beq\label{eq:proof2}
\frac{1}{\CWeyl C_{\rm dis} \Ccond} k^{-M-M'-(d-1)} \lesssim \delta \lesssim 1.
\eeq
By the bound \eqref{eq:discon2} from Assumption A0 and then \eqref{eq:upper},
\beq\label{eq:proof3}
\N{\bfM^{-1}\matrixD}_{2} \leq \Capproxone \N{A_k'}_{\LtGt}\lesssim  \Capproxone k^{1/2}\log (k+2),
\eeq
where the omitted constant depends only on $\Oi$.
Using  \eqref{eq:proof2} and \eqref{eq:proof3}  in  \eqref{eq:mlowerbound} (and recalling the asymptotics \eqref{eq:Wed1}), we obtain that 
if $m$ satisfies 
\begin{align}
 m\geq \CWeyl k^{d-1}+  C_1 \N{A_k'}_{\LtGt}
 \left(
\sum_{j=1}^\ell
 \log\dfrac{1}{\lvert \lambda_j \rvert}
 + 
 \log (\varepsilon^{-1})+ C\log k +
C_2 k^{d-1}\log k
 \right)
 \label{eq:proof4}
    \end{align}
          then the $m$th GMRES residual satisfies $\|\bfr_m \|_2/\|\bfr_0 \|_2 \leq \eps$; in \eqref{eq:proof4}, $C_1$ depends on $\Capproxone$, $\Oi$ and $d$, $C_2$ depends on $\CWeyl, \Capproxone$, $\Oi$, and $d$, and $C$ depends on $\CWeyl$, $\Ccond$, $C_{\rm dis}$, $M$, $M'$, and $d$. 
     The bound on $m$ \eqref{eq:mlowerbound2} now follow from absorbing the $\log k$ term into the $k^{d-1}\log k$ term, and modifying the definition of $C_2$ appropriately, so that it now depends also on  $\Ccond$, $C_{\rm dis}$, $M$, $M'$, and $d$. 

\bre[The analogue of Theorem \ref{thm:main2} for $B_{k, {\rm reg}}$]\label{rem:Neumann}
If $B_{k, {\rm reg}}$ satisfies appropriate analogues of Assumptions A0-A3, then Theorem \ref{thm:main2} holds with $A_k'$ replaced by $B_{k, {\rm reg}}$ (and the Galerkin matrices modified accordingly). Indeed, the bound 
\beq\label{eq:upperN}
\N{B_{k, {\rm reg}}}_{\LtGt} \leq C_1 k^{1/4}\log (k+2) 
\quad\tfa k\geq k_0
\eeq
for piecewise-smooth $\Oi$ is proved in  \cite{GaMaSp:21N} using results from \cite[Appendix A]{HaTa:15}, \cite{GaSm:15},  \cite[Chapter 4]{Ga:19} and results about semiclassical pseudodifferential operators (to bound $S_{\ri k}$.
The proof of Theorem \ref{thm:main2} therefore goes through for $B_{k,{\rm reg}}$ with \eqref{eq:upperN} replacing \eqref{eq:upper}.

In \S\ref{sec:discussion}, we recalled that $\|(A_k')^{-1}\|_{\LtGt}\lesssim \exp(\alpha k)$ when $\Oi$ is smooth, with this proved in 
 \cite[Equation 1.35 and Lemma 6.2]{ChSpGiSm:20} by combining bounds on the exterior Dirichlet problem from \cite[Theorem 1.2]{Bu:98} and bounds on the interior impedance problem in \cite[Theorem 1.8 and Corollary 1.9]{BaSpWu:16}.
The analogous bound for $B_{k, {\rm reg}}$ is proved in \cite{GaMaSp:21N} by expressing $(B_{k, {\rm reg}})^{-1}$ in terms of appropriate exterior and interior solution operators using \cite[Lemma 6.1, Equation 83]{BaSpWu:16}, 
and then using the bounds on the exterior Neumann problem from \cite{Bu:98}/\cite{Vo:00} and the relevant interior problem (an interior impedance-like problem involving $S_{\ri k}$ in the boundary condition) from the combination of \cite[Theorem 4.6]{GaLaSp:21} and  \cite[Lemma 3.2]{GaMaSp:21}.
\ere

\section{Weyl asymptotics, Assumption A2, and why the large cavity has more near-zero eigenvalues than the small cavity}
\label{sec:Weyl_new}

Recall that Assumption A2 is that the number of near-zero eigenvalues of $\bfM^{-1}\matrixD$ is $\lesssim k^{d-1}$, where the near-zero eigenvalues are defined as those in the set $k$-independent open set $\cN$ in O1/A1.

In this section we use Weyl asymptotics to give heuristic arguments about why this assumption holds, and also why the large cavity has more near-zero eigenvalues than the small cavity. We begin by giving numerical evidence that Assumption A2 holds.

\subsection{Numerical evidence for Assumption A2}\label{sec:numA2}

\begin{figure}[h!]
    \centering
    \includegraphics[width=0.9\textwidth]{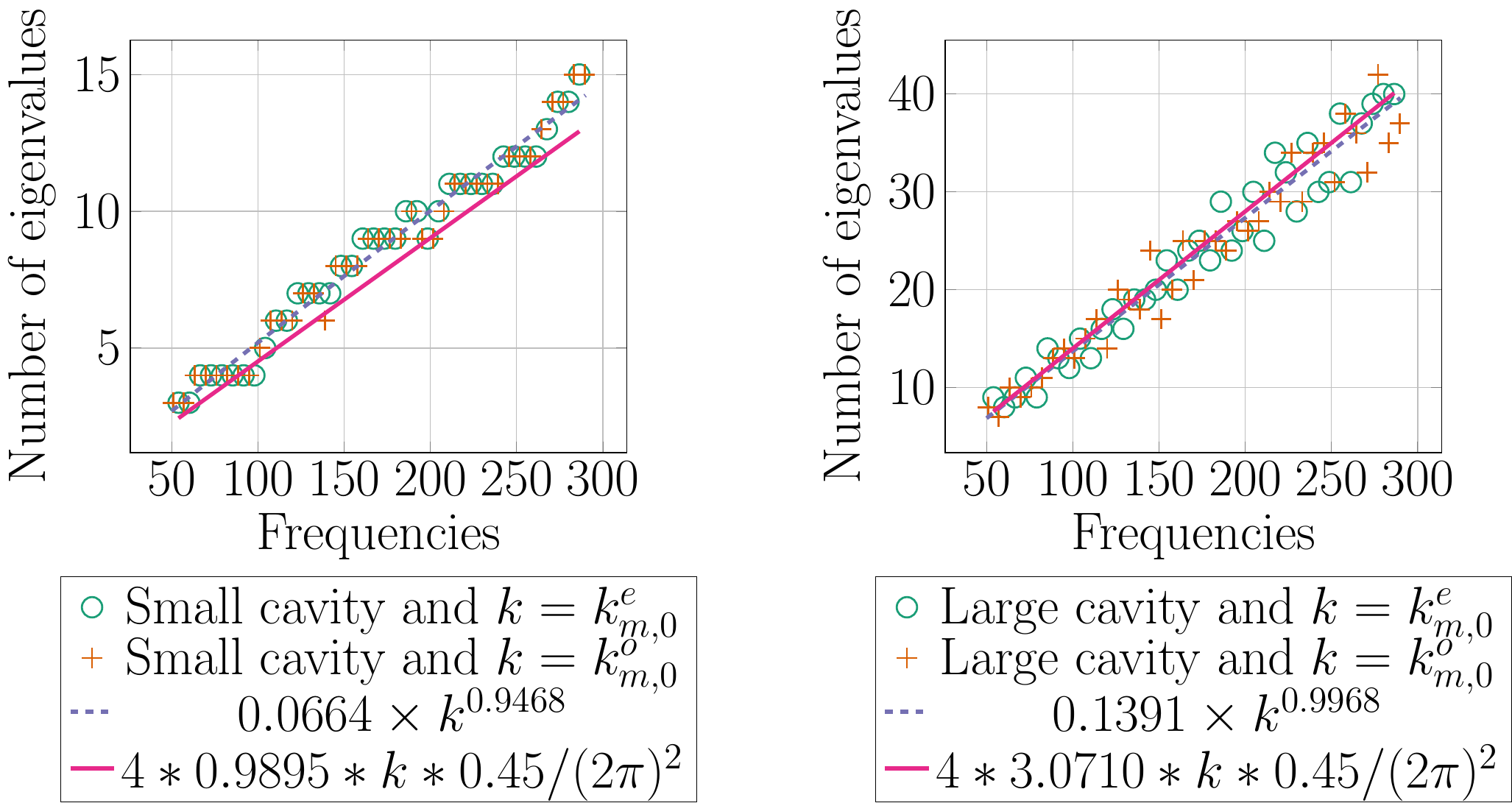}
    \caption{The number of eigenvalues of $\bfM^{-1}\matrixD$ in the rectangle $[-0.1, 0.1]\times [-0.6,0.6]$  plotted against $k=k_{m,0}^e$ and $k=k^o_{m,n}$ for both the small and large cavities.}\label{fig:counting_function}
\end{figure}

Figure \ref{fig:counting_function} plots the the number of eigenvalues of $\bfM^{-1}\matrixD$ in the rectangle  $[-0.1, 0.1]\times [-0.6,0.6]$ as a function of $k$ (plotted at $k=k_{m,0}^e$ and $k=k^o_{m,n}$) for both the small and large cavities in 2-d \footnote{Since the number of eigenvalues is discrete, the rectangle needs to be sufficiently large so that it includes enough eigenvalues for the Weyl asymptotics to hold.}.
The dashed lines are best-fit lines fitted using the nonlinear least-squares Marquardt-Levenberg algorithm; 
we see the number of eigenvalues growing very close to linearly with $k$. 
The same experiments for the Galerkin matrices (preconditioned with the mass matrix) of the operators $B_k$ and $B_{k, {\rm reg}}$ result in very similar plots; we do not show them here, but the exponents in the best fit lines are displayed in Table \ref{tab:kdependence}.

The solid lines in Figure \ref{fig:counting_function} show how many Laplace eigenfunctions of the ellipse $E$ \eqref{eq:ellipse} are localised in the respective cavity; we explain in \S\ref{sec:smalllarge} below how we calculate this.
These localised eigenfunctions produce quasimodes (in the sense of Definition \ref{def:quasimodes}) with small quality.
Figure \ref{fig:counting_function} therefore gives strong evidence for 
near-zero eigenvalues of $\bfM^{-1}\matrixD$ in $\cN$ correspond to localised eigenfunctions of the Laplacian in the ellipse.

\subsection{Recap of Weyl asymptotics for the number of eigenvalues of the Laplacian}\label{sec:Weyl_recap}
On a bounded domain $\Omega\subset \mathbb{R}^d$, the standard Weyl law  states that, for either the Dirichlet or Neumann problem in $\Omega$,
\begin{equation}
\label{eq:WeylStandard}
N_{\Omega}(\Lambda)= \frac{\Lambda^{d/2}\omega_d\operatorname{Vol}(\Omega)}{(2\pi )^d}+O\big(\Lambda^{(d-1)/2}\big) 
 \quad\tas\Lambda\tendi;
\end{equation}
see e.g.~\cite[Theorem 1.2.1]{SaVa:92}, \cite{Se:78,Ph:81}.
In fact, if the periodic billiard trajectories on $\Omega$ with speed one 
form a set of zero measure, 
\cite{Iv:80} (see also~\cite[Corollary 29.1.6]{Ho:85a}~\cite[Theorem 1.6.1]{SaVa:92}) proved that 
\beq\label{eq:Weylimproved}
N_{\Omega}(\Lambda)= \frac{\Lambda^{d/2}\omega_d\operatorname{Vol}(\Omega)}{(2\pi )^d}\pm \frac{\omega_{d-1}}{4(2\pi )^{d-1}}\Lambda^{(d-1)/2}\operatorname{Vol}(\partial\Omega)+o\big(\Lambda^{(d-1)/2}\big) \quad\tas\Lambda\tendi,
\eeq
where $\omega_d$ is the volume of the unit ball in $\mathbb{R}^d$, the plus sign is taken for the Dirichlet problem, and the minus sign for the Neumann problem.

\subsection{Two non-rigorous arguments about why we expect A2 to hold}

The first uses Observations O2(a)-(c) from \S\ref{sec:observations} along with results from \cite{St:99} linking the number of quasimodes to the number of resonances, and Weyl-type bounds on the number of resonances from \cite{PeZw:99}.

The second uses the Weyl asymptotics for the eigenfunctions of the Laplacian on $\Gamma$ (i.e.~the surface Laplacian) along with results from \cite[Chapter 4]{Ga:19} about the properties of the boundary-integral operators $D_k'$ and $S_k$ as semiclassical pseudodifferential operators. 

\subsubsection{Why we expect A2 to hold when there exist quasimodes with $\QMC(k)=O(k^{-\infty})$.}\label{sec:heuristic1}

We use the notation that 
$a = O(k^{-\infty})$ as $k\tendi$ if, given $N>0$, there exists $C_N$ and $k_0$ such that $|a|\leq C_N k^{-N}$ for all $k\geq k_0$, i.e.~$a$ decreases superalgebraically in $k$.

The steps in this argument are as follows:
\ben
\item We assume that there exists a one-to-one mapping between eigenvalues of $\bfM^{-1}\matrixD$ in $\cN$ and quasimode frequencies $k_j$, such that an eigenvalue corresponding to frequency $k_{j^*}$ is closest to the origin when $k=k_{j^*}$.
\item By Point 1, O1/A1, and O2(c) (the eigenvalues move at $O(1)$ speed), the number of eigenvalues in $\cN$ at $k$ equals the number of quasimode frequencies $k_j$ in 
an interval $[k+p,k-p]$ with $p$ independent of $k$. 
\item If the quality of the quasimode is $O(k^{-\infty})$, then the number of quasimode frequencies $k_j$ in 
an interval $[k+p,k-p]$, with $p$ independent of $k$, is $\lesssim k^{d-1}$.
\item By Points 2 and 3, the number of eigenvalues in $\cN$ at $k$ is $\lesssim k^{d-1}$.
\een
Regarding Point 2: Recall that $\cN$ is independent of $k$ by O1/A1. By O2(c), the eigenvalues in $\cN$ move at $k$-independent speed, and we assume that the paths of the eigenvalues are similar to those in Figure \ref{fig:outliers} and \ref{fig:flow_eig} (i.e.,~the eigenvalues don't move, e.g., in circles in $\cN$).
Therefore, the number of eigenvalues in $\cN$ at $k$ is equal to the number of eigenvalues that pass close to zero in an interval $[k+p,k-p]$ with $p$ independent of $k$. By Point 1, this number is equal to the number of $k_j$ in $[k+p,k-p]$.

Regarding Point 3:
If the quality of the quasimode is $O(k^{-\infty})$, then, by \cite[Theorem 2 and Corollary 2]{St:99}, the number of quasimodes is bounded by the number of resonances 
in an $O(k^{-\infty})$ neighbourhood below the real axis. 
Using the notation 
\beq\label{eq:countingN}
N(K):= \big|\big\{ k_j \leq K\, :\, k_j \text{ is a quasimode frequency}\big\}\big|,
\eeq
where $|\cdot|$ denotes cardinality of a set, the bound 
\beq\label{eq:counting2}
N(k+p) - N(k-p) \lesssim k^{d-1}
\eeq
follows since, by
\cite[Proposition 2]{PeZw:99}, the counting function of the number of resonances in an $O(1)$ neighbourhood of the real axis satisfies \eqref{eq:counting2}.
(Note that the assumption~\cite[Equation 1.9]{PeZw:99} about the Weyl asymptotics for the reference operator in the black-box framework holds by the results recapped in \S\ref{sec:Weyl_recap}).

\bre[When do there exist quasimodes with $\QMC(k)=\mathcal{O}(k^{-\infty})$?]
There exist quasimodes with $\QMC(k)=\mathcal{O}(k^{-\infty})$ if $\Oi$ satisfies the assumptions of Theorem \ref{thm:ellipse}, and also in the following two situations
by \cite[Theorem 1]{CaPo:02} and \cite[Theorem 1]{St:00} respectively.

(i) $\Gamma$ has zero Dirichlet boundary conditions and $\Omega_+$ contains an elliptic-trapped ray such that (a) $\Gamma$ is analytic in a neighbourhood of the ray and (b) the ray satisfies the stability condition \cite[(H1)]{CaPo:02}. In this situation, if  $q >11/2$ when $d=2$ and $q>2d+1$ when $d\geq 3$, then there exists a family of quasimodes (in the sense of Definition \ref{def:quasimodes}) with 
\beqs
\QMC(k)= C_1 \exp( - C_2 k^{1/q})
\eeqs
for some $C_1, C_2>0$  and independent of $k$.

(ii) There exists a sequence of resonances $\{\lambda_\ell\}_{\ell=1}^\infty$ of the exterior Dirichlet/Neumann problem with
\beqs
0\leq -\Im \lambda_\ell = \mathcal{O}\big(|\lambda_\ell|^{-\infty}\big)  \quad\tand \quad \Re \lambda_\ell \tendi \quad\tas\quad \ell \tendi
\eeqs
(recall that the resonances of the exterior Dirichlet/Neumann problem are the poles of the meromorphic continuation of the solution operator 
from $\Im k\geq 0$ to $\Im k<0$; see, e.g., \cite[Theorem 4.4. and Definition 4.6]{DyZw:19}).
\ere

\subsubsection{A second argument why we expect Assumption A2 to hold.}
We consider the case when $\Gamma\in C^\infty$ since a great deal of information is then available about the structure of $D_k'$ and $S_k$ (see~\cite[Chapter 4]{Ga:19}). In particular, these operators have the following two important features.
\begin{itemize}
\item[(i)]  For any $\eps>0$, there is $M>0$ such that if $v$ is a function with frequency $\geq Mk$, then
$$
\|D_k'v\|_{L^2} +k\|S_k v\|_{L^2}\leq \eps \|v\|_{L^2}.
$$
\item[(ii)] $S_k$ and $D_k'$ almost map the Hilbert space of functions with frequency $\geq Mk$ to itself.
\end{itemize}
(When we say ``a function $v$ with frequency $\geq Mk$'' we mean that $v= \sum_{\lambda_j \geq M k} a_{\lambda_j}\phi_{\lambda_j}$ for some $a_{\lambda_j}\in\mathbb{C}$, where $(-\Delta_g -\lambda_j^2 )\phi_{\lambda_j} =0$ are the eigenfunctions of $-\Delta_g$ on $\Gamma$.)

If (ii) were exactly true (i.e., $S_k$ and $D_k'$ exactly preserve the space of functions with frequency larger than $Mk$) then we could decompose $L^2=F_{\leq Mk}\oplus F_{>Mk}$ where $F_{\leq MK}$ denotes the Hilbert space of functions with frequency $\leq Mk$, and  $F_{> Mk}$ its orthogonal complement. In particular, since $F_{\leq Mk}$ and $F_{>Mk}$ would be invariant under the action of $A_k'$, the eigenvalues of $A_k'$ would then be the union of the eigenvalues of
\beqs
A'_{<}:=\Pi_{F_{\leq Mk}}A'_k\Pi_{F_{\leq Mk}}:F_{\leq Mk}\to F_{\leq Mk}
\eeqs
 and
 \beqs
 A'_>:=\Pi_{F>Mk}A'_k\Pi_{F>Mk}:F_{>Mk}\to F_{>Mk}.
 \eeqs
Then, by (i) all of the eigenvalues of $A'_{>}$ would be contained in $B(\frac{1}{2},\eps)$ and hence eigenvalues outside this ball would correspond to eigenvalues of $A'_{<}$.  By Weyl asymptotics, $\dim_{F_{\leq Mk}}\leq Ck^{d-1}$; therefore the number of eigenvalues of $A'_k$ outside the ball of radius $\eps$ around $\frac{1}{2}$ would be $\lesssim k^{d-1}$.

We now make this heuristic argument slightly more precise. Let $-\Delta_g$ be the Laplacian on $\Gamma$, and $\chi\in C_c^\infty((-2,2);[0,1])$ with $\chi \equiv 1$ on $[-1,1]$. Then, writing
$$
R_k:=(D_k'-\ri k S_k)(I-\chi(-k^{-2}\Delta_g)),
$$
we have that $R_k$ is semiclassical pseudodifferential operator of order $-1$. Furthermore, inspecting the semiclassical principle symbols of $D_k'$ and $S_k$, we see that there is $f$ such that $|\langle t\rangle f(t)|\leq C$ and 
\begin{align*}
&\big\|\langle -k^{-2}\Delta_g\rangle^{1/2}(R_k-f(-k^{-2}\Delta_g)))\big\|_{\LtGt}\\
&\qquad\qquad+\big\|(R_k-f(-k^{-2}\Delta_g)))\langle -k^{-2}\Delta_g\rangle^{1/2}\big\|_{\LtGt}\leq Ck^{-1}.
\end{align*}
Furthermore, for any $\eps>0$ there exists $M$ large enough such that
\begin{equation}
\label{e:oneAndTwo}
\begin{gathered}
\big\|(D_k'-\ri k S_k)(I-1_{[-M^2,M^2]}(-k^{-2}\Delta_g))\big\|_{\LtGt}< \eps,\\
k\big\|\big[(D_k'-\ri k S_k),1_{[-M^2,M^2]}(-k^{-2}\Delta_g)\big]\big\|_{\LtGt}<\eps,
\end{gathered}
\end{equation}
and then both (i) and (ii) above follow from \eqref{e:oneAndTwo}.

Define $F_{>Mk}$ to be the cokernel of $1_{[-M^2,M^2]}(-k^{-2}\Delta_g)$ in $L^2(\Gamma)$ and then $F_{\leq Mk}$ its orthogonal complement. Then let $\Pi=1_{[-M^2,M^2]}(-k^{-2}\Delta_g)$ be the orthogonal projector onto $F_{\leq Mk}$. 
By \eqref{e:oneAndTwo}, 
$$
D_k'-\ri k S_k=\widetilde{D}'_k-ik\widetilde{S}_k+O(k^{-1}\eps)_{L^2\to L^2},\qquad \widetilde{D}_k'-ik\widetilde{S}_k:= (I-\Pi)(D_k'-\ri k S_k)(I-\Pi)+\Pi (D_k'-\ri k S_k)\Pi.
$$

We now argue with $D_k'-\ri k S_k$ replaced by $\widetilde{D}_k-ik\widetilde{S}_k$ and choose $\eps<\frac{1}{4}$ in~\eqref{e:oneAndTwo}. In this case,  we can orthogonally decompose $L^2(\Gamma)$ into the subspaces $F_{\leq Mk}$ and $ F_{>Mk}$ which are invariant under application of $\tilde{A}_k':=\frac{1}{2}I+\widetilde{D}'_k-ik\widetilde{S}_k$. Then, since
$$
(\tfrac{1}{2}-z)I+(I-\Pi)(D_k-\ri k S_k)(I-\Pi)
$$
is invertible for $|z-\frac{1}{2}|> \frac{1}{4}$, $\tilde{A}_k'$ has at most $\dim F_{\leq Mk}$ eigenvalues in $|z-\frac{1}{2}|>\frac{1}{4}$. 

By the Weyl law on $\Gamma$ (which follows from \cite{Ho:68,Le:52,Av:56} since $\Gamma$ has no boundary),
$$
\dim F_{\leq MK}=\big|\{ \lambda_j\leq Mk\,:\, \lambda_j^2\text{ is an eigenvalue of }-\Delta_g\}\big|= \frac{\operatorname{Vol}_g(\Gamma) \omega_{d-1}}{(2\pi )^{d-1}}(Mk)^{d-1}+O(k^{d-2}),
$$
and, in particular, $\tilde{A}_k'$ has at most $Ck^{d-1}$ eigenvalues in $|z-\frac{1}{2}|>\frac{1}{4}$. 
\begin{remark}
Although the difference is small, replacing $D_k'-\ri k S_k$ by $\widetilde{D}_k'-\ri k\widetilde{S}_k$ as we did in the arguments above is a \emph{serious} simplification and a more sophisticated argument would be needed to obtain a genuine bound on the number of eigenvalues away from $z=1/2$.
\end{remark}

\subsection{Why $\bfM^{-1}\matrixD$ has more near-zero eigenvalues for the large cavity than the small cavity}\label{sec:smalllarge}

\paragraph{How the pink lines in Figure \ref{fig:counting_function} were determined.}

Figure \ref{fig:counting_function} shows that the number of eigenvalues of $\bfM^{-1}\matrixD$ in $\cN= [-0.1,0.1]\times[-0.6,0.6]$ grows with $k$ for both the small and large cavities, but the rate of growth is higher for the large cavity than the small cavity.
Recall that the pink lines in Figure \ref{fig:counting_function} show how many eigenfunctions of the Laplacian in the ellipse are localised in the respective cavity.

In \S\ref{sec:heuristic1} we assumed that all the near-zero eigenvalues of $\bfM^{-1}\matrixD$ correspond to quasimode frequencies $k_j$, and we argued that 
\beq\label{eq:Sunday1}
\big| \big\{ \lambda \text{ eigenvalue of } \bfM^{-1}\matrixD  \, :\, \lambda \in \cN \big\}\big|
= N(k+p)- N(k-p),
\eeq
for an appropriate $p$, where $N(k)$ is given by \eqref{eq:countingN}; i.e., $N(k)$ is the counting function of the quasimode frequencies.

We assume further that all the quasimode frequencies correspond to eigenvalues of Laplacian in the ellipse 
$E$ \eqref{eq:ellipse} whose eigenfunctions localised about the minor axis, so that
\beq\label{eq:Sunday2}
N(k+p)- N(k-p)= N_{\rm loc}(k+p)- N_{\rm loc}(k-p),
\eeq
where $N_{\rm loc}$ is the counting function of these eigenvalues of the Laplacian.

We now use a microlocal version of the Weyl asymptotics \eqref{eq:Weylimproved} to determine the asymptotics of $N_{\rm loc}(k)$.
Assume that the ellipse is cut at $(x_{\rm cut}, y_{\rm cut})$, so that  the small cavity corresponds to 
$x_{\rm cut}=-\cos(3\pi/10)=\cos(7\pi/10)$, and the large cavity corresponds to $x_{\rm cut}= -\cos(\pi/10)=\cos(9\pi/10)$; see Figure \ref{fig:geometries}. 
Let 
\beq\label{eq:alphacut}
\alpha_{\rm cut}:= - \left( 1- \frac{x_{\rm cut}^2}{a_1^2}\right) = -\frac{y_{\rm cut}^2}{a_2^2}
\eeq
and let $a:= \sqrt{a_1^2-a_2^2}$ (as in Appendix \ref{sec:Mathieu}).
We show below that the asymptotics of $N_{\rm loc}$ for eigenfunctions localised in the cut ellipse is given by 
\beq\label{eq:Nloc}
N_{\rm loc}(K) = \frac{V_{\rm loc}(\alpha_{\rm cut})}{(2\pi)^2} K^d + c_1 K^{d-1} + o\big(K^{d-1}\big) \quad\tas K\tendi,
\eeq
where
\begin{align}
V_{\rm loc}(\alpha):=8\int^{\pi/2}_{\arcsin(\sqrt{-\alpha})}\int_0^{\cosh^{-1}(\frac{a_1}{a})}
\Phi(\alpha,\omega,\theta)\, a^2(\sinh^2(\omega)+\sin^2(\theta)) \rd\omega \rd\theta,
\label{eq:Vloc}
\end{align}
where
\beq\label{eq:Phi}
\Phi(\alpha,\omega,\theta):=
\begin{cases}\arcsin\left(\sqrt{ \frac{1}{\sinh^2(\omega)+\sin^2(\theta)}\Big( \alpha+\sin^2\theta\Big)}\right)& \text{ if } \sin^2(\theta)\geq -\alpha\\0&\text{ otherwise.}
\end{cases}
\eeq
Calculating the integral in \eqref{eq:Vloc}, we find $V_{\rm loc}= 0.9895$ for the small cavity and $V_{\rm loc} = 3.0710$ for the large cavity. 

Then, combining \eqref{eq:Sunday1} and \eqref{eq:Sunday2}, we find that
\begin{align}\nonumber
\big| \big\{ \lambda \text{ eigenvalue of } \bfM^{-1}\matrixD  \, :\, \lambda \in \cN \big\}\big|
&= N(k+p)- N(k-p)\\
&= N_{\rm loc}(k+p)- N_{\rm loc}(k-p) = 2 p d \frac{V_{\rm loc}}{(2\pi)^2} K^{d-1}  + o\big(K^{d-1}\big).
\label{eq:blueline}
\end{align}
We now determine an appropriate value of $p$ when $\cN= [-0.1,0.1]\times[-0.6,0.6]$ (since this is the $\cN$ we chose in Figure \ref{fig:counting_function}).
Figure \ref{fig:outliers} indicates that the eigenvalues all move on roughly the same trajectory through $\cN$. 
Since the eigenvalues move with speed observed numerically to be approximately one, the appropriate $p$ is half the length of the portion of the curve that intersects $\cN$.
With the imaginary part the $x$ variable and the real part the $y$ variable, we fit a polynomial of degree two in $x$ to 
this portion of the curve, and find its length to be $0.90$, i.e., we take $p=0.45$.
The pink lines in Figure \ref{fig:counting_function} are then the linear function of $k$ on the right-hand side \eqref{eq:blueline} with $p=0.45$ and $d=2$. 
As mentioned above, the fact the these pink lines match so well the number of eigenvalue of $\bfM^{-1}\matrixD$ in $\cN$ give 
strong evidence for the assumptions that (i) all near-zero eigenvalues of $\bfM^{-1}\matrixD$ correspond to quasimodes, and 
(ii) the majority of quasimodes correspond to localised eigenfunctions of the Laplacian in the ellipse.

\paragraph{How we obtained \eqref{eq:Nloc} and \eqref{eq:Vloc}.}
In \S\ref{sec:Weyl_recap} we recapped the standard \eqref{eq:WeylStandard} and improved \eqref{eq:Weylimproved} Weyl asymptotics  for the number of eigenvalues of the Laplacian on a bounded domain. Furthermore, there are microlocal versions of~\eqref{eq:Weylimproved} (see e.g., \cite[Theorems 1.8.5, 1.8.7]{SaVa:92}) that can be integrated to state, roughly, that the total $L^2$ mass of the normalised eigenfunctions with eigenvalue $\leq \Lambda$ in any subset, $U\Subset \{x\in \Omega,|\xi|\leq 1\}$ is given by 
\beq\label{eq:WeylimprovedMicrolocal}
M_{\rm eig}(U,\Lambda) = \frac{\int 1_{U}(x,\xi)\rd x\rd\xi}{(2\pi)^d} \Lambda^{d/2} +c_1(U)\Lambda^{(d-1)/2}+o\big(\Lambda^{(d-1)/2}\big) \quad\tas \Lambda\tendi,
\eeq
where $x$ is the position variable, and $\xi$ the momentum variable.
(For domains without boundary, these estimates can be recovered from 
\cite{DuGu:75} and the full statement together with more quantitative versions can be found in~\cite[Theorem 6]{CaGa:20}.)
We now apply \eqref{eq:WeylimprovedMicrolocal} to the ellipse. 
We could not find a proof that the periodic billiard trajectories with speed one on an ellipse  
form a set of zero measure, under which the improved Weyl asymptotics \eqref{eq:Weylimproved}/\eqref{eq:WeylimprovedMicrolocal} hold. However, the results of \cite[\S4]{Co:10} indicate that the counting function of eigenvalues of the ellipse does indeed satisfy the improved Weyl asymptotics \eqref{eq:Weylimproved}/\eqref{eq:WeylimprovedMicrolocal}.

When $\Omega\subset \mathbb{R}^2$ is the ellipse, the Laplacian is quantum 
completely integrable (see, e.g., \cite{GaTo:20} for the definition of quantum complete integrability); one consequence of this is the separation of variables used in Section~\ref{sec:Mathieu}.  Moreover, this complete integrability implies the existence of a basis of eigenfunctions that concentrate along integrable tori. In particular, one expects that for $U$ a union of integrable tori, 
\begin{equation}
\label{e:localisedWeyl}
M_{\rm eig}(U,\Lambda)\sim \big|\big\{\lambda_j^2\leq \Lambda\,:\, u_{\lambda_j}\text{ is localised inside }U\big\}\big|
\end{equation}
(given $N>0$, we say that $u_{\lambda_j}$ is localised inside $U$ if, for all symbols $a\in C_0^\infty(\overset{\circ}{U^c})$, $\|{\rm Op}_{\lambda_j^{-1}}(a) u_{\lambda_j} \|_{L^2} \leq \lambda_j^{-N}\|u\|_{L^2}$, where ${\rm Op}_{\lambda_j^{-1}}$ is semiclassical quantisation; see, e.g., \cite[Chapter 4]{Zw:12}, \cite[Page 543]{DyZw:19}).

In addition to the the fact that we could not find a reference for the improved Weyl law on the ellipse, this last step is also non-rigorous; a sophisticated analysis of the quantum completely integrable system would be required to justify~\eqref{e:localisedWeyl}.

The integrable tori on the ellipse correspond to billiard trajectories that remain tangent to some confocal conic; see, e.g., \cite[Page 8]{Ly:19}. The eigenfunctions that localise inside the elliptic cavity correspond precisely to eigenfunctions localised on integrable tori generated by confocal hyperbole that intersect the boundary of the ellipse where it has not been truncated; see Figure~\ref{fig:integrableTori}. In Appendix~\ref{sec:Weyl}, we compute the volume, $V_{\rm loc}$, in phase space of the integrable tori contained entirely inside the 
small and large cavities and show that \eqref{eq:Vloc} holds.
The fact that these are \emph{all} the eigenfunctions localised inside the cavity (by \eqref{e:localisedWeyl}) then implies that \eqref{eq:Nloc} holds.
(We note that a similar phase-space volume calculation in an integrable setting occurs in \cite[Appendix B]{BaBe:07}.)

  \begin{figure}
    \centering
    
    \begin{tikzpicture}
    \def \a{1.3};
    \def \b{.6};
    \def \t{45};
    \def \d{1}
    \def \ta{65}
    \begin{scope}[scale=2 ]
    \draw[thick,black] (0,0) ellipse (1 and .5);
    \draw[thick,black] (0,0) ellipse ({\a} and {\b});
	\fill[white] ({cos(126)},-.7)-- ({cos(126)},.7)-- (-\a-.1,.7)--(-\a-.1,-.7)--cycle;
    \draw[thick,black] ({cos(126)},{.5*sqrt(1-cos(126)^2)})--({cos(126)},{\b*sqrt(1-cos(126)^2/(\a)^2)});
    \draw[thick,black] ({cos(126)},-{.5*sqrt(1-cos(126)^2)})--({cos(126)},-{\b*sqrt(1-cos(126)^2/(\a)^2)});
    
    \draw[scale=1, domain=-\d:\d, smooth, variable=\x, blue, dashed,thick] plot ({sqrt(.75)*cos(\t)*(pow(2,\x)+pow(2,-\x))/2)},{sqrt(.75)*sin(\t)*(pow(2,\x)-pow(2,-\x))/2} );
    \draw[scale=1, domain=-\d:\d, smooth, variable=\x, blue, dashed,thick] plot (-{sqrt(.75)*cos(\t)*(pow(2,\x)+pow(2,-\x))/2)},{sqrt(.75)*sin(\t)*(pow(2,\x)-pow(2,-\x))/2} );
    \draw[scale=1, domain=-\d:\d, smooth, variable=\x, red, dotted,thick] plot ({sqrt(.75)*cos(\ta)*(pow(2,\x)+pow(2,-\x))/2)},{sqrt(.75)*sin(\ta)*(pow(2,\x)-pow(2,-\x))/2} );
    \draw[scale=1, domain=-\d:\d, smooth, variable=\x, red,dotted, thick] plot (-{sqrt(.75)*cos(\ta)*(pow(2,\x)+pow(2,-\x))/2)},{sqrt(.75)*sin(\ta)*(pow(2,\x)-pow(2,-\x))/2} );

      \end{scope}
    \end{tikzpicture}
    \caption{The small cavity along with two confocal hyperbole. Eigenfunctions of the ellipse associated to trajectories tangent to the red (dotted) hyperbola localise inside the cavity and hence produce good quasimodes. Those associated to trajectories tangent to the blue (dashed) hyperbola do not localise in the cavity and hence do not produce good quasimodes.}\label{fig:integrableTori}
\end{figure}
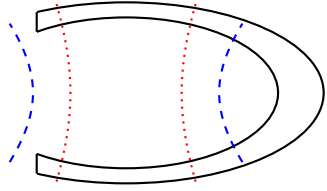

\section{Numerical experiments.}
\label{sec:experiments}

\subsection{Description of the set-up used for the experiments}

The BIEs \eqref{eq:direct_combined_dir}, \eqref{eq:direct_combined_neu}, and \eqref{eq:direct_combined_neu2}, involving the operators 
$A_k'$, $B_k$, and $B_{k, {\rm reg}}$ respectively, were discretised using the BEM with ten points per wavelength (as described in \S\ref{sec:statement:BEM}).
The eigenvalues and singular values of the resulting Galerkin matrices preconditioned with the mass matrix in 2-d were computed 
using the library BemTool\footnote{\url{https://github.com/xclaeys/BemTool}} and LAPACK~\cite{Anderson1999}. 
The largest matrices were around \(4,500\times 4,500\), and no distributed memory parallelisation was used. 

The results about the Galerkin error, GMRES residual, and GMRES error were obtained using the libraries PETSc~\cite{petsc-efficient,petsc-user-ref}, BemTool, Htool\footnote{\url{https://github.com/htool-ddm/htool}}, and SuperLU\_DIST~\cite{lidemmel03} via the software FreeFEM~\cite{Hecht2012}. Note that no compression was used and the largest matrices were around \(4,500\times 4,500\) in 2-d, and \(450,000\times 450,000\) in 3-d.

In some of the figures we plot best-fit lines; these are computed with the nonlinear least-squares Marquardt-Levenberg algorithm (the basis of the `fit' command in gnuplot). 

\subsection{Experiments about F3(a) and F3(b) in 2-d}
\begin{figure}[h!]
\begin{subfigure}{\textwidth}
    \centering
    \includegraphics[width=0.9\textwidth]{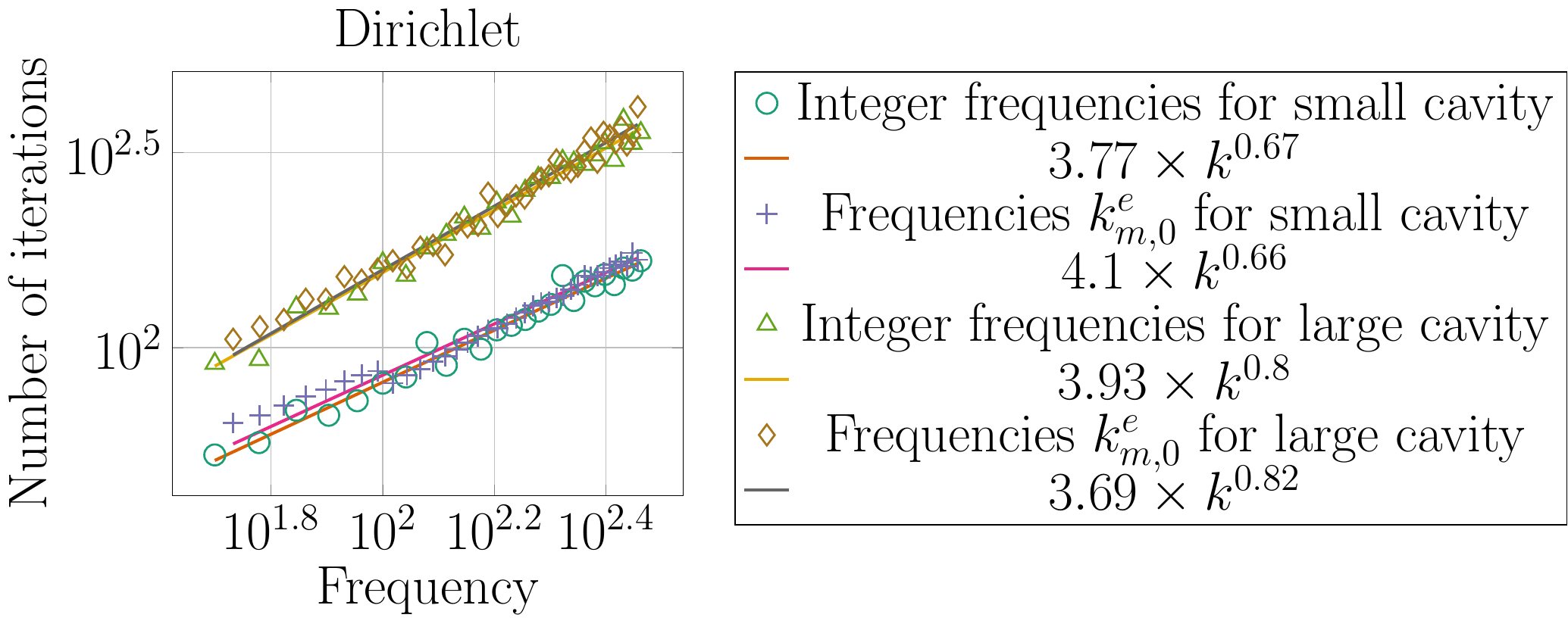}
    \caption{Number of iterations for the Dirichlet BIE involving $A_k'$}
    \label{fig:iterations_main_dir}
    \end{subfigure}
%
\begin{subfigure}{\textwidth}
    \centering
    \includegraphics[width=0.9\textwidth]{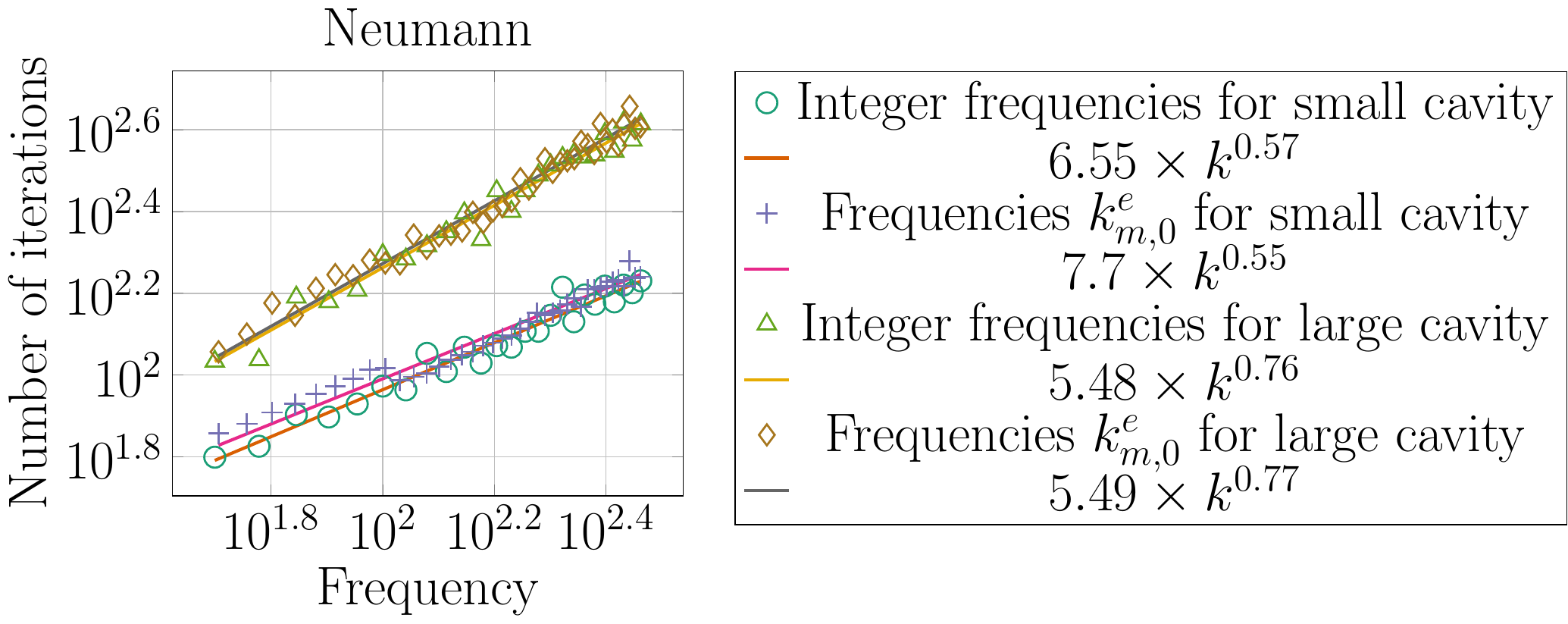}
    \caption{Number of iterations for the Neumann BIE involving $B_k$}
    \label{fig:iterations_main_neu}
\end{subfigure}
\begin{subfigure}{\textwidth}
    \centering
    \includegraphics[width=0.9\textwidth]{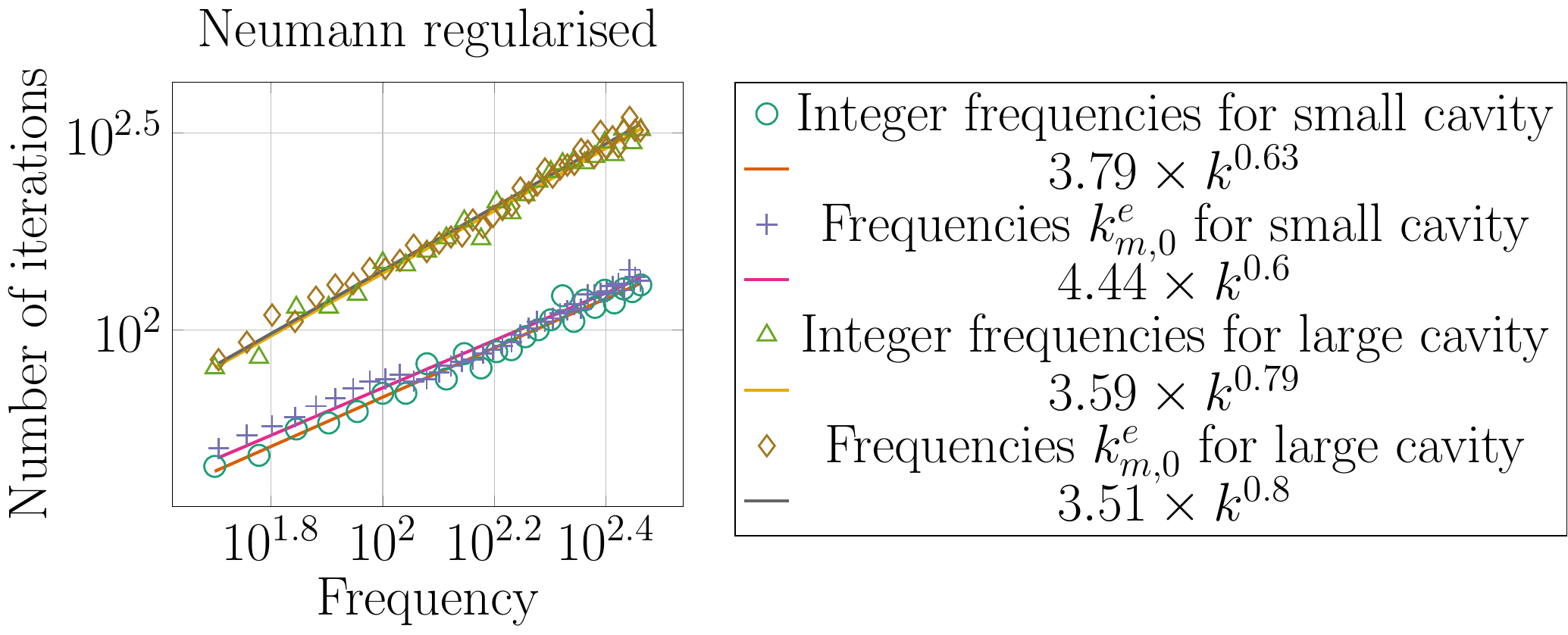}
    \caption{Number of iterations for the regularised Neumann BIE involving $B_{k, {\rm reg}}$}
    \label{fig:iterations_main_neu2}
    \end{subfigure}
    \caption{Number of iterations for the BIEs in 2-d with \(\theta=4\pi/10\) and $\Oi$ the small or large cavity}\label{fig:iterations_main}
\end{figure}

Figure \ref{fig:iterations_main} shows experiments about F3(a) ($k$-dependence of number of iterations) and F3(b) (dependence of number of iterations on the cavity size) for the Dirichlet BIE \eqref{eq:direct_combined_dir}, the Neumann BIE \eqref{eq:direct_combined_neu}, and the regularised Neumann BIE \eqref{eq:direct_combined_neu2}, through both integer frequencies and $k=k^e_{m,0}$. We saw a subset of these for the Dirichlet BIE in Figures \ref{fig:comparison_spectrum_it} and \ref{fig:comparison_size_cavit}.

The key point is that the growth of the number of iterations is the same through both sets of frequencies and for all three BIEs -- around $O(k^{0.8})$ for the large cavity and around $O(k^{0.6})$ for the small cavity -- even though 

(i) $B_{k}$ has a very different distribution of eigenvalues to $A'_k$ and $B_{k, {\rm reg}}$ as shown in \S\ref{sec:Nevalue} below, and 

(ii) the $k$-dependence of the norms of $B_k$ and $B_{k, {\rm reg}}$ is very different from that of $A_k'$ -- see \S\ref{sec:numbounds} below.

\subsection{Experiments about F3(a) in 3-d}\label{sec:3-dexpts}

\begin{figure}[h!]
    \centering
    \includegraphics[width=0.9\textwidth]{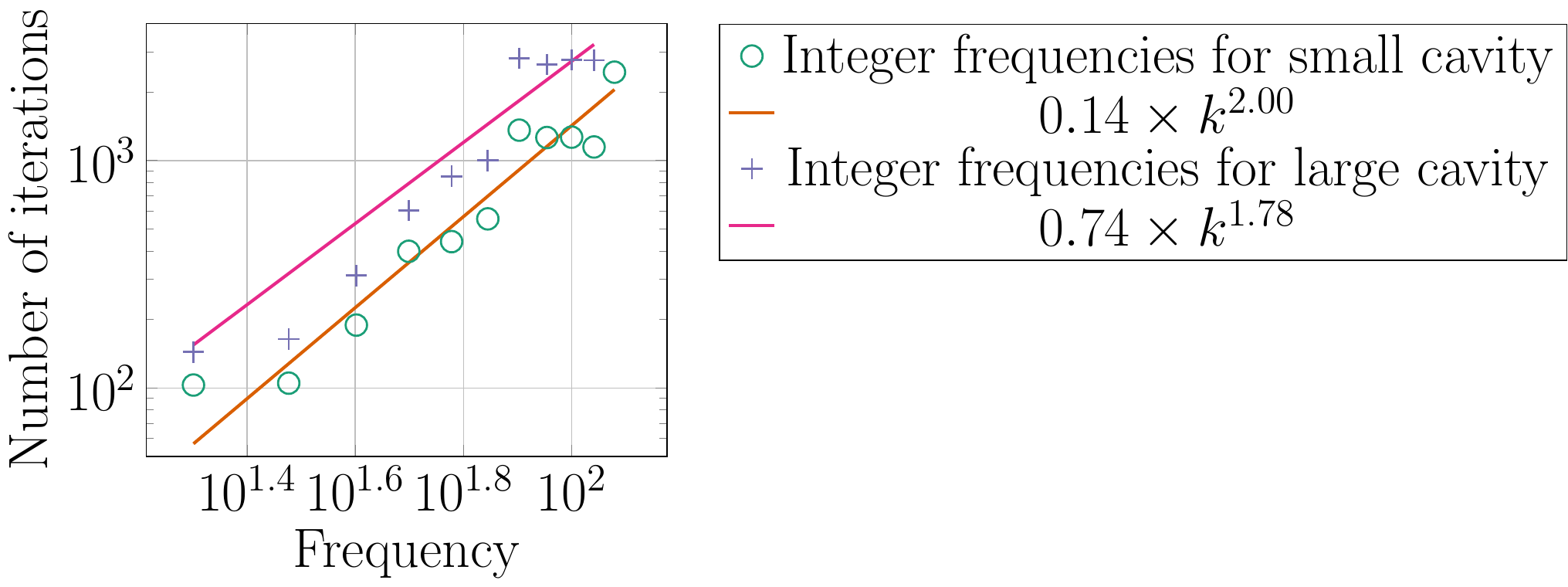}
    \caption{Number of iterations for the Dirichlet BIE involving $A_k'$ with \(\theta=4\pi/10\) and $\Oi$ the small and large cavities in 3-d}\label{fig:3d}
\end{figure}

\begin{figure}[h!]
    \centering
    \includegraphics[width=0.8\textwidth]{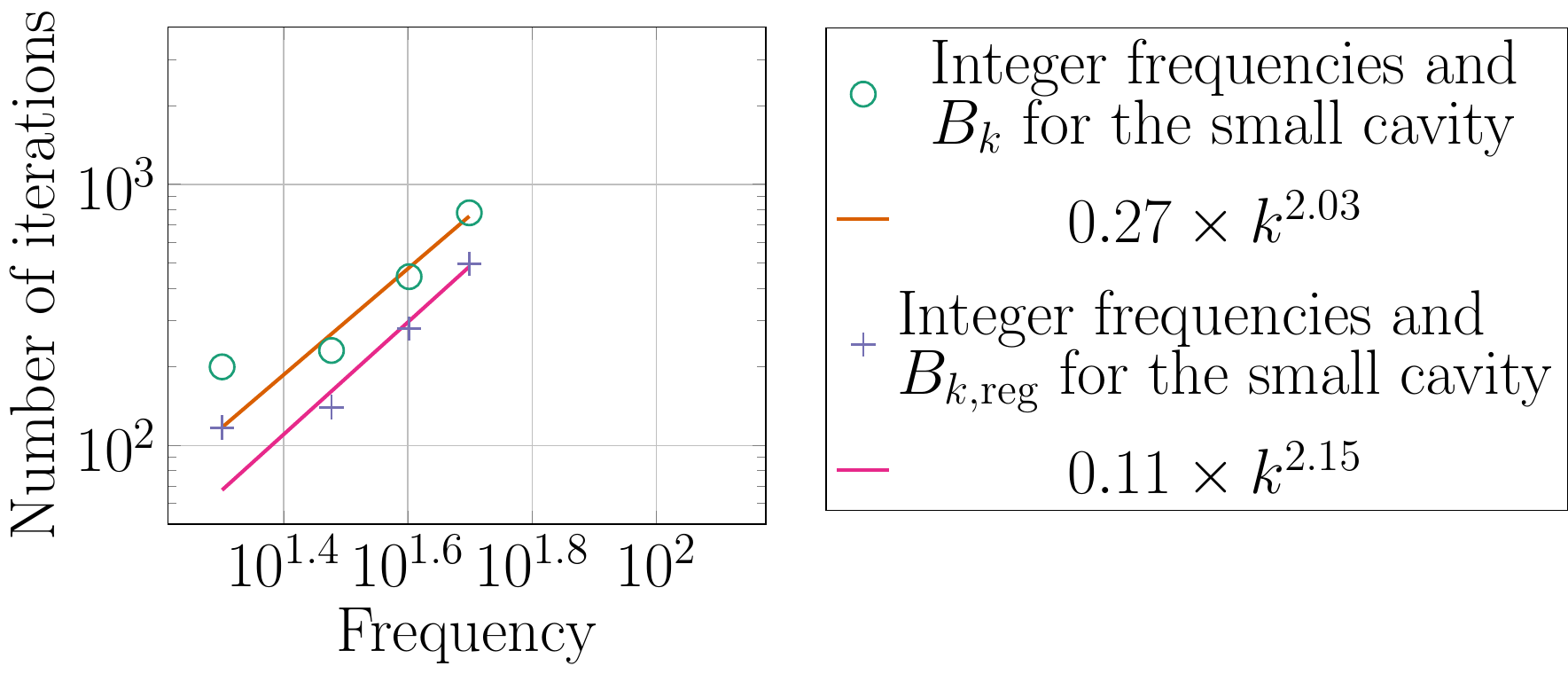}
    \caption{Number of iterations for the Neumann BIEs involving $B_k$ and $B_{k,{\rm reg}}$ with \(\theta=4\pi/10\) and $\Oi$ the small cavity in 3-d}\label{fig:3d_neu}
\end{figure}

Figure \ref{fig:3d} shows the number of GMRES iterations when the Dirichlet BIE \eqref{eq:direct_combined_dir} is solved with $\Oi$ the 3-d analogues of the small and large cavities (as described in \S\ref{sec:smalllargegeometries})
Similarly, Figure \ref{fig:3d_neu} 
(plotted on the same scale as Figure \ref{fig:3d} for ease of comparison) shows the number of GMRES iterations when the Neumann BIEs \eqref{eq:direct_combined_neu} and 
\eqref{eq:direct_combined_neu2}
are solved with $\Oi$ the 3-d analogue of the small cavity. For the Neumann BIEs we were unable to go up to as high a frequency as for the Dirichlet BIE because of memory issues.

The fact that in all cases the number of iterations grows roughly like $O(k^2)$, in contrast to roughly like $O(k)$ in 2-d, is consistent with the factors of $k^{d-1}$ appearing in Theorem \ref{thm:main2}. 

These 3-d experiments only consider the number of GMRES iterations at integer frequencies. Our definitions of the 3-d small and large cavities are such that the 3-d analogue of the ellipse $E$ \eqref{eq:ellipse} is now a prolate spheroid. 
Since the Laplacian is separable in prolate spheroidal coordinates (see, e.g., \cite[\S30.13]{Di:21}), 
the eigenvalues of the Laplacian, and hence the corresponding quasimode frequencies (analogous to $k_{m,n}^{e/o}$), can be computed in a similar way to those in 2-d (as described in Appendix \ref{sec:Mathieu}), but we have not pursued this here.

\subsection{Plots of the eigenvalues of $B_{k}$ and $B_{k, {\rm reg}}$ in 2-d.}\label{sec:Nevalue} 

\begin{figure}[h!]
    \centering
    \includegraphics[width=0.9\textwidth]{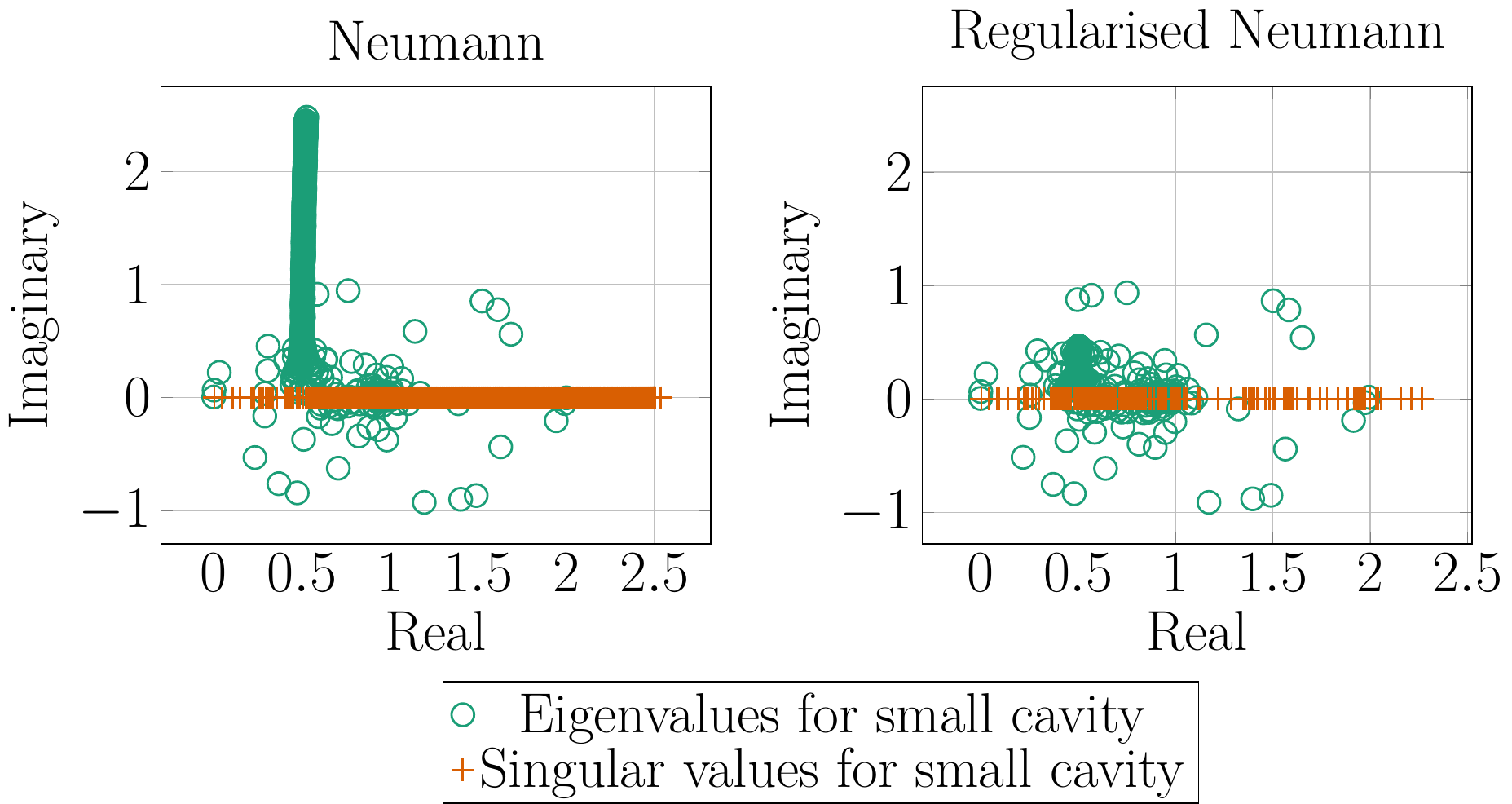}
    \caption{The eigenvalues and singular values for the discretisations of the Neumann BIEs involving $B_k$ (``Neumann'') and $B_{k, {\rm reg}}$ (``Regularised Neumann'') for the small cavity at \(k=k^e_{0,0}\).}\label{fig:eigenvalues_Neumann}
\end{figure}

Figure \ref{fig:eigenvalues_Neumann} plots the eigenvalues and singular values of the discretisations of $(2/\ri) B_k$ (``Neumann'') and $(2/\ri)B_{k, {\rm reg}}$ (``regularised Neumann'') for $\Oi$ the small cavity at $k=k^e_{0,0}$;
the division by $\ri$ rotates the eigenvalues so that they are in the same half-plane as the eigenvalues of $A_k'$. The accumulation point in the right-hand plot of Figure \ref{fig:eigenvalues_Neumann} is at $(1+\ri)/2$, which is consistent with $B_{k, {\rm reg}}$ being a compact perturbation of $\ri (1+\ri)/4$ when $\Gamma$ is $C^1$ by \eqref{eq:Breg_compact_pert}.

These plots show the effect of regularising the hypersingular operator $H_k$. 
Indeed, the discretisation of $B_k$ contains a vertical ``tail'' of eigenvalues that $\sim 1/h$ as $h\to 0$ for fixed $k$ 
\cite[Exercise 4.5.2]{SaSc:11}, \cite[Lemma 12.9]{St:08}. Since $H_k: H^1(\Gamma)\to \LtG$ and $\LtG\to H^{-1}(\Gamma)$, $H_k$ restricted to low frequencies maps to $\LtG$,
and the eigenvalues in the tail are associated with high frequencies. Furthermore, the direction of the tail can be explained by the symbol of $H_k$ as a pseudodifferential operator on high frequencies; see \cite[\S4.3]{Ga:19}, \cite[\S4.1]{GaMaSp:21N}.
This tail of eigenvalues is not present for $B_{k, {\rm reg}}$ because $S_{\ri k}H_k:\LtG\to\LtG$.

\subsection{Experiments about quantities in the bound of Theorem \ref{thm:main1}.}\label{sec:numbounds}

\begin{figure}[h!]
    \centering
    \includegraphics[width=0.9\textwidth]{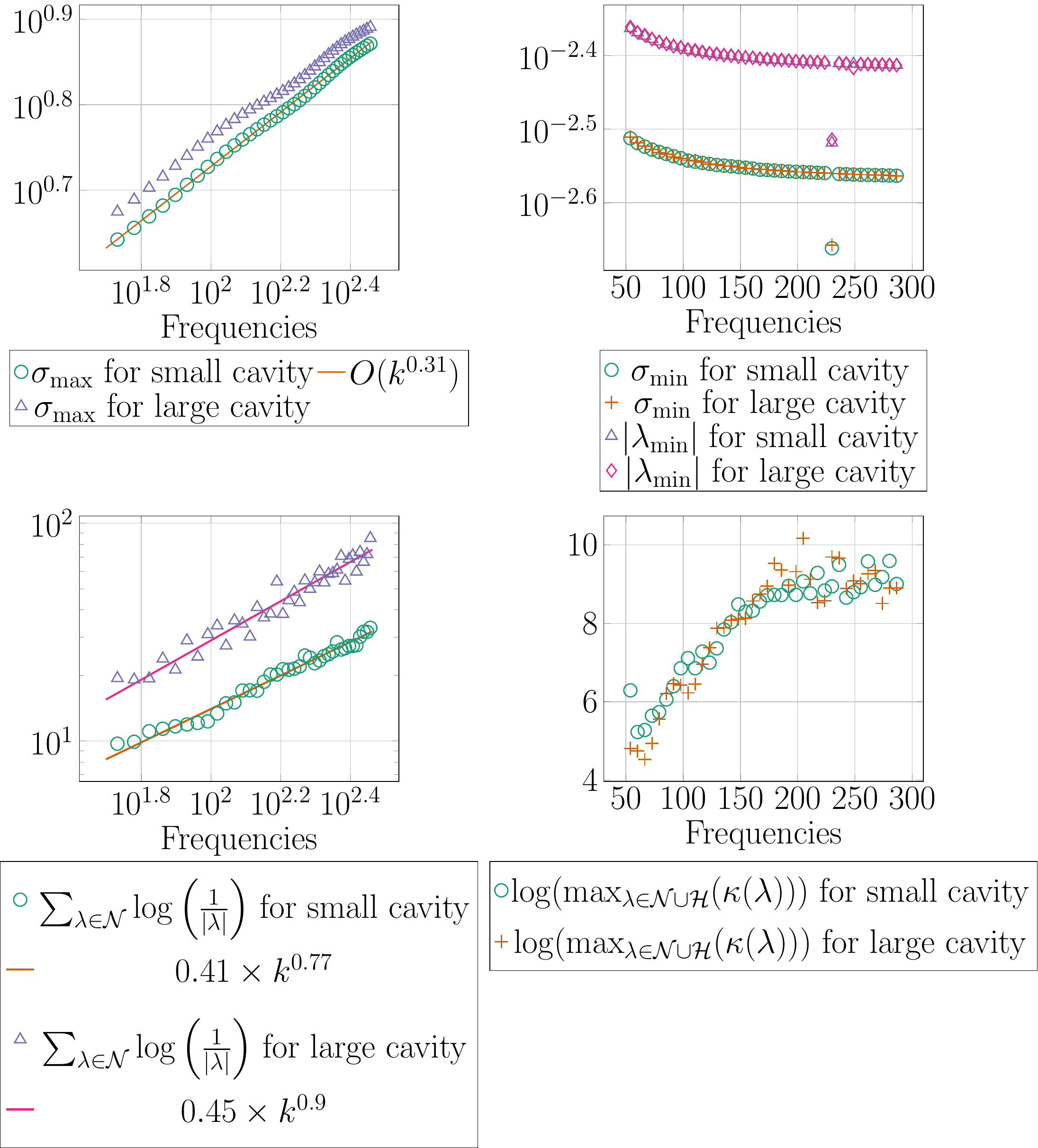}
    \caption{Properties of the discretisation of the Dirichlet BIE involving $A_k'$ for  \(k=k^e_{m,0}\)
    and $\cN= [-0.1,0.1]\times[-0.6,0.6]$.   
    }\label{fig:Dirichlet}
\end{figure}

\begin{figure}[h!]
    \centering
    \includegraphics[width=0.9\textwidth]{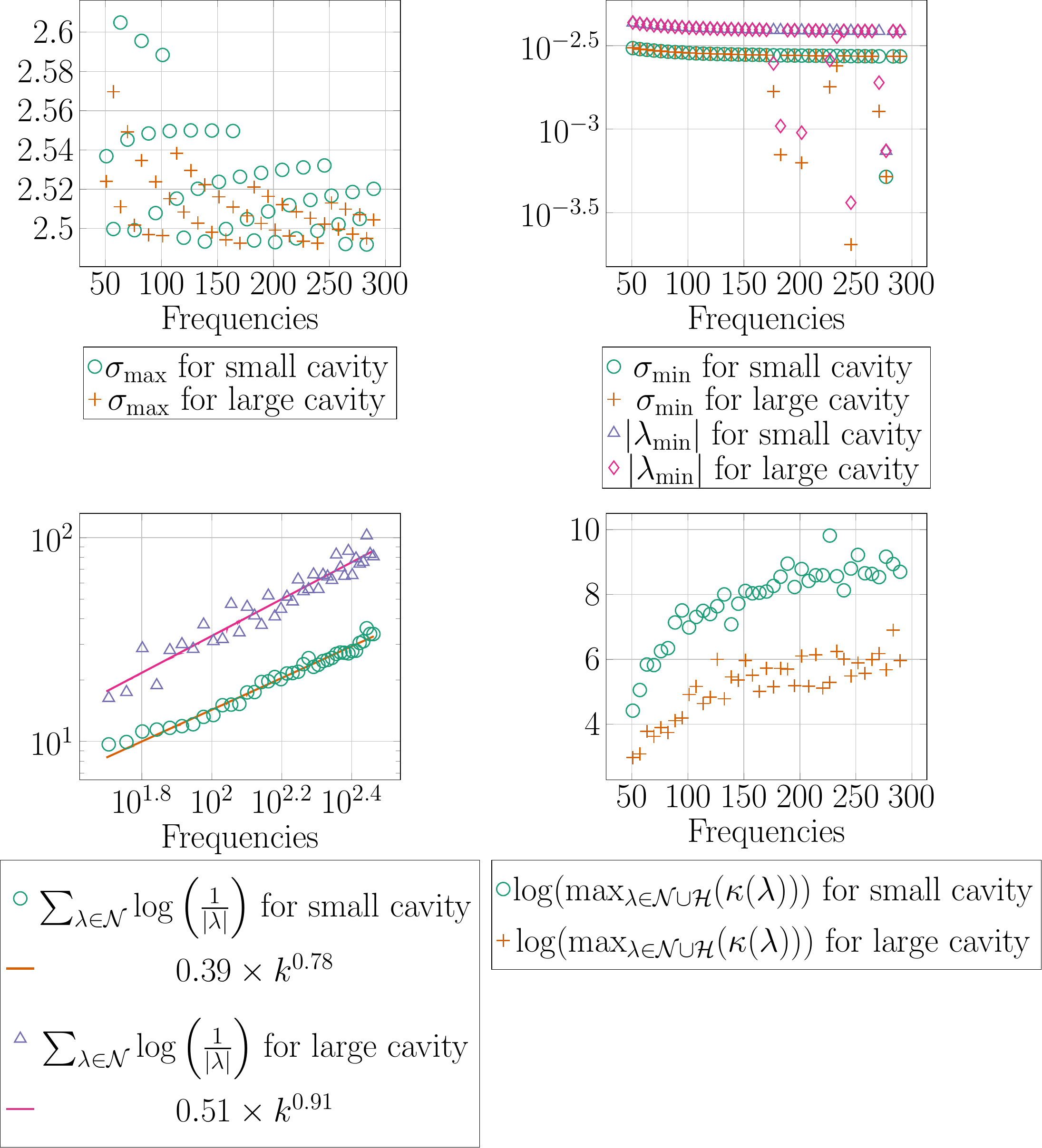}
    \caption{Properties of the discretisation of the Neumann BIE involving $B_k$ for  \(k=k^e_{m,0}\)
        and $\cN= [-0.1,0.1]\times[-0.6,0.6]$.
    }\label{fig:Neumann}
\end{figure}

\begin{figure}[h!]
    \centering
    \includegraphics[width=0.9\textwidth]{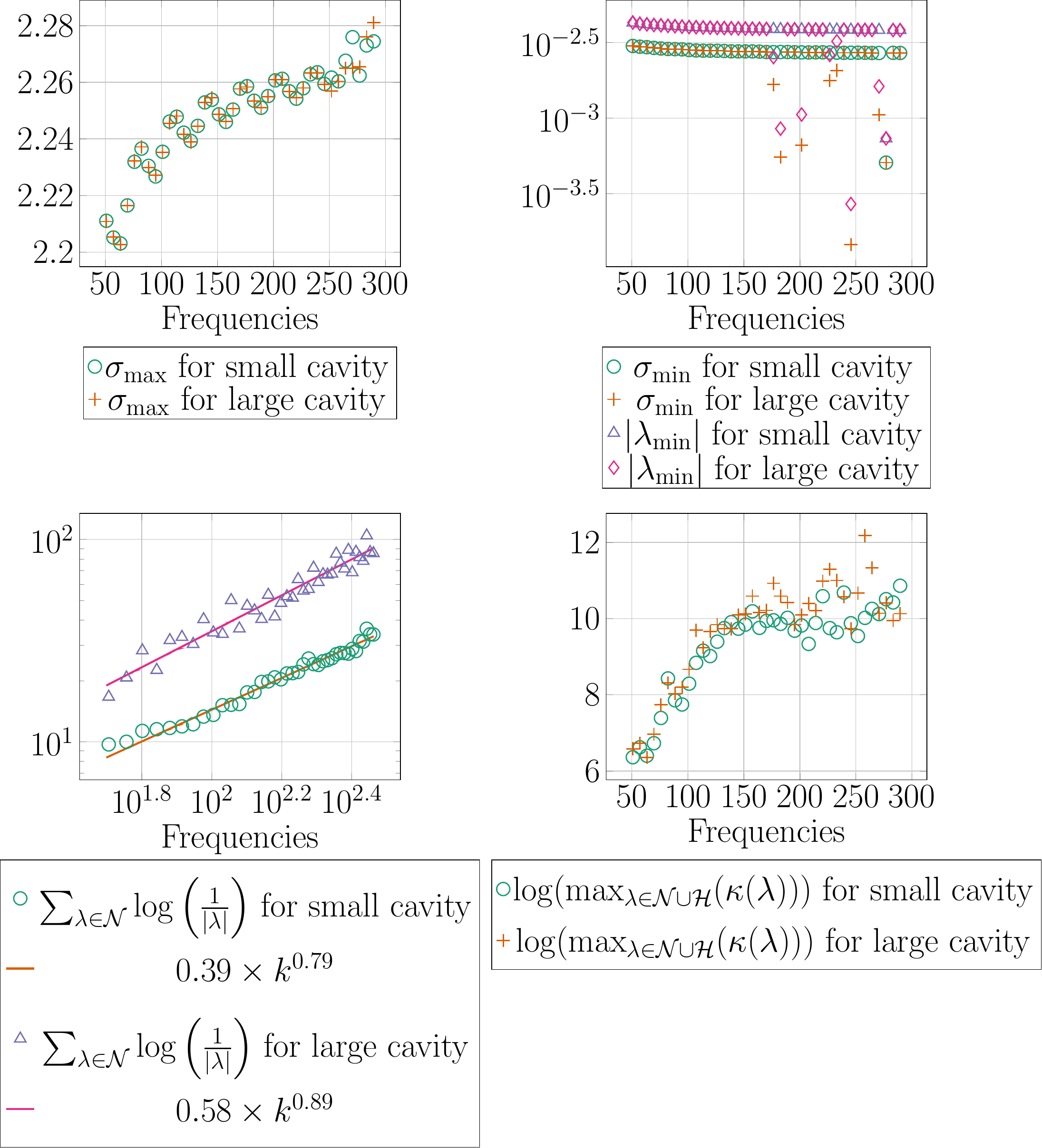}
    \caption{Properties of the discretisation of the regularised Neumann BIE involving $B_{k, {\rm reg}}$ for \(k=k^e_{m,0}\)
    and $\cN= [-0.1,0.1]\times[-0.6,0.6]$.
    }\label{fig:reg_Neumann}
\end{figure}

Figures \ref{fig:Dirichlet}, \ref{fig:Neumann}, and \ref{fig:reg_Neumann} plot the following quantities for discretisations of each of the operators $A_k', B_k,$ and $B_{k, {\rm reg}}$ through  \(k=k^e_{m,0}\).
\bit
\item Top-left plot: the maximum singular value of $\bfM^{-1}\bfA$, i.e., $\|\bfM^{-1}\bfA\|_2$, where $\bfA$ is the respective Galerkin matrix, 
\item Top-right plot: the minimum singular value of $\bfM^{-1}\bfA$, i.e., $\|(\bfM^{-1}\bfA)^{-1}\|_2$, and the eigenvalue of $\bfM^{-1}\bfA$ with the smallest modulus.
\item Bottom-left plot: the quantity
\beq\label{eq:keyquant}
\sum_{\lambda \in \cN}\log \left(\frac{1}{|\lambda|}\right)
\eeq
where $\cN = [-0.1,0.1]\times [-0.6,0.6]$ (as in Figures \ref{fig:outliers} and \ref{fig:counting_function}).
\item Bottom-right plot: $\log \big(\max_{\lambda}(\kappa(\lambda))\big)$, where $\kappa(\lambda)$ is the eigenvalue condition number defined by \eqref{eq:eigenvaluecondition}.
\eit

\paragraph{Regarding the top-left plots:} recalling that $\|\bfM^{-1}\matrixD\|_2$ approximates $\|A_k'\|_{\LtGt}$ (see Lemma \ref{lem:discon} and Assumption A0), we 
 used this information about $\|A_k'\|_{\LtGt}$ and $\|B_{k, {\rm reg}}\|_{\LtGt}$ in our discussion in \S\ref{sec:discussion} about the sharpness of Theorem \ref{thm:main2}.
 
Figure \ref{fig:Dirichlet} shows  $\|A_k'\|_{\LtGt}$ growing approximately like $k^{1/3}$. 
How the geometry of $\Oi$ affects the $k$-dependence of $\|A_k'\|_{\LtGt}$
is now well-understood thanks to the results of \cite{ChGrLaLi:09}, \cite[Appendix A]{HaTa:15}, \cite{GaSm:15},  \cite[Chapter 4]{Ga:19}, and \cite{GaSp:19}.
In fact, these results show that the $k$-dependence of $\|A_k'\|_{\LtGt}$ as $k\tendi$ is dominated by the $k$-dependence of $k\|S_k\|_{\LtGt}$, and this $\sim k^{1/3}$ on curved parts of $\Gamma$ and $\sim k^{1/2}$ on flat parts, with the omitted constants dependent on the surface measure of these parts of the boundary. For both the small and large cavities, the surface measure of the flat parts of $\Oi$ is much smaller than the surface measure of the curved parts of $\Oi$ (see Figure \ref{fig:geometries}), and this is the reason why we only see the $k^{1/3}$ growth for the range $k\in (50,290)$ in  Figure \ref{fig:Dirichlet}.

Similarly, Figure \ref{fig:reg_Neumann} shows $\|B_{k, {\rm reg}}\|_{\LtGt}$ being essentially constant for the range of $k$ considered, although,
 at least in 2-d, $\|B_{k, {\rm reg}}\|_{\LtGt} \gtrsim k^{1/4}$ for large enough $k$; indeed, \cite[Theorem 4.6]{ChGrLaLi:09} shows that $\|D_k'\|_{\LtGt}\gtrsim k^{1/4}$ for a certain class of 2-d domains (to see that the elliptic cavity falls in this class, take the points $x_1$ and $x_2$ in the statement of \cite[Theorem 4.6]{ChGrLaLi:09} to lie on one of the flat ends of the cavity, with $x_2$ in the middle of this end, and $x_1$ at one of the corners)
and \cite[Theorems 4.6 and 4.8]{GaMaSp:21N} show that $\|S_{\ri k }H_k\|_{\LtGt}\lesssim (\log k)^{3/2}$.
 
\paragraph{Regarding the top-right plots:} these show both 

(i) the feature F2, i.e.~that while the norms of the inverses of the boundary-integral operators grow exponentially through $k^e_{m,0}$, 
and thus the smallest singular values should decrease exponentially, this growth/decay stagnates, and 

(ii) that the smallest eigenvalue modulus is very close the smallest singular value, 
giving indirect evidence for Assumption A2, i.e., that at $k^e_{m,0}$ (for large enough $m$), the matrix has both a small singular value and a near-zero eigenvalue.

\paragraph{Regarding the bottom-left plots:}
these plots show the quantity \eqref{eq:keyquant} growing differently for the small and large cavities, and we used this information in our discussion in \S\ref{sec:discussion} about the sharpness of Theorem \ref{thm:main2}.

\paragraph{Regarding the bottom-right plots:} these verify Assumption A3, i.e., that the maximum eigenvalue condition number does not grow exponentially with $k$ (at least for the range of $k$ considered, i.e., $k\in (50, 290)$).

\section{Conclusions}\label{sec:conclusions}

In \S\ref{sec:features}, we stated that the main goals of this paper were to explain the feature F3(a) (i.e., why the number of GMRES iterations grows algebraically with $k$, with no worse growth through quasimode frequencies)
and, to a certain extent,
F3(b) (i.e., why the number of iterations depends on whether $\Oi$ is the small or large cavity).

Theorem \ref{thm:main2} addresses the $k$-dependence of the number of iterations (i.e., F3(a)). Although the bound \eqref{eq:mlowerbound2} in Theorem \ref{thm:main2} 
does not directly distinguish between the small and large cavities, and hence does not explain F3(b), the coefficient of the highest-order terms in the bound \eqref{eq:mlowerbound2} depends on the number of the near-zero eigenvalues (via the constant $\CWeyl$), and the arguments in \S\ref{sec:smalllarge} then explain heuristically the difference between this number for the small and large cavities.

\

For future investigations of GMRES applied to Helmholtz trapping scenarios, we have the following conclusions/messages.

\paragraph{The difference between $k= k_j$ and $k\neq k_j$ (where $k_j$ is a frequency in a quasimode) is not important.}

We saw from Figures \ref{fig:comparison_spectrum_it}
and \ref{fig:iterations_main} that the growth in the number of GMRES iterations did not depend on whether the frequency was in a quasimode (in contrast to the condition number, which does depend strongly on this).

\paragraph{The important quantities are (i) the rate of growth of the cluster, and (ii) the number of near-zero eigenvalues (governed by the number of quasimode frequencies).}

\

Regarding (i): the norm is a proxy for this, but comparing the experiments for $A_k'$ (the BIE \eqref{eq:direct_combined_dir} for the Dirichlet problem) in Figure \ref{fig:Dirichlet} and $B_{k, {\rm reg}}$ (the regularised BIE \eqref{eq:direct_combined_neu2} for the Neumann problem) 
 in Figure \ref{fig:reg_Neumann} we see one norm growing with $k$ (i.e.~$\|A_k'\|_\LtGt$), the other norm remaining constant (i.e.~$\|B_{k, {\rm reg}}\|_\LtGt$), but the 
number of GMRES iterations for both growing at the same rate -- see Figure \ref{fig:iterations_main}.

Regarding (ii): the arguments in \S\ref{sec:Weyl_new} show that this number is governed by the Weyl law, and hence depends on dimension. 
We highlight that, once the frequency is high enough, the density of these near-zero eigenvalues becomes too high for them to be considered as true ``outliers'' -- see Figure \ref{fig:example_spectrum} -- but the bounds of Theorems \ref{thm:main1} and \ref{thm:main2} still hold.

We advocate that these two quantities (i) and (ii) should play the role for Helmholtz trapping problems that the condition number plays in \emph{both} understanding the behaviour of the conjugate gradient method (CG) \emph{and} designing preconditioners for symmetric, positive-definite matrices. Indeed, if \(\genmatrix\) is symmetric positive-definite, it is well-known that 
\begin{align}\label{eq:CG}
    \dfrac{\lVert \bfr_m (\genmatrix,\bfb,\bfx_0) \rVert_{\genmatrix}}{\lVert \bfr_0 (\genmatrix,\bfb,\bfx_0) \rVert_{\genmatrix}} \leq2  \left(\sqrt{\dfrac{\kappa(\genmatrix)-1}{\kappa(\genmatrix)+1}}\right)^m,
\end{align}
where \(\lVert \cdot \rVert_{\genmatrix}\) denotes the norm induced by \(\genmatrix\). 
In a similar way to the bounds on GMRES, the bound \eqref{eq:CG} can overestimate the number of iterations because it does not take into account either the right-hand side, or the fact that CG (like GMRES) can have a superlinear convergence, which depends on the particular distribution of eigenvalues. 
Therefore, CG does not converge with the same speed for all matrices with the same maximum/minimum eigenvalues (and therefore for which the right-hand side of \eqref{eq:CG} is the same).

Despite these drawbacks, the bound \eqref{eq:CG} is useful in at least the two following ways:

(a) It indicates that preconditioners should be designed with the goal of decreasing the condition number and guarantees a reduction in the number of iterations if the resulting condition number is sufficiently small.

(b) If one can show that the condition number of $\genmatrix$ is independent of a certain parameter, it shows that the number of GMRES iterations to achieve a prescribed tolerance can be bounded independently of this parameter; this fact is used in, e.g., domain-decompositions methods where the number of iterations must be independent of the number of subdomains for the method to scale.

We advocate the use of the quantities (i) and (ii) above in a similar way. In particular the bounds of Theorems \ref{thm:main1} and \ref{thm:main2} show that a sufficient condition for a preconditioner to be robust at high frequency is for it to mitigate against the quantities (i) and (ii) growing with frequency.

\appendix

\section{Definitions of layer potentials and boundary-integral operators}

The single-layer and double-layer potentials, $\cS_k$ and $\cD_k$ respectively, are defined for $\phi\in L^1(\Gamma)$ by 
\begin{align}\label{eq:SLP}
    \calS_k \varphi (\bx) &= \int_{\Gamma} \Phi_k (\bx,\by) \varphi (\by) \dif s (\by) \quad\tfa \bx \in \bbR^d \setminus \Gamma, \quad\tand \\
    \calD_k \varphi (\bx) &= \int_{\Gamma} \dfrac{\partial \Phi_k (\bx,\by)}{\partial n(\by)} \varphi (\by) \dif s (\by) \quad \tfa \bx \in \bbR^d \setminus \Gamma,\nonumber
\end{align}
where the fundamental solution \(\Phi_k\) is defined by 
\begin{align*}
    \Phi_k (\bx,\by) := 
            \dfrac{\ri}{4} H^{(1)}_0(k \lvert \bx - \by \rvert ),\quad d=2,\qquad
:=            \dfrac{\re^{\ri k \lvert \bx-\by \rvert}}{4 \pi \lvert \bx - \by \rvert},&\quad d=3,
\end{align*}
where \(H^{(1)}_0\) is the Hankel function of the first kind and order zero. 
If $u$ is the solution to the scattering problem \eqref{eq:Helmholtz}, then Green's second identity implies that 
\begin{align}\label{eq:Green}
    u(\bx) = u^I(\bx) +\big(\calD_k u\big)(\bx)- \big(\calS_k \partial_{n} u\big) (\bx) \quad \tfor \bx \in \Oe,
\end{align}
see, e.g., \cite[Theorem 2.43]{ChGrLaSp:12}.

The single-layer, adjoint-double-layer, double-layer, and hypersingular operators are defined for $\phi\in \LtG$ and $\psi\in H^1(\Gamma)$ by 
\begin{align}\label{eq:SD'}
&S_k \phi(\bx) := \int_\bound \Phi_k(\bx,\by) \phi(\by)\,\rd s(\by), \qquad
D'_k \phi(\bx) := \int_\bound \frac{\partial \Phi_k(\bx,\by)}{\partial n(\bx)}  \phi(\by)\,\rd s(\by),\\
&D_k \phi(\bx) := \int_\bound \frac{\partial \Phi_k(\bx,\by)}{\partial n(\by)}  \phi(\by)\,\rd s(\by), 
\quad 
H_k \psi(\bx) := \pdiff{}{n(\bx)} \int_\bound \frac{\partial \Phi_k(\bx,\by)}{\partial n(\bx)}  \psi(\by)\,\rd s(\by),
 \label{eq:DH}
\end{align}
for $\bx \in \Gamma$. When $\Gamma$ is Lipschitz, the integrals defining $D_k$ and $D_k'$ must be understood as Cauchy principal values (see, e.g., \cite[Equation 2.33]{ChGrLaSp:12}), and the integral defining $H_k$ is understood as a non-tangential limit (see, e.g., \cite[Equation 2.36]{ChGrLaSp:12}), but we do not need the details of these definitions in this paper.

\section{Bounds on the Galerkin matrix in terms of the continuous operator.}

\ble
\label{lem:discon}
Let $V_n\subset \LtG$ be a finite-dimensional space with real basis $\{\phi_j\}_{j=1}^n$. Let $\bfM$ be defined by \eqref{eq:massmatrix}.
Given $A:\LtGt$, let $\bfA$ be defined by the first equation in \eqref{eq:Galerkinmatrix} (with $A'_k$ replaced by $A$).
Let $P_h:\LtG\rightarrow V_h$ be the orthogonal projection, and let 
\beqs
\widetilde{A}:= P_h A|_{V_h}.
\eeqs

(i) 
\beq\label{eq:disconresult1}
(\cond(\bfM))^{-1/2}\, \big\|\widetilde{A}\big\|_{\LtGt} \leq \N{\bfM^{-1}\bfA}_2 \leq (\cond(\bfM))^{1/2} \,\big\|\widetilde{A}\big\|_{\LtGt}
\eeq
where $\cond(\bfM):= \|\bfM\|_2 \|\bfM^{-1}\|_2$, and if $(\bfM^{-1}\bfA)^{-1}$ exists, then
\beq\label{eq:disconresult2}
 (\cond(\bfM))^{-1/2}\, \big\|\widetilde{A}^{-1}\big\|_{\LtGt}\leq \N{(\bfM^{-1}\bfA)^{-1}}_2 \leq  (\cond(\bfM))^{1/2}
\, \big\|\widetilde{A}^{-1}\big\|_{\LtGt}.
\eeq

(ii) If $P_h \phi \rightarrow \phi$ as $h\tendo$ for all $\phi\in\LtG$, then
\beq\label{eq:disconresult3}
\big\|\widetilde{A}\big\|_{\LtGt} \rightarrow \N{A}_{\LtGt} \quad\tas h\tendo;
\eeq
if, in addition, $A = a I + K$, where $a\neq 0$ and $K$ is compact, then 
\beq\label{eq:disconresult4}
\big\|\widetilde{A}^{-1}\big\|_{\LtGt} \rightarrow \big\|A^{-1}\big\|_{\LtGt} \quad\tas h\tendo.
\eeq
\ele

\bre
For standard BEM spaces, $\cond(\bfM)$ is bounded independently of $h$; see \cite[Theorem 4.4.7 and Remark 4.5.3]{SaSc:11} and \cite[Corollary 10.6]{St:08}.
\ere

\bre\label{rem:compact}
If $\Gamma$ is $C^1$, then both $A_k'$ and $B_{k, {\rm reg}}$ can be written as $a I + K$, where $a\neq 0$ and $K$ is compact. 
For $A_k'$ this follows since $S_k:\LtG\rightarrow \HoG$ when $\Gamma$ is Lipschitz (see, e.g., \cite[Theorem 2.17]{ChGrLaSp:12}) and $D_k'$ is compact on $\LtG$ when $\Gamma$ is $C^1$ by \cite[Theorem 1.2]{FaJoRi:78}.
\footnote{It was recently shown in \cite{ChSp:21-1}, however, that there exist star-shaped Lipschitz polyhedra such that $A_k'$ cannot be written as the sum of a coercive operator and a compact operator for any $k>0$, and thus cannot be written as $a I + K$, where $a\neq 0$ and $K$ is compact; see \cite[Corollary 1.6]{ChSp:21-1}.}
For $B_{k, {\rm reg}}$, the Calder\'on relations (see, e.g., \cite[Equation 2.56]{ChGrLaSp:12}) imply that
\begin{align}\label{eq:Breg_compact_pert}
    B_{k,{\rm reg}}&= \dfrac{\ri}{2} \left(\dfrac{I}{2}-D_k\right) +  S_{\ri k}\big(H_k- H_{\ri k}\big) + S_{\ri k} H_{\ri k}= \dfrac{\ri}{2} \left(\dfrac{I}{2}-D_k\right) +  S_{\ri k}\big(H_k- H_{\ri k}\big) + \left(-\dfrac{I}{4}+D_{\ri k}^2\right).
\end{align}
By bounds on the fundamental solution $\Phi_k$ appearing in, e.g., \cite[Equation 2.25]{ChGrLaSp:12}, the kernel of the integral operator $H_k - H_{\ri k}$ is weakly singular, and thus the operator is compact on $\LtG$ by, e.g., the combination of \cite[Part 3 of the theorem on Page 49]{Pr:91} and Young's inequality. 
\ere

\bpf[Proof of Lemma \ref{lem:discon}]

Part (ii) is proved in \cite{BeChGrLaLi:11}, with \eqref{eq:disconresult3} proved in \cite[Equation 3.3]{BeChGrLaLi:11} and 
 \eqref{eq:disconresult4} proved in \cite[Text below Equation 3.3]{BeChGrLaLi:11}.
 We note that similar results are contained in \cite[\S2.4]{Ki:10}.
 
For Part (i), let $\{\psi_j\}_{j=1}^n$ be an real orthonormal basis of $V_n$ so that 
\beq\label{eq:disconresult5}
\phi_i = \sum_{j=1}^n \big(\phi_i,\psi_j\big)_{\LtG} \psi_j.
\eeq
Let the matrices $\bfB$ and $\bfD$ be defined by 
\beqs
    (\bfB)_{i,j}:= \int_{\Gamma} \big(A\psi_j\big) (\bx)\psi_i (\bx) \dif s(\bx)
\quad\tand\quad
    (\bfD)_{i,j}:= \int_{\Gamma} \phi_i (\bx)\,\psi_j (\bx) \dif s(\bx);
\eeqs
i.e.~$\bfB$ is the Galerkin matrix with respect to the basis $\{\psi_j\}_{j=1}^n$, and $\bfD$ is a change of basis matrix from 
$\{\psi_j\}_{j=1}^n$ to $\{\phi_j\}_{j=1}^n$.
Since $\{\psi_j\}_{j=1}^n$ is orthonormal,
\beq\label{eq:disconresult5a}
\N{\bfB}_2 = \big\|\widetilde{A}\big\|_{\LtGt} \quad\tand\quad \N{\bfB^{-1}}_2 = \big\|\widetilde{A}^{-1}\big\|_{\LtGt}
\eeq
(provided that $\widetilde{A}^{-1}$ exists) by, e.g., \cite[Equation 3.2]{BeChGrLaLi:11}.

The definitions of $\bfB$ and $\bfD$, combined with \eqref{eq:disconresult5}, imply that $\bfM = \bfD \bfD^T$ and $\bfA = \bfD \bfB\bfD^T$
so that 
\beq\label{eq:disconresult6}
\bfM^{-1}\bfA = \bfD^{-T} \bfB \bfD^T \quad\tand\quad(\bfM^{-1}\bfA )^{-1}= \bfD^{-T} \bfB^{-1} \bfD^T.
\eeq
Since $\bfM$ is symmetric positive definite, the definition $\|\bfD\|_2 = \sqrt{\lambda_{\max}(\bfD^T \bfD)}$ implies that $\|\bfD\|_2 =
\|\bfD^T\|_2=
 \|\bfM\|_2^{1/2}$ and similarly $\|\bfD^{-1}\|_2 =
\|\bfD^{-T}\|_2=
 \|\bfM^{-1}\|_2^{1/2}$
The results \eqref{eq:disconresult1} and \eqref{eq:disconresult2} then follow from combining these results about norms with 
\eqref{eq:disconresult5a} and 
\eqref{eq:disconresult6}.
\epf

\bre[Bounds on norms without preconditioning by $\bfM^{-1}$.]
The proof of Lemma \ref{lem:discon} also implies that
\beq\label{eq:Galerkinmatrixbounds2}
\|\bfA\|_2 \leq \|\bfM\|_2 
\big\|\widetilde{A}\big\|_{\LtGt}
\quad\tand\quad
\big\|\bfA^{-1}\big\|_2 \leq \big\|\bfM^{-1}\big\|_2 
 \big\|\widetilde{A}^{-1}\big\|_{\LtGt}
\eeq
so that 
\beqs
\cond(\bfA) \leq \cond(\bfM) \cond\big(\widetilde{A}\big);
\eeqs
see, e.g., \cite[Equation 3.6.166]{At:97}.
\ere

\section{Numerical experiments about F1 and F4}\label{app:F1F4}

Recall from \S\ref{sec:features} the features F1 and F4: 
\bit
\item[F1] When the incoming plane wave enters the cavity, one needs a larger number of points per wavelength for accuracy of the Galerkin solutions 
than when the wave doesn't enter the cavity. 
\item[F4] The GMRES residual being small does not necessarily mean that the error is small, and the relative sizes of the residual and error depend on both $k$ and the direction of the plane wave.
\eit

\paragraph{Numerical experiments demonstrating F1.} These experiments involve the Galerkin solutions on four different meshes with $\Oi$ the small cavity. To create these meshes  
we start with a mesh with \(10\) points by wavelength, and then we refine splitting the segments in two, three, and six.
The Galerkin solutions are computed using LU factorisations computed by SuperLU~\cite{lidemmel03}. 
We use the Galerkin solution on the finest mesh as a proxy for the true solution (with this mesh denoted by \(\Gamma_{\text{ref}}\)),
and let the three other Galerkin solutions be
 \(u_{h_0},u_{h_1},u_{h_2}\) (with $h_0 > h_1 >h_2$). 
Because all the meshes are obtained by refinement from the same mesh, we can interpolate a finite element function from one mesh to the other, and we let \(I_{\text{ref}}\) denote interpolation to the reference mesh.
The relative \(L^2\) Galerkin errors are defined as \(\lVert I_{\text{ref}}(u_{h_i})-u_{\text{ref}}\rVert_{L^2(\Gamma_{\text{ref}})}/\lVert u_{\text{ref}}\rVert_{L^2(\Gamma_{\text{ref}})}\) for \(i=0,1,2\).

\begin{figure}[h!]
    \centering
    \includegraphics[width=0.8\textwidth]{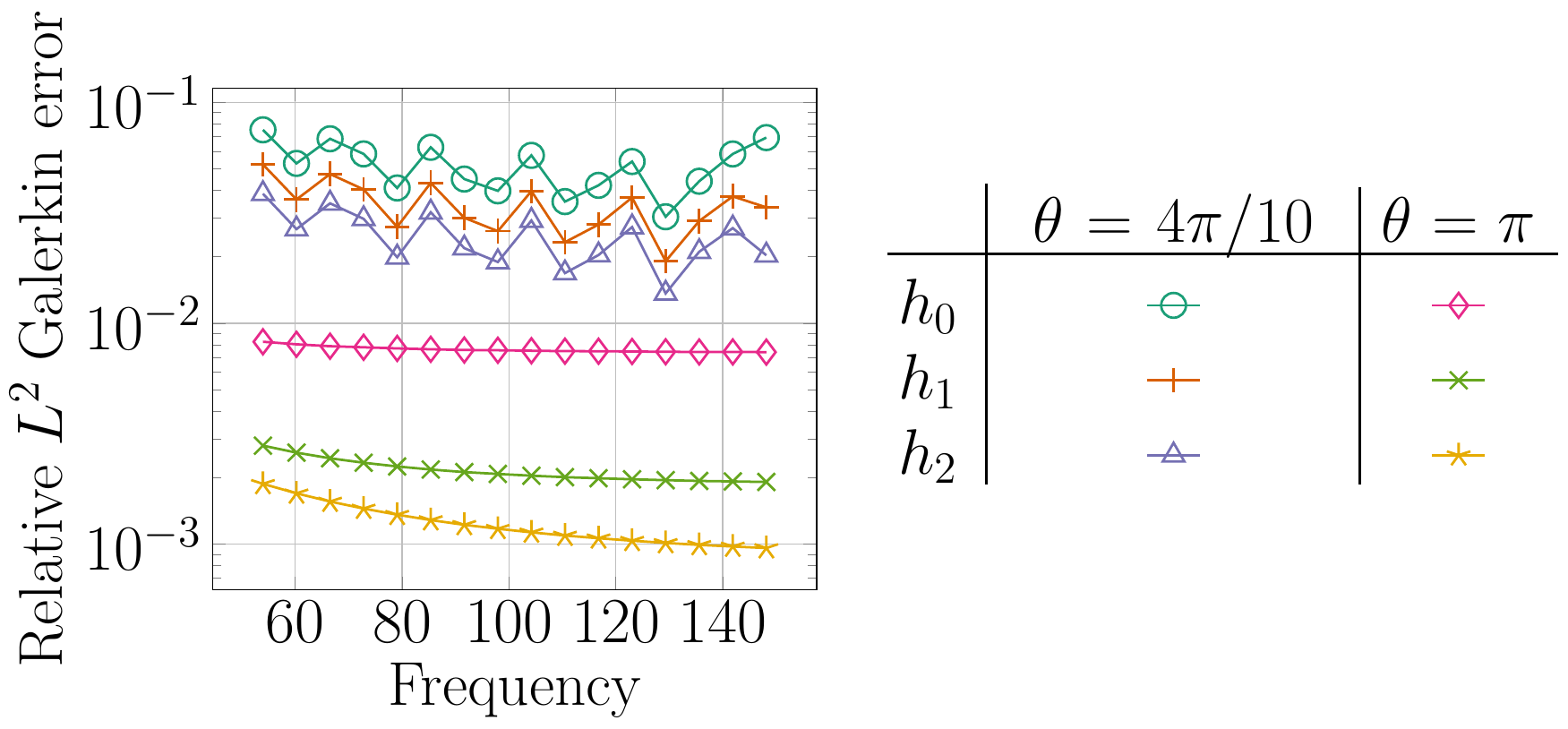}
    \caption{Relative \(L^2\) Galerkin error for Galerkin solutions on three nested meshes for $k= k_{m,0}^e$ and for two different plane-wave directions.}\label{fig:error_Galerkin}
\end{figure}

Figure \ref{fig:error_Galerkin} plots these errors for $k= k_{m,0}^e$ for $k \in (50,150)$
 and for two different choices of the plane-wave direction: 
 $a=(\cos\theta, \sin\theta)$ with 
$\theta =4\pi/10$ (the plane wave is almost vertical and enters the cavity) and $\theta=\pi$ (the plane wave is horizontal, travelling from the right, and thus does not enter the cavity). In this figure, we see exactly the feature F1.
 
\begin{figure}[h!]
    \centering
    \includegraphics[width=0.4\textwidth]{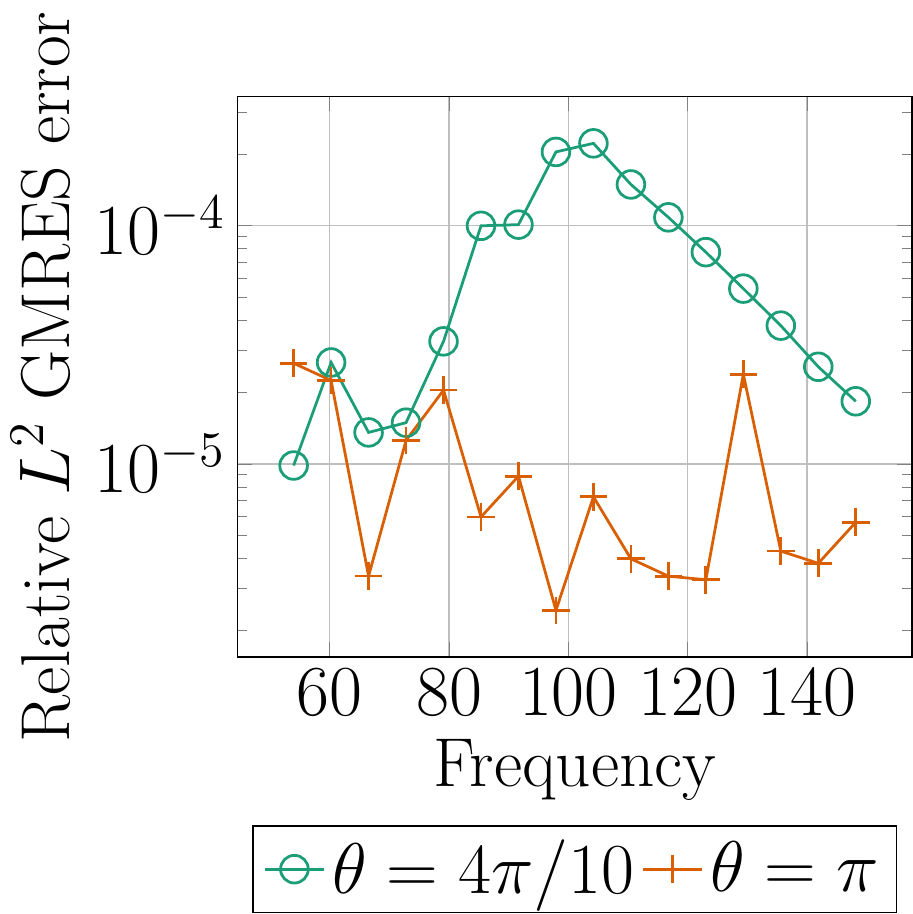}
    \caption{The relative \(L^2\) error between the GMRES solution of the Galerkin equations and the solution via a direct solver 
     for $k= k_{m,0}^e$ and for two different plane-wave directions.}\label{fig:error_GMRES}
\end{figure}

\paragraph{Numerical experiments demonstrating F4.} 
These experiments work only on the mesh with \(10\) points by wavelength used in the experiments for F1. 
We compute the error between the solutions of the Galerkin equations computed by (i) SuperLU~\cite{lidemmel03}, denoted by $u_{h_0}$ as above, and (ii) GMRES with a relative tolerance of \(1\times 10^{-6}\). 
We then normalise this difference by  \(\lVert u_{h_0}\rVert_{L^2(\Gamma)}\).
This relative $L^2$ GMRES error is plotted in Figure \ref{fig:error_GMRES}, both for $\theta=4\pi/10$ and $\pi$, and we see the dependence on angle as stated in F4.

\section{Discussion of the results of \cite{Titley2014} on how the GMRES residual depends on the right-hand-side vector (relevant for F3(c))}
\label{app:GMRES_rhs}

We consider solving the linear system  $\genmatrix \bfx = \bfb$ with GMRES, as described in \S\ref{sec:GMRES}.
We assume $\genmatrix$ is diagonalisable, so that
\beq
\label{eq:diag}
\genmatrix= \bfV {\bf \Lambda} \bfV^{-1},
\eeq
where the columns of $\bfV$, denoted by $(\bfv_1, \bfv_2, \dots, \bfv_n)$,  are the right eigenvectors of $\genmatrix$ corresponding to the eigenvalues $\lambda_1, \lambda_2, \dots, \lambda_n$, respectively. 

We expand $\bfb$ as a linear combination of right eigenvectors with corresponding components $\beta_j$, i.e.,
\beqs
\bfb = \sum^n_{j=1} \beta_j\bfv_j = \bfV \bfbeta,
\eeqs
and let $\bfbeta':= \bfbeta/\|\bfb\|_2$.
The result \cite[Theorem 2.2]{Titley2014} states that 
\beq
\label{eq:weightedbd}
\frac{\|\bfr_m\|_2}{\|\bfr_0\|_2} \le \N{{ \bf V}}_2 \min_{p \in \bbP_m, \; p(0)=1} \left ( \sum_{j=1}^{n} |\beta'_j|^2\;|p(\lambda_j)|^2 \right )^{1/2};
\eeq
i.e.,  the relative residual is bounded above by $\|\bfV\|_2$ times the residual of a polynomial least-squares approximation problem on the spectrum of $\genmatrix$, weighted by the scaled components $\beta'_j$.

We now consider the special case that the spectrum of $\genmatrix$ has a single outlier near zero, plus a cluster bounded well-away from zero; i.e.  
$| \lambda_1| \approx 0$ and, for $j=2,\ldots,n$, $\lambda_j$ is such that $|\lambda_j-c| < \rho$, with $\rho\ll |c|$.
In this case, the bound (\ref{eq:weightedbd}) becomes
\begin{align*}
\frac{\|r_m\|_2}{\|r_0\|_2} &\leq \|V\|_2 \min_{p \in \bbP_m, p(0)=1} \left ( |\beta'_1|^2|p(\lambda_1)|^2 + \sum_{j=2}^{n} |\beta'_j|^2\;|p(\lambda_j)|^2 \right )^{1/2},\\
&\approx \|V\|_2 \min_{p \in \bbP_m, p(0)=1} \left ( |\beta'_1|^2|p(\lambda_1)|^2 + \sum_{j=2}^{n} |\beta'_j|^2\;|p(c)|^2 \right )^{1/2}.
\end{align*}

In its initial stages, GMRES tries to construct a polynomial that is one at zero and is very small at $\lambda_1$ before dealing with the values $p(\lambda_j)$ for the eigenvalues in the cluster. This feature of GMRES has often been remarked on, see for example, \cite[Discussion after Theorem 1]{Meurant2015} and \cite[\S4 and Figure 1]{ElSiWa:02}. Now, if $\bfb$ is varied such that the relative coefficient $\beta'_1$ increases, then the weight on $p(\lambda_1)$ increases, making it harder for GMRES to make the term $|\beta'_1|^2|p(\lambda_1)|^2$ very small. Therefore, in this special scenario of a single outlier near zero, one would expect the number of GMRES iterations to depend significantly  on the size of $\beta'_1$, with the number of iterations increasing as $|\beta'_1|$ increases (and decreasing if $|\beta'_1|$ decreases).

Though the argument above is heuristic and applies to a very simple situation, 
it illustrates that the size of the component of the right-hand-side vector in the direction of an eigenvector with corresponding eigenvalue very close to zero is likely to significantly influence the bound \eqref{eq:weightedbd} on the residual in GMRES.

\section{Eigenvalues and eigenfunctions of the Laplacian in an ellipse in terms of Mathieu functions}\label{sec:Mathieu}

The  eigenvalue problem for the Dirichlet/Neumann Laplacian in the ellipse $E$ \eqref{eq:ellipse} is 
\begin{align}\label{eq:Helmholtz_ellipse}
        \Delta u + k^2 u = 0 \text{ in } E, \qquad \text{ either }
        u = 0 \text{ or } \partial_n u=0 \text{ on } \partial E.
\end{align}
It is customary to call $\{ (x_1,0) : |x_1|\leq a_1\}$ the \emph{major axis}, $\{ (0,x_2) : |x_2|\leq a_2\}$ the \emph{minor axis}, \(\epsilon=\sqrt{1-\frac{a_2^2}{a_1^2}}\)  the \emph{eccentricity},  and \(a:=\sqrt{a_1^2-a_2^2}\)
the \emph{linear eccentricity}. 

We use the following change of variables, introduced in~\cite{Ma:68}, 
\begin{align}\label{eq:change_variables}
    \left\{
    \begin{aligned}
        x_1 &= a \cosh (\mu) \cos (\nu), \\
        x_2 &= a \sinh (\mu) \sin (\nu), \\
    \end{aligned}    
    \right.
\end{align}
so that 
\beqs
E=\Big\{(a \cosh (\mu) \cos (\nu), a \sinh (\mu) \sin (\nu))\in \bbR^2 \, \,:\, \, 0 \leq \mu \leq \mu_0\text{,}\, 0\leq \nu <2 \pi \Big\},
\eeqs
where \(\mu_0:= \cosh^{-1} (a_1/a) = \sinh^{-1}(a_2/a)\). 
(Note that we have used the same notation as in \cite[Appendix A]{BeChGrLaLi:11} for variable names etc.)

Substituting \(u(x_1,x_2)=M(\mu)N(\nu)\) into \eqref{eq:Helmholtz_ellipse}, we find
\begin{empheq}[left=\empheqlbrace]{align}
    N''(\nu)+(\alpha-2 q \cos (2\nu))N(\nu)&=0 &\text{(standard Mathieu equation)}\label{eq:standard_Mathieu_equations}\\
    M''(\mu)-(\alpha-2 q \cosh (2\mu))M(\mu)&=0 &\text{(modified Mathieu equation)}\label{eq:modified_Mathieu_equations}
\end{empheq}
where \(\alpha\) is the separation constant and 
\begin{align}\label{eq:q}
    q=\frac{(k a)^2}{4}.
\end{align}
Since \eqref{eq:standard_Mathieu_equations} is symmetric in $\nu$, if \(N(\nu)\) is solution of \eqref{eq:standard_Mathieu_equations}, then so is \(N(-\nu)\); we therefore restrict attention to 
solutions of \eqref{eq:standard_Mathieu_equations} that are even or odd.

We therefore seek solutions of \eqref{eq:standard_Mathieu_equations} and \eqref{eq:modified_Mathieu_equations}, with $N$ even or odd, satisfying 
    \begin{equation}\label{eq:periodic_boundary_conditions}
        N(0)=N(2\pi) \quad \text{and} \quad N'(0)=N'(2\pi),
    \end{equation}
    to ensure periodicity in $\nu$, and 
    \begin{equation}\label{eq:dirichlet_condition}
\text{ either }      \quad  M(\mu_0)=0 \quad\text{ or }\quad M'(\mu_0)=0,
    \end{equation}
    to ensure the zero Dirichlet/Neumann boundary condition on $\partial E$.
Furthermore, to obtain a well-defined solution at \(\mu=0\),~\cite{Mc:51} shows that 
we also need $M$ to satisfy
    \begin{subequations}\label{eq:zero_condition}
    \begin{empheq}[left=\empheqlbrace]{align}
        M'(0)&=0 \quad\text{ if }N\text{ is even,}\label{eq:zero_condition_even}\\
        M(0)&=0 \quad\text{ if }N\text{ is odd.}\label{eq:zero_condition_odd}
    \end{empheq}
\end{subequations}

In analogy with polar coordinates, the boundary value problem for $N$ \eqref{eq:standard_Mathieu_equations}, \eqref{eq:periodic_boundary_conditions} is called the \emph{angular problem}, while the boundary value problem for $M$ defined by \eqref{eq:modified_Mathieu_equations}, \eqref{eq:dirichlet_condition}, and \eqref{eq:zero_condition} is called the \emph{radial problem}. An eigenmode \(u(x_1,x_2)=M(\mu)N(\nu)\) with  \(N\) is even is called an \emph{even eigenmode}, 
and one with $N$ odd is called an \emph{odd eigenmode}.

The multiparametric spectral problems defined by \eqref{eq:standard_Mathieu_equations}, \eqref{eq:periodic_boundary_conditions}, \eqref{eq:modified_Mathieu_equations}, and \eqref{eq:zero_condition_even} for even modes, and \eqref{eq:standard_Mathieu_equations}, \eqref{eq:periodic_boundary_conditions}, \eqref{eq:modified_Mathieu_equations}, and \eqref{eq:zero_condition_odd} for odd modes are well-defined by \cite{Ne:10} \footnote{Note that \cite{Ne:10} only considers the Dirichlet problem.}. Indeed,
\begin{itemize}
    \item for \((m,n)\in \{0,1,2,\ldots\}^2\), there exists a unique pair \((\alpha_{m,n}^e,q^e_{m,n})\in \bbR \times (0,\infty)\) such that \eqref{eq:standard_Mathieu_equations}, \eqref{eq:periodic_boundary_conditions}, \eqref{eq:modified_Mathieu_equations}, and \eqref{eq:zero_condition_even} have non trivial solutions \(M_{m,n}^e(\mu)\) and \(N_{m,n}^e(\nu)\) with respectively \(m\) zeros in \((0,\mu_0)\) and \(n\) zeros in \([0,\pi)\),
    \item for \((m,n)\in \{0,1,2,\ldots\}\times \{1,2,\ldots\}\), there exists a unique pair \((\alpha_{m,n}^o,q^o_{m,n})\in \bbR \times (0,\infty)\) such that \eqref{eq:standard_Mathieu_equations}, \eqref{eq:periodic_boundary_conditions}, \eqref{eq:modified_Mathieu_equations}, and \eqref{eq:zero_condition_odd} have a non trivial solutions \(M_{m,n}^o(\mu)\) and \(N_{m,n}^o(\nu)\) with respectively \(m\) zeros in \((0,\mu_0)\) and \(n\) zeros in \([0,\pi)\),
\end{itemize}

Recall that $q$ and $k$ are related by \eqref{eq:q}.
The frequencies associated with \(q_{m,n}^e\) and \(q_{m,n}^o\) 
are denoted by \(k_{m,n}^e\) and \(k_{m,n}^o\), 
respectively, and the associated eigenfunctions are denoted by \(u_{m,n}^e\) and \(u_{m,n}^o\).
By~\cite[Equations (21) and (22)]{Ne:10}, \(k^e_{m,n}\) and \(k^o_{m,n}\) both increase with \(m\) and \(n\), which is consistent with the fact the only accumulation point in the spectrum of the Laplacian is infinity.

\paragraph{References for the proof of Theorem \ref{thm:ellipse}.}

The result follows by combining the following three ingredients.

(i) The results of \cite[Equation A.16]{BeChGrLaLi:11} and \cite[Theorem 3.1]{NgGr:13}  that the eigenfunctions associated with $k_{m,n}^{e/0}$ exponentially localise about the minor axis 
as $m\tendi$ for fixed $n$.

(ii) The arguments in \cite[Proof of Theorem 2.8]{BeChGrLaLi:11} that construct quasimodes of the exterior Dirichlet problem from Dirichlet eigenfunctions of the ellipse that are exponentially localised (note that these arguments also apply to the Neumann problem). \footnote{In \cite{BeChGrLaLi:11}, $\Omega_+$ is assumed to contain the whole ellipse $E$. However, inspecting the proof, we see that the result remains unchanged if $E$ is replaced with the convex hull of the neighbourhoods of $(0,\pm a_2)$.}

Regarding (i): 
for \(\nu_0\in (0,\pi/2)\), let 
\begin{align*}
&    E_{\nu_0}:=\Big\{(a \cosh (\mu) \cos (\nu), a \sinh (\mu) \sin (\nu))\in \bbR^2 \, \,:\,\, 0 \leq \mu \leq \mu_0,
    \\
 &   \hspace{5cm}
    \text{ and either} \quad 0\leq \lvert\nu\rvert <\nu_0 \text{ or } \lvert\pi - \nu\rvert <\nu_0 \Big\};
\end{align*}
and 
\begin{align*}
    \rho_{\nu_0}^{e/o}(m,n):= \left(\dfrac{\int_{E_{\nu_0}} (u_{m,n}^{e/o})^2}{\int_E (u_{m,n}^{e/o})^2}\right)^{1/2}
\end{align*}    
i.e., $E_{\nu_0}$ corresponds to the ``wings'' of the ellipse, away from the minor axis, and $\rho_{\nu_0}^{e/o}(m,n)$ measures the mass of $u^{e/o}_{m,n}$ in these regions.
By \cite[Equation (A.16)]{BeChGrLaLi:11}, there exists \(K^e(\nu_0)>0\) such that \(\rho^e_{\nu_0}(m,0) \lesssim e^{-k_{m,0}^e}\) for any \(k^e_{m,0}>K^e(\nu_0)\).
By \cite[Theorem 3.1]{NgGr:13}, for \(n\) fixed, there exists \(K^{e/o}_n(\nu_0)>0\) such that, if \(k_{m,n}^{e/o}>K^{e/o}_n(\nu_0)\), then \(\rho^{e/o}_{\nu_0}(m,n) \lesssim e^{-k_{m,n}^{e/o}}\). Note that, although the inequality \cite[Equation 3.7]{NgGr:13} in \cite[Theorem 3.1]{NgGr:13} is stated for all $m$ and $n$ sufficiently large, the factor $D_n$ on the right-hand side of 
 \cite[Equation 3.7]{NgGr:13} 
blows up if $n\tendi$, and thus exponential localisation is proved in \cite{NgGr:13} for fixed $n$ as $m\tendi$.

\section{Calculating the constant $V_{\rm loc}$ in the Weyl asymptotics \eqref{eq:Nloc} for $N_{\rm loc}$.}\label{sec:Weyl}

Recall from \S\ref{sec:smalllarge} that we need to compute the volume, $V_{\rm loc}$, in phase space of the integrable tori contained entirely inside the small and large cavities and show that \eqref{eq:Vloc} holds.

Let 
$$
p((x,y),(\xi,\eta))=\sqrt{\xi^2+\eta^2}.
$$
As in \S\ref{sec:Mathieu}, let $a:=\sqrt{a_1^2-a_2^2}$. We change variables $\theta\in[0,2\pi]$ and $\omega \in (0,\cosh^{-1}(a_1/a))$
$$
x= a\cosh( \omega)\cos (\theta) ,\qquad y=a\sinh(\omega)\sin(\theta)
$$
(this is the same change of variables as \eqref{eq:change_variables} but with different variable names).

We now make a symplectic change of variables following \cite[\S2.3, Example 3 and Theorem 2.6]{Zw:12}. If
$$
\gamma(\omega,\theta)=\big(a\cosh(\omega)\cos \theta, a\sinh(\omega)\sin(\theta)\big),
$$
then
$$
\partial\gamma(\omega,\theta)=a\begin{pmatrix} \sinh(\omega)\cos \theta&\cosh(\omega)\sin(\theta)\\
-\cosh(\omega)\sin(\theta)&\sinh(\omega)\cos(\theta)\end{pmatrix},
$$
and hence
$$
\partial\gamma(\omega,\theta)^{-1}=\frac{1}{a\big(\sinh^2(\omega)+\sin^2(\theta)\big)}\begin{pmatrix} \sinh(\omega)\cos(\theta)&-\cosh(\omega)\sin(\theta)\\
\cosh(\omega)\sin(\theta)&\sinh(\omega)\cos \theta\end{pmatrix}
$$
$$
(\partial\gamma(\omega,\theta)^{-1})^t=\frac{1}{a\big(\sinh^2(\omega)+\sin^2(\theta)\big)}\begin{pmatrix} \sinh(\omega)\cos(\theta)&\cosh(\omega)\sin(\theta)\\-\cosh(\omega)\sin(\theta)
&\sinh(\omega)\cos \theta\end{pmatrix},
$$
Therefore by \cite[Theorem 2.6]{Zw:12}, the corresponding symplectomorphism is given by $\kappa((\omega,\theta),(\omega^*,\theta^*))=((x,y),(\xi,\eta))$ with
\begin{align}\label{eq:xy}
x&=a\cosh(\omega)\cos(\theta), \qquad y=a\sinh(\omega)\sin(\theta),\\ \nonumber
\xi&=\frac{1}{a\big(\sinh^2(\omega)+\sin^2(\theta)\big)}\big(\sinh(\omega)\cos(\theta)\omega^*+\cosh(\omega)\sin(\theta)\theta^*\big), \quad\tand\\ \nonumber
\eta&=\frac{1}{a\big(\sinh^2(\omega)+\sin^2(\theta)\big)}\big(-\cosh(\omega)\sin(\theta)\omega^*+\sinh(\omega)\cos(\theta)\theta^*\big).
\end{align}
Then
$$
\widetilde{p}:=p\circ \kappa = \frac{\sqrt{(\omega^*)^2+(\theta^*)^2}}{a\sqrt{\sinh^2(\omega)+\sin^2(\theta)}},
$$
so that, by \cite[Theorem 2.10]{Zw:12},
\begin{align}\label{eq:dotomega}
&\dot{\omega}= \partial_{\omega^*}\widetilde{p}= \frac{\omega^*}{a\sqrt{(\omega^*)^2+(\theta^*)^2}\sqrt{\sinh^2(\omega)+\sin^2(\theta)}},\quad\tand
\\
&\dot{\theta}=\partial_{\theta^*}\widetilde{p}= \frac{\theta^*}{a\sqrt{(\omega^*)^2+(\theta^*)^2}\sqrt{\sinh^2(\omega)+\sin^2(\theta)}}.\label{eq:dottheta}
\end{align}
One can then easily check that 
$$
L:=\frac{-\sin^2(\theta)(\omega^*)^2+\sinh^2(\omega)(\theta^*)^2}{(\theta^*)^2+(\omega^*)^2}
$$
is conserved by the Hamiltonian flow of $\widetilde{p}$; i.e., $\dot{L}=0$.

Now, \eqref{eq:dotomega} and \eqref{eq:dottheta} imply that, on a billiard trajectory, $\dot{\theta}=0$ if and only if $\theta^*=0$, and $\dot{\omega}=0$ if and only if $\omega^*=0$. 
That is, a trajectory is tangent to a curve $\{\omega=\omega_0\}$ when $\omega^*=0$, and to $\{\theta=\theta_0\}$ when $\theta^*=0$. Note that  
$$
\theta^*=0 \quad \Rightarrow \quad \sin^2(\theta)=-L,\qquad\tand\qquad \omega^*=0\quad \Rightarrow \quad \sinh^2(\omega)=L.
$$
Next, observe that curves of constant $\theta$ are confocal hyperbole, and curves of constant $\omega$ are confocal ellipses. 
Since every trajectory is tangent to (possibly degenerate) confocal conic and $L$ is constant, a trajectory is tangent to a (possibly degenerate) confocal ellipse if and only if $L\geq 0$ and to a confocal hyperbola if and only if $L<0$. Since we are interested in the volume of phase space occupied by trajectories trapped near the minor axis, we consider the case $L<0$. In that case, on $\widetilde{p}=1$, the confocal hyperbola is given by
\beq\label{eq:hyperbola}
\frac{x^2}{a^2(1+L)}+\frac{y^2}{a^2L}=1,
\eeq
where we have used the fact that $\sin^2(\theta)=-L$ and \eqref{eq:xy}.
We therefore want to find the volume of 
$$
A:=\big\{L<\alpha<0\big\}\cap \big\{\widetilde{p}\leq 1\big\}
$$
 in the $(\omega,\theta,\omega^*,\theta^*)$ variables, i.e., the volume of the set 
$$
\frac{-\sin^2(\theta)(\omega^*)^2+\sinh^2(\omega)(\theta^*)^2}{(\omega^*)^2+(\theta^*)^2}<\alpha,\qquad  \frac{\sqrt{(\omega^*)^2+(\theta^*)^2}}{a\sqrt{\sinh^2(\omega)+\sin^2(\theta)}}\leq 1.
$$
We change variables in $(\omega^*,\theta^*)$ by letting
$$
(\omega^*,\theta^*)=a\sqrt{\sinh^2(\omega)+\sin^2(\theta)}(r\cos\phi,r\sin\phi)
$$
with $\phi\in [0,2\pi]$, $r\in [0,\infty)$.
Then $\widetilde{p}=r$, 
$$
\rd\omega^*\,\rd\theta^*\,\rd\omega \,\rd\theta =a^2(\sinh^2(\omega)+\sin^2(\theta))r \,\rd r \,\rd\phi\, \rd \omega\, \rd\theta,
$$
and
$$
A:=\Big\{r\leq 1,\quad \big(-\sin^2(\theta)\cos^2\phi+\sinh^2(\omega)\sin^2\phi\big)\leq \alpha\Big\}.
$$
We now observe that
$$(-\sin^2(\theta)\cos^2\phi+\sinh^2(\omega)\sin^2\phi)\leq \alpha$$
if and only if
$$
(-\sin^2(\theta)+(\sinh^2(\omega)+\sin^2(\theta))\sin^2\phi)\leq \alpha
$$
if and only if
$$
\sin^2\phi\leq\frac{1}{(\sinh^2(\omega)+\sin^2(\theta))}\Big( \alpha+\sin^2\theta\Big).
$$
So
$$
A:=\big\{ r\leq 1\big\} \cap \big\{|\phi|\leq \Phi(\omega,\theta)\big\}\cup \big\{|\phi-\pi|\leq \Phi(\alpha,\omega, \theta)\big\},
$$
where $\Phi$ is defined by \eqref{eq:Phi}.
Therefore, the volume of $A$ equals
\begin{align*}
&\int_A a^2(\sinh^2(\omega)+\sin^2(\theta))r \,\rd r\,\rd\phi\, \rd\omega\, \rd\theta\\
&=\int_0^{2\pi}\int_0^{\cosh^{-1}(\frac{a_1}{a})}\int_0^1\left(\int_{-\Phi}^{\Phi}\rd\phi+\int_{\pi-\Phi}^{\pi+\Phi}\rd\phi\right) a^2\big(\sinh^2(\omega)+\sin^2(\theta)\big)r\,\rd r\, \rd\omega \,\rd\theta\\
&=\int_0^{2\pi}\int_0^{\cosh^{-1}(\frac{a_1}{a})}\int_0^14\Phi(\alpha,\omega,\theta) \, a^2(\sinh^2(\omega)+\sin^2(\theta))r\, \rd r\, \rd\omega\, \rd\theta\\
&=\int_0^{2\pi}\int_0^{\cosh^{-1}(\frac{a_1}{a})}2\Phi(\alpha,\omega,\theta) \,a^2(\sinh^2(\omega)+\sin^2(\theta)) \rd\omega\, \rd\theta,
\end{align*}
which is \eqref{eq:Vloc}.

Finally, to determine the relevant $\alpha$, we find the points where the boundary of the ellipse $E$ \eqref{eq:ellipse} meets the hyperbola \eqref{eq:hyperbola}, i.e., 
\beq\label{eq:twoineq}
1=\frac{x^2}{a_1^2}+\frac{y^2}{a_2^2}=\frac{x^2}{(a_1^2-a_2^2)(1+L)}+\frac{y^2}{(a_1^2-a_2^2)L}.
\eeq
Rearranging the second inequality in \eqref{eq:twoineq}, we find that
\beqs
\frac{y^2}{x^2}=\frac{a_2^2}{a_1^2}\left(\frac{-L}{1+L}\right)
\eeqs
and so using the first inequality in \eqref{eq:twoineq}, we find that
\beqs
\frac{y^2}{a_2^2}\left(\frac{1+L}{-L}\right)=\frac{x^2}{a_1^2}= 1- \frac{y^2}{a_2^2}, 
\eeqs
which implies that $L= - y^2 /a_2^2$.
Therefore, 
if the ellipse is cut at $y_{\rm cut}$, then $\alpha_{\rm cut}= -y_{\rm cut}^2/a_2^2$, i.e., \eqref{eq:alphacut} holds,
and the volume of the relevant piece of phase space is indeed given by \eqref{eq:Vloc}.

\section*{Acknowledgements}

EAS  gratefully acknowledges discussions with Alex Barnett (Flatiron Institute) that started his interest in eigenvalues of discretisations of the Helmholtz equation under strong trapping.
In addition, all the authors thank Barnett for giving them insightful comments on an earlier version of this paper. 
PM thanks Pierre Jolivet (Institut de Recherche en Informatique de Toulouse, CNRS) and Pierre-Henri Tournier (Sorbonne Universit\'e, CNRS) for their help with the software FreeFEM.
The authors thank the referees for their careful reading of the paper and numerous suggestions for improvement.
This research made use of the Balena High Performance Computing (HPC) Service at the University of Bath.
PM and EAS were supported by EPSRC grant EP/R005591/1.

\footnotesize{
\bibliographystyle{plain}
\bibliography{biblio_combined_sncwadditions.bib}
}

\end{document}